\input amstex
\input amsppt.sty
\input epsf.sty

\magnification 1200
\voffset -1cm

\TagsOnRight
\NoBlackBoxes
\NoRunningHeads

\define\GT{\Bbb G \Bbb T}

\define\R{\Bbb R}
\define\Z{\Bbb Z}
\define\C{\Bbb C}

\define\X{\frak X}
\define\al{\alpha}
\define\be{\beta}
\define\ga{\gamma}
\define\la{\lambda}
\define\Om{\Omega}
\define\om{\omega}
\define\x{\goth X}
\define\cl{\Cal L}
\define\inr{\operatorname{in}}
\define\out{\operatorname{out}}
\define\xin{\x_{\inr}}
\define\xout{\x_{\out}}
\define\Dim{\operatorname{Dim}}
\define\psiin{\psi_{\inr}}
\define\psiout{\psi_{\out}}
\define\Rout{R_{\out}}
\define\Sout{S_{\out}}
\define\Rin{R_{\inr}}
\define\Sin{S_{\inr}}
\define\wt{\widetilde}

\define\const{\operatorname{const}}
\define\p{\goth p}

\define\Ga{\Gamma}
\define\s{\Sigma}

\define\pp{\goth P}
\define\sgn{\operatorname{sgn}}
\define\Conf{\operatorname{Conf}}
\define\ep{\varepsilon}
\define\tht{\thetag}
\define\A{\Cal A}
\define\tP{\wt{\Cal P}}
\define\tK{\wt{K}^{(N)}}
\define\Kcirc{\overset\circ\to K}

\define\N{\Cal N}
\define\E{\Bbb E}
\define\de{\delta}
\define\cd{{\operatorname{CD}}}
\define\Conff{\Conf_{\operatorname{fin}}}

\define\r{\underset{\zeta=x}\to{\operatorname{Res}}\,}
\define\ru{\underset{u=x}\to{\operatorname{Res}}\,}
\define\Si{\Sigma}

\define\ze{\zeta}

\define\Treg{\Cal T_{\operatorname{reg}}}
\define\un{\underline}
\define\tr{\operatorname{tr}}
\define\ain{A_{\operatorname{in}}}
\define\aout{A_{\operatorname{out}}}
\define\n{{(N)}}
\define\Zp{\Cal Z_{\operatorname{princ}}}
\define\Zc{\Cal Z_{\operatorname{compl}}}
\define\Zd{\Cal Z_{\operatorname{degen}}}
\define\Zdm{\Cal Z_{\operatorname{degen},m}}
\define\Zdk{\Cal Z_{\operatorname{degen},k}}
\define\Zdl{\Cal Z_{\operatorname{degen},l}}
\define\zw{{z,z',w,w'}}
\define\Dadm{\Cal D_{\operatorname{adm}}}

\topmatter

\title
Harmonic analysis on the infinite--dimensional unitary group
and determinantal point processes
\endtitle
\author Alexei Borodin and Grigori Olshanski
\endauthor

\abstract
The infinite--dimensional unitary group $U(\infty)$ is the inductive
limit of growing compact unitary groups $U(N)$. In this paper we
solve a problem of harmonic analysis on $U(\infty)$ stated in
\cite{Ol3}. The problem consists in computing spectral
decomposition for a remarkable 4--parameter family of characters of
$U(\infty)$. These characters generate representations which should
be viewed as analogs of nonexisting regular representation of
$U(\infty)$.

The spectral decomposition of a character of $U(\infty)$ is described
by the spectral measure which lives on an infinite--dimensional space
$\Omega$ of indecomposable characters. The key idea which allows us
to solve the problem is to embed $\Omega$ into the space of point
configurations on the real line without 2 points. This turns the
spectral measure into a stochastic point process on the real line.
The main result of the paper is a complete description of the
processes corresponding to our concrete family of characters. We
prove that each of the processes is a determinantal point process.
That is, its correlation functions have determinantal form with a
certain kernel. Our kernels have a special `integrable' form and are
expressed through the Gauss hypergeometric function.

{}From the analytic point of view, the problem of computing the
correlation kernels can be reduced to a problem of evaluating uniform
asymptotics of certain discrete orthogonal polynomials studied
earlier by Richard Askey and Peter Lesky. One difficulty lies in the
fact that we need to compute the asymptotics in the oscillatory
regime with the period of oscillations tending to 0. We do this by
expressing the polynomials in terms of a solution of a discrete
Riemann--Hilbert problem and computing the (non--oscillatory)
asymptotics of this solution.

{}From the point of view of statistical physics, we study thermodynamic
limit of a discrete log--gas system. An interesting feature of this
log--gas is that its density function is asymptotically equal to the
characteristic function of an interval. Our point processes describe
how different the random particle configuration is from the typical
`densely packed' configuration.

In simpler situations of harmonic analysis on infinite symmetric
group and harmonic analysis of unitarily invariant measures on
infinite hermitian matrices similar results were obtained in our
papers \cite{BO1}, \cite{BO2}, \cite{BO4}.
\endabstract

\toc
\widestnumber\head{15}
\head {} Introduction \endhead
\head 1. Characters of the group $U(\infty)$
\endhead
\head 2. Approximation of spectral measures\endhead
\head 3. ZW--measures \endhead
\head 4. Two discrete point processes\endhead
\head 5. Determinantal point processes. General theory\endhead
\head 6. $\tP^{(N)}$ and $\Cal P^{(N)}$ as determinantal point
processes
\endhead
\head 7. The correlation kernel of the process $\tP^{(N)}$\endhead
\head 8. The correlation kernel of the process $\Cal P^{(N)}$\endhead
\head 9. The spectral measures and continuous point processes \endhead
\head 10. The correlation kernel of the process $\Cal P$\endhead
\head 11. Integral parameters $z$ and $w$\endhead
\head {} Appendix\endhead
\head {} Picture \endhead
\endtoc

\endtopmatter

\document

\head Introduction \endhead

\subhead (a) Preface \endsubhead
We tried to make this work accessible and interesting for a wide
category of readers. So we start with a brief explanation of the
concepts that enter in the title.

The purpose of harmonic analysis is to
decompose natural representations of a given group on irreducible
representations. By natural representations we mean those
representations that are produced, in a natural way, from the group
itself. For instance, this can be the regular representation, which
is realized in the $L^2$ space on the group, or a
quasiregular representation, which is built from the action of the
group on a homogeneous space.

In practice, a natural representation often comes together with a
distinguished cyclic vector. Then the decomposition into irreducibles
is governed by a measure, which may be called the {\it spectral
measure.\/} The spectral measure lives on the dual space to the
group, the points of the dual being the irreducible unitary
representations. There is a useful analogy in analysis: expanding a
given function on eigenfunctions of a self--adjoint operator. Here
the spectrum of the operator is a counterpart of the dual space.

If our distinguished
vector lies in the Hilbert space of the representation, then the
spectral measure has finite mass and can be made a probability measure.
\footnote{It may well happen that the distinguished vector belongs to
an extension of the Hilbert space (just as in analysis, one may well be
interested in expanding a function which is not square integrable).
For instance, in the case of the regular representation of a Lie
group one usually takes the delta function at the unity of the group,
which is not an element of $L^2$. In such a situation the spectral
measure is infinite. However, we shall deal with finite spectral
measures only.}

Now let us turn to {\it point processes\/} (or random point fields), which
form a special class of stochastic processes. In general, a
stochastic process is a discrete or continual family of random
variables, while a point process (or random point field) is a random
point configuration. By a (nonrandom) point configuration we mean an
unordered collection of points in a locally compact space $\frak X$.
This collection may be finite or countably infinite, but it cannot have
accumulation points in $\frak X$. To define a point process on $\frak
X$, we have to specify a probability measure on $\Conf(\frak X)$, the
set of all point configurations.

The classical example is the Poisson process, which is employed in
a lot of probabilistic models and constructions. Another important
example (or rather a class of examples) comes from random matrix
theory. Given a probability measure on a space of $N\times N$
matrices, we pass to the matrix eigenvalues and get in this way a
random $N$--point configuration. In a suitable scaling limit
transition (as $N\to\infty$), it turns into a point process living on
infinite point configurations.

As long as we are dealing with `conventional' groups (finite groups,
compact groups, real or $p$--adic reductive groups, etc.),
representation theory seems to have nothing in common with point
processes. However, the situation drastically changes when we turn to
`big' groups whose irreducible representations depend on infinitely
many parameters. Two basic examples are the infinite symmetric group
$S(\infty)$ and the infinite--dimensional unitary group $U(\infty)$,
which are defined as the unions of the ascending chains of finite or
compact groups
$$
S(1)\subset S(2)\subset S(3)\subset\dots, \qquad
U(1)\subset U(2)\subset U(3)\subset\dots,
$$
respectively. It turns out that for such groups, the clue to the problem
of harmonic analysis can be found in the theory of point processes.

The idea is to convert any infinite collection of parameters,
which corresponds to an irreducible representation, to a point
configuration. Then the spectral measure defines a point process, and
one may try to describe this process (hence the initial measure) using
appropriate probabilistic tools.

In \cite{B1}, \cite{B2}, \cite{BO1}, \cite{P.I--V} we applied this
approach to the group $S(\infty)$.
In the present paper we study the more complicated group $U(\infty)$.
Notice that the point processes arising from the spectral measures
do not resemble the Poisson process but are close to the processes of
the random matrix theory.

Now we proceed to a detailed description of the content of the paper.

\subhead (b) From harmonic analysis on $U(\infty)$ to a random
matrix type asymptotic problem \endsubhead
Here we summarize the necessary preliminary results established in
\cite{Ol3}. For a more detailed review see \S\S1--3 below.

The conventional definition of the regular representation is not
applicable to the group $U(\infty)$: one cannot define the $L^2$
space on this group, because $U(\infty)$ is not locally compact and
hence does not possess an invariant measure. To surpass this
difficulty we embed $U(\infty)$ into a larger space $\frak U$, which
can be defined as a {\it projective limit\/} of the spaces $U(N)$ as
$N\to\infty$. The space $\frak U$ is no longer a group but is still a
$U(\infty)\times U(\infty)$--space. That is, the two--sided action of
$U(\infty)$ on itself can be extended to an action on the space
$\frak U$. In contrast to $U(\infty)$, the space $\frak U$ possesses
a biinvariant finite measure, which should be viewed as a
substitute of the nonexisting Haar measure. Moreover, this
biinvariant measure is included into a whole family
$\{\mu^{(s)}\}_{s\in\C}$ of measures with good transformation
properties.
\footnote{The idea to enlarge an infinite--dimensional space in
order to build measures with good transformation properties is well
known. This is a standard device in measure theory on linear spaces,
but there are not so much works where it is applied to `curved'
spaces (see, however, \cite{Pi1}, \cite{Ner}). For the history of the
measures $\mu^{(s)}$ we refer to \cite{Ol3} and \cite{BO4}. A
parallel construction for the symmetric group case is given in
\cite{KOV}.} 
Using the measures $\mu^{(s)}$ we explicitly construct a family
$\{T_{zw}\}_{z,w\in\C}$ of representations, which seem to be a good
substitute of the nonexisting regular representation.
\footnote{More precisely, the $T_{zw}$'s are representations of the
group $U(\infty)\times U(\infty)$. Thus, they are a
substitute of the {\it biregular\/} representation. The reason why we
are dealing with the group $U(\infty)\times U(\infty)$ and not
$U(\infty)$ is explained in \cite{Ol1}, \cite{Ol2}. We also give in \cite{Ol3}
an alternative construction of the representations $T_{zw}$.}
In our understanding, the $T_{zw}$'s are `natural representations', and we
state the problem of harmonic analysis on $U(\infty)$ as follows:

\demo{Problem 1} Decompose the representations $T_{zw}$ on
irreducible representations.
\enddemo

This initial formulation then undergoes a few changes.

The first step follows a very general principle of representation
theory: reduce the spectral decomposition of representations to the
decomposition on extreme points in a convex set $\Cal X$
consisting of certain positive definite functions on the group.

In our concrete situation, the elements of the set $\Cal X$ are positive
definite functions on $U(\infty)$, constant on conjugacy classes
and taking the value 1 at the unity. These functions are called {\it
characters\/} of $U(\infty)$. The extreme points of $\Cal X$, or {\it
extreme characters,\/} are known. They are in a
one--to--one correspondence, $\chi^{(\om)}\leftrightarrow\om$, with the
points $\om$ of an infinite--dimensional region $\Om$ (the set $\Om$
and the extreme characters $\chi^{(\om)}$ are described in \S1 below).
An arbitrary character $\chi\in \Cal X$ can be written in the form
$$
\chi=\int_\Om \chi^{(\om)} P(d\om),
$$
where $P$ is a probability measure on $\Om$. The measure $P$ is defined
uniquely, it is called the {\it spectral measure\/} of the character
$\chi$.

Now let us return to the representations $T_{zw}$. We focus
on the case when the parameters $z,w$ satisfy the condition
$\Re(z+w)>-\tfrac12$. Under this restriction, our construction
provides a distinguished vector in $T_{zw}$. The matrix
coefficient corresponding to this vector can be viewed as a character
$\chi_{zw}$ of the group $U(\infty)$. The spectral
measure of $\chi_{zw}$ is also the spectral measure of
the representation $T_{zw}$ provided that $z$ and $w$ are not
integral.
\footnote{If $z$ or $w$ is integral then the distinguished vector is
not cyclic, and the spectral measure of$\chi_{zw}$ governs the
decomposition of a proper subrepresentation of $T_{zw}$.}

Furthermore, we remark that the explicit expression of
$\chi_{zw}$, viewed as a function in four parameters $z$,
$z'=\bar z$, $w$, $w'=\bar w$, correctly defines a character
$\chi_\zw$ for a wider set $\Dadm\subset\C^4$ of `admissible'
quadruples $(\zw)$. The set $\Dadm$ is defined by the inequality
$\Re(z+z'+w+w')>-1$ and some extra restrictions, see Definition 3.4
below. Actually, the `admissible' quadruples depend on 4 {\it real\/}
parameters.

This leads us to the following reformulation of Problem 1:

\demo{Problem 2} For any $(\zw)\in\Dadm$, compute the spectral
measure of the character $\chi_\zw$.
\enddemo

To proceed further we need to explain in what form we express the
characters. Rather than write them directly as functions on the group
$U(\infty)$ we prefer to work with their `Fourier coefficients'. Let
us explain what this means.

Recall that the irreducible representations of the compact group
$U(N)$ are labeled by the dominant highest weights, which are
nothing but $N$--tuples of nonincreasing integers
$\la=(\la_1\ge\dots\ge\la_N)$. For the reasons which are explained in
the text we denote the set of all these $\la$'s by $\GT_N$ (here `GT'
is the abbreviation of `Gelfand--Tsetlin'). For each $\la\in\GT_N$ we
denote by $\wt\chi^{\,\la}$ the normalized character of the
irreducible representation with highest weight $\la$. Here the term
`character' has the conventional meaning, and normalization means
division by the degree, so that $\wt\chi^{\,\la}(1)=1$. Given a
character $\chi\in \Cal X$, we restrict it to the subgroup $U(N)\subset
U(\infty)$. Then we get a positive definite function on $U(N)$, constant on
conjugacy classes and normalized at $1\in U(N)$. Hence it can be
expanded on the functions $\wt\chi^{\,\la}$, where the coefficients
(these are the `Fourier coefficients' in question) are nonnegative
numbers whose sum equals 1:
$$
\chi\mid_{U(N)}=\sum_{\la\in\GT_N}
P_N(\la)\wt\chi^{\,\la}; \qquad
P_N(\la)\ge0, \quad \sum_{\la\in\GT_N}P_N(\la)=1;
\quad N=1,2,\dots\,.
$$
Thus, $\chi$ produces, for any $N=1,2,\dots$, a probability measure
$P_N$ on the discrete set $\GT_N$. This fact plays an important role
in what follows.

For any character $\chi=\chi_\zw$ we dispose of an exact expression
for the `Fourier coefficients' $P_N(\la)=P_N(\la\mid\zw)$:
$$
\gathered
P_N(\la\mid \zw)=(\text{\rm normalization constant})\cdot
\prod_{1\le i<j\le N}(\la_i-\la_j-i+j)^2\\
\times \prod_{i=1}^N
\frac1{\Gamma(z-\la_i+i)\Gamma(z'-\la_i+i)
\Gamma(w+N+1+\la_i-i)\Gamma(w'+N+1+\la_i-i)}\,.
\endgathered  \tag0.1
$$
Hence we explicitly know the
corresponding measures $P_N=P_N(\,\cdot\,\mid\zw)$ on the sets
$\GT_N$. Formula \tht{0.1} is the starting point of the present paper.

In \cite{Ol3} we prove that for any character $\chi\in \Cal X$, its
spectral measure $P$ can be obtained as a limit of the measures $P_N$
as $N\to\infty$. More precisely, we define
embeddings $\GT_N\hookrightarrow\Omega$ and we show
that the pushforwards of the $P_N$'s weakly converge to $P$.
\footnote{The definition of the embeddings
$\GT_N\hookrightarrow\Omega$ is given in \S2(c) below.}

By virtue of this general result, Problem 2 is now
reduced to the following

\demo{Problem 3} For any `admissible' quadruple of parameters
$(\zw)$, compute the limit of the measures $P_N(\,\cdot\,\mid\zw)$,
given by formula \tht{0.1}, as  $N\to\infty$.
\enddemo

This is exactly the problem we are dealing with in the present paper.
There is a remarkable analogy between Problem 3 and asymptotic
problems of random matrix theory. We think this fact is important, so
that we discuss it below in detail. From now on the reader
may forget about the initial representation--theoretic motivation: we
switch to another language.

\subhead (c) Random matrix ensembles, log--gas systems, and determinantal
processes \endsubhead

\noindent Assume we are given a sequence of measures $\mu_1,\mu_2,\dots$ on
$\R$ and a parameter $\be>0$. For any $N=1,2,\dots$, we introduce a
probability distribution $P_N$ on the space of ordered $N$--tuples of real
numbers $\{x_1>\dots>x_N\}$ by
$$
P_N\left(\prod_{i=1}^N[x_i,x_i+dx_i]\right)
=(\text{\rm normalization constant})\cdot
\prod_{1\le i<j\le N}|x_i-x_j|^\be\cdot
\prod_{i=1}^N\mu_N(dx_i). \tag0.2
$$

Important examples of such distributions come from random matrix
ensembles $(E_N,\mu_N)$, where $E_N$ is a vector space of matrices
(say, of order $N$) and $\mu_N$ is a probability measure on $E_N$.
Then $x_1,\dots,x_N$ are interpreted as the eigenvalues of
an $N\times N$ matrix, and the distribution $P_N$ is induced by the
measure $\mu_N$. As for the parameter $\be$, it takes
values $1,2,4$, depending on the base field.

For instance, in the {\it Gaussian ensemble,\/} $E_N$ is the
space of real symmetric, complex Hermitian or quaternion Hermitian
matrices of order $N$, and $\mu_N$ is a Gaussian measure invariant
under the action of the compact group $O(N)$, $U(N)$ or $Sp(N)$,
respectively. Then $\be=1,2,4$, respectively.

If $\mu_N$ is absolutely continuous with respect to the Lebesgue
measure then the distribution \tht{0.2} is also absolutely
continuous, and its density can be written in the form
$$
F_N(x_1,\dots,x_N)=(\text{\rm constant})\cdot
\exp\left\{-\be\left(\sum_{1\le i<j\le N}\log|x_i-x_j|^{-1}+
\sum_{i=1}^N V_N(x_i)\right)\right\}. \tag0.3
$$
This can interpreted as the Gibbs measure of a system of $N$
repelling particles interacting through a logarithmic Coulomb potential
and confined by an external potential $V_N$. In mathematical physics
literature such a system is called a {\it log--gas system,\/}
see, e.g., \cite{Dy}.

Given a distribution of form \tht{0.2} or \tht{0.3}, one is interested
in the statistical properties of the random configuration $x=(x_i)$
as $N$ goes to infinity. A typical question concerns the asymptotic
behavior of the correlation functions. The {\it $n$--particle correlation
function\/}, $\rho_n^{(N)}(y_1,\dots,y_N)$, can be defined as the
density of the probability to find a `particle' of the random
configuration in each of $n$ infinitesimal intervals  $[y_i, y_i+dy_i]$.
\footnote{This is an intuitive definition only. In a rigorous
approach one defines the correlation {\it measures\/}, see, e.g.
\cite{Len, DVJ} and also the beginning of \S4 below.}

One can believe that under a suitable limit
transition the $N$--particle system `converges' to a point
process --- a probability distribution on infinite configurations of
particles. The limit distribution cannot be given by a formula of
type \tht{0.2} or \tht{0.3}. However, it can be characterized by its
correlation functions, which presumably are limits of the functions
$\rho_n^{(N)}$ as $N\to\infty$. The limit transition is usually
accompanied by a scaling (a change of variables depending on $N$),
and the final result may depend on the scaling. See, e.g.,
\cite{TW}.

The special case $\be=2$ offers much more possibilities for analysis
than the general one. This is due to the fact that for $\be=2$, the
correlation functions before the limit transition are readily
expressed through the orthogonal
polynomials $p_0, p_1, \dots$ with weight $\mu_N$. Namely, let
$S^{(N)}(y',y'')$ denote the $N$th Christoffel--Darboux kernel,
$$
\align
S^{(N)}(y',y'')&=\sum_{i=0}^{N-1}\frac{p_i(y')p_i(y'')}
{\Vert p_i\Vert^2}\\
&=(\text{\rm a constant})\cdot
\frac{p_N(y')p_{N-1}(y'')-p_{N-1}(y')p_N(y'')}{y'-y''}\,,
\qquad y',y''\in\R,
\endalign
$$
and assume (for the sake of simplicity only) that $\mu_N$ has a
density $f_N(x)$. Then the correlation functions are given by a
simple determinantal formula
$$
\rho_n^{(N)}(y_1,\dots,y_n)
=\det\left[S^{(N)}(y_i,y_j)
\sqrt{f_N(y_i)f_N(y_j)}\right]_{1\le i,j\le n}\,, \qquad n=1,2,\dots\,.
$$

If the kernel $S^{(N)}(y',y'')\sqrt{f_N(y')f_N(y'')}$ has a limit
$K(x',x'')$ under a scaling limit transition then the limit
correlation functions also have a determinantal form,
$$
\rho_n(x_1,\dots,x_n)=\det\left[K(x_i,x_j)\right]_{1\le i,j\le n}\,,
\qquad n=1,2\dots\,.  \tag0.4
$$
The limit kernel can be evaluated if one disposes of an appropriate
information about the asymptotic properties of the orthogonal
polynomials.

A point process whose correlation functions have the form \tht{0.4}
is called {\it determinantal,\/} and the corresponding kernel $K$ is
called the {\it correlation kernel.\/} Finite log--gas systems and their
scaling limits are examples of determinantal point processes. In
these examples, the correlation kernel is symmetric, but this
property is not necessary. Our study leads to processes with
nonsymmetric correlation kernels (see (k) below). A comprehensive
survey on determinantal point processes is given in \cite{So}.

\subhead (d) Lattice log--gas system defined by (0.1) \endsubhead
Remark that the expression \tht{0.1} can be transformed to the form
\tht{0.2}. Indeed, given $\la\in\GT_N$, set $l=\la+\rho$, where
$$
\rho=(\tfrac{N-1}2,\tfrac{N-3}2,\dots,-\tfrac{N-3}2,-\tfrac{N-1}2)
$$
is the half--sum of positive roots for $GL(N)$. That is,
$$
l_i=\la_i+\tfrac{N+1}2-i, \qquad i=1,\dots,N.
$$
Then $\Cal L=\{l_1,\dots,l_N\}$ is an $N$--tuple of distinct numbers
belonging to the lattice
$$
\frak X^{(N)}=\cases \Z, & \text{$N$ odd}, \\ \Z+\tfrac12, &
\text{$N$ even.}\endcases
$$
The measure \tht{0.1} on $\la$'s induces a probability measure on
$\Cal L$'s such that
$$
(\text{\rm Probability of $\Cal L$})=(\text{\rm a constant})\cdot
\prod_{1\le i<j\le N}(l_i-l_j)^2\cdot\prod_{i=1}^N f_N(l_i), \tag0.5
$$
where, for any $x\in\frak X^{(N)}$,
$$
f_N(x)=\frac 1{\Gamma\left(z-x+\frac{N+1}2\right)\Gamma\left(z'-
x+\frac{N+1}2\right)\Gamma
\left(w+x+\frac{N+1}2\right)\Gamma\left(w'+x+\frac{N+1}2\right)}\,.
\tag0.6
$$

Now we see that \tht{0.5} may be viewed as a {\it discrete log--gas
system\/} living on the lattice $\frak X^{(N)}$.

\subhead (e) Askey--Lesky orthogonal polynomials \endsubhead
The orthogonal polynomials that are defined by the weight function
\tht{0.6} on $\frak X^{(N)}$ are rather interesting. To our knowledge,
they appeared for the first time in Askey's paper \cite{As}. Then
they were examined in Lesky's papers \cite{Les1}, \cite{Les2}. We
propose to call them the {\it Askey--Lesky polynomials.\/} More
precisely, we reserve this term for the orthogonal polynomials that
are defined by a weight function on $\Z$ of the form
$$
\frac 1{\Gamma(A-x)\Gamma(B-x)\Gamma(C+x)\Gamma(D+x)}\,, \tag0.7
$$
where $A,B,C,D$ are any complex parameters such that \tht{0.7} is
nonnegative on $\Z$.

The Askey--Lesky polynomials are orthogonal polynomials of
hypergeometric type in the sense of \cite{NSU}. That is, they
are eigenfunctions of a difference analog of the hypergeometric
differential operator.

In contrast to classical orthogonal polynomials, the Askey--Lesky
polynomials form a {\it finite\/} system. This is caused by the fact
that (for nonintegral parameters $A,B,C,D$) the weight function has
slow decay as $x$ goes to $\pm\infty$, so that only finitely many
moments exist.

The Askey--Lesky polynomials admit an explicit expression in terms of
the generalized hypergeometric series ${}_3F_2(a,b,c;e,f;1)$ with
unit argument: the
parameters $A,B,C,D$ are inserted, in a certain way, in the indices
$a,b,c,e,f$ of the series. This allows us to explicitly express the
Christoffel--Darboux kernel in terms of the ${}_3F_2(1)$
series.

\subhead (f) The two--component gas system
\endsubhead
We have just explained how to reduce \tht{0.1} to a lattice log--gas system 
\tht{0.5},
for which we are able to evaluate the correlation functions. To solve
Problem 3, we must then pass to the large $N$ limit. However, the limit
transition that we need here is qualitatively different from typical scaling
limits of Random Matrix Theory. It can be shown that, as $N$ gets large,
almost all $N$ particles occupy positions inside $(-\frac N2, \frac N2)$.
Note that there are exactly $N$ lattice points in this interval, hence,
almost all of them are occupied by particles. More precisely, for any
$\varepsilon>0$, as $N\to\infty$,  the number of particles outside
$\left(-(\frac 12+\varepsilon)N,(\frac 12+\varepsilon)N\right)$ remains
finite almost surely.  In other words, this means that the density function
of our discrete log--gas is asymptotically equal to the characteristic
function of the $N$--point set of lattice points inside $(-\frac N2,\frac
N2)$.

At first glance, this picture looks discouraging. Indeed, we know that
in the limit all the particles are densely packed inside $(-\frac N2,\frac
N2)$, and there seem to exist no nontrivial limit point process. However,
the representation theoretic origin of the problem leads to the following
modification of the model which possesses a meaningful scaling limit.

Let us divide the lattice $\x^{(N)}$ into two parts, which will be
denoted by $\x^{(N)}_{\inr}$ and $\x^{(N)}_{\out}$:
$$
\gathered
\x^{(N)}_{\inr}=\left\{-\tfrac{N-1}2,-\tfrac{N-3}2,\dots,\tfrac{N-3}2,
\tfrac{N-1}2\right\},\\
\x^{(N)}_{\out}=\left\{\dots,-\tfrac{N+3}2,-\tfrac{N+1}2\right\}
\cup\left\{\tfrac{N+1}2,\tfrac{N+3}2,\dots\right\}.
\endgathered
$$
Here $\x^{(N)}_{\inr}$, the `inner' part, consists of $N$ points of the
lattice that lie on the interval $(-\tfrac N2,\tfrac N2)$, while
$\x^{(N)}_{\out}$, the `outer' part, is its complement in
$\x^{(N)}$, consisting of the points outside this interval. .

Given a configuration $\Cal L$ of $N$ particles sitting at points
$l_1,\dots,l_N$ of the lattice $\x^{(N)}$, we assign to it another
configuration, $X$, formed by the particles in
$\x^{(N)}_{\out}$ and the {\it holes\/} (i.e., the unoccupied
positions) in $\x^{(N)}_{\inr}$. Note that $X$ is a finite
configuration, too. Since the `interior' part consists of exactly $N$
points, we see that in $X$, there are equally many particles and
holes. However, there number is no longer fixed, it varies between 0
and $2N$, depending on the mutual location of $\Cal L$ and
$\x^{(N)}_{\inr}$. For instance, if these two sets coincide then $X$
is the empty configuration, and if they do not intersect then
$|X|=2N$.

Under the correspondence $\Cal L\mapsto X$ our random $N$--particle
system turns into a random system of particles and holes. Note that
$\Cal L\mapsto X$ is reversible, so that both systems are equivalent.

Rewriting \tht{0.5} in terms of the configurations $X$ one sees that
the new system can be viewed as a {\it discrete two--component log--gas
system\/} consisting of oppositely signed charges. Systems of such a
type were earlier investigated in the mathematical physics
literature (see \cite{AF},
\cite{CJ1}, \cite{CJ2}, \cite{G}, \cite{F1--3} and references therein). However, 
the
known concrete models are quite different from our system.

{}From what was said above it follows that all but finitely many particles of
the new system concentrate, for large $N$, near the points $\pm \frac N2$.
This suggests that if we shrink our phase space $\x^{(N)}$ by the factor of
$N$ (so that the points $\pm \frac N2$ turn into $\pm\frac 12$) then our
two--component log--gas should have a well--defined scaling limit. We prove
that such a limit exists and it constitutes a point process on $\Bbb
R\setminus\{\pm \frac 12\}$ which we will denote by $\Cal P$.

As a matter of fact, the process $\Cal P$ can be defined directly from the
spectral measure $P$ of the character $\chi_{z,z',w,w'}$ as we explain in
\S9. Moreover, knowing $\Cal P$ is almost equivalent to knowing $P$, see
the discussion before Proposition 9.7. Thus, we may restate Problem 3 as

\demo{Problem 4} Describe the point process $\Cal P$.
\enddemo

It turns out that the most convenient way to describe this point process is
to compute its correlation functions. Since the correlation functions
of $\Cal P$ define $\Cal P$ uniquely, we will be solving

\demo{Problem 4$\,'$} Find the correlation functions of $\Cal P$.
\enddemo

\subhead (g) Two correlation kernels of the two--component log--gas
\endsubhead 
There are two ways of computing the correlation functions of the
two--component log--gas system introduced above. The first one is based
on the {\it complementation principle}, see \cite{BOO, Appendix} and \S5(c) 
below, which says that
if we have a determinantal point process defined on a discrete set
$\frak Y=\frak Y_1\sqcup \frak Y_2$ then a new process whose point
configurations consist of particles in $\frak Y_1$ and holes in $\frak Y_2$,
is also determinantal. Furthermore, the correlation kernel of this new
process is easily expressed through the correlation kernel of the original
process. Thus, one way to obtain the correlation functions for the
two--component log--gas is to apply the complementation principle to the
(one--component) log--gas \tht{0.1}, whose correlation kernel is,
essentially, the Christoffel--Darboux kernel for Askey--Lesky orthogonal
polynomials.  Let us denote by $K_{compl}^{(N)}$ the correlation kernel for
the two--component log--gas obtained in this way.

Another way to compute the correlation functions of out two--component
log--gas is to notice that this system belongs to the class of point
processes with the following property:

\noindent The probability of a given point
configuration $X=\{x_1,\dots,x_n\}$ is given by
$$
\operatorname{Prob}\{X\}=\operatorname{const} \cdot
\det[L^{(N)}(x_i,x_j)]_{i,j=1}^n
$$
where $L^{(N)}$ is a $\x^{(N)}\times \x^{(N)}$ matrix, see \S6.
A simple general theorem shows that any point process with this property is
determinantal, and its correlation kernel $K^{(N)}$ is given by the relation
$K^{(N)}=L^{(N)}(1+L^{(N)})^{-1}$.

Thus, we end up with two correlation kernel $K_{compl}^{(N)}$ and $K^{(N)}$
of the same point process. Of course, these two kernels must not coincide. For
example, they may be related by conjugation:
$$
K_{compl}^{(N)}(x,y)=\frac{\phi(x)}{\phi(y)}\, K^{(N)}(x,y)
$$
where $\phi(\,\cdot\,)$ is an arbitrary nonvanishing function on $\x^{(N)}$.
(The determinants of the form $\det[K(x_i,x_j)]$ for two conjugate kernels
are always equal.) We show that this is indeed the case, and that the
function $\phi$ takes values $\pm 1$. Moreover, we prove this statement in
a more general setting of a two--component log--gas system obtained in a
similar way by particles--holes exchange from an {\it arbitrary} $\beta=2$
discrete log--gas system on the real line.

\subhead (h) Asymptotics \endsubhead
In our concrete situation the function $\phi$ is identically equal to 1 on
the set $\xout^{(N)}$ and is equal to $(-1)^{x-\frac{N-1}2}$ on the set
$\xin^{(N)}$. This means that if we want to compute the scaling limit of the
correlation functions of our two--component log--gas system as $N\to\infty$,
then only one of the kernels $K_{compl}^{(N)}$ and $K^{(N)}$ may be used for
this purpose, because the function $\phi$ does not have a scaling limit.
It is not hard to guess which kernel is `the right one' from the asymptotic
point of view.

Indeed, it is easy to verify that the kernel $L^{(N)}$
mentioned above has a well--defined scaling limit which we will denote by
$L$. It is a (smooth) kernel on $\Bbb R\setminus \{\pm \frac 12\}$. It is then
quite natural to assume that the kernel $K^{(N)}=L^{(N)}(1+L^{(N)})^{-1}$
also has a scaling limit $K$ such that $K=L(1+L)^{-1}$. Although this
argument is only partially correct (the kernel $L$ does not always define a
bounded operator in $L^2(\Bbb R)$), it provides good intuition. We prove that
for all admissible values of the parameters $z,z',w,w'$, the kernel $K^{(N)}$
has a scaling limit $K$, and this limit kernel is the correlation kernel
for the point process $\Cal P$.

Explicit computation of the kernel $K$ is our main result, and we state it in
\S10.

\subhead (i) Overcoming technical difficulties: Riemann--Hilbert approach
\endsubhead
The task of computing the scaling limit of $K^{(N)}$ as $N\to\infty$ is by no
means easy. As was explained above, this kernel coincides, up to a sign, with
$K_{compl}^{(N)}$ which, in turn, is easily expressible through the
Christoffel--Darboux kernel for the Askey--Lesky orthogonal polynomials.
Thus, Problem 4 (or 4$'$) may be restated as

\demo{Problem 5} Compute the asymptotics of the Askey--Lesky orthogonal
polynomials
\enddemo

Since it is known how to express these polynomials through the ${}_3F_2$
hypergeometric series, one might expect that the remaining part is rather
smooth and is similar to the situation arising in most $\beta=2$ random matrix
models.  That is, in the chosen scaling the polynomials converge with all the
derivatives to nice analytic functions (like sine or Airy) which then enter
in the formula for the limit kernel. As a matter of fact, this is indeed how
things look like on $\xout^{(N)}$. The limit kernel $K$ is not hard to
compute and it is expressed through the Gauss hypergeometric function
${}_2F_1$.

The problem becomes much more complicated when we look at
$\xin^{(N)}$. The basic reason is that this is the oscillatory zone for our
orthogonal polynomials, and in the scaling limit that we need the period of
oscillations tends to zero. Of course, one cannot expect to see any
uniform convergence in this situation.

Let us recall, however, that all we need is the asymptotics {\it on the
lattice}. This remark is crucial. The way we compute the asymptotics on the
lattice is, roughly speaking, as follows. We find meromorphic functions
with poles in $\xout^{(N)}$ which coincide, up to a sign, with our
orthogonal polynomials on $\xin^{(N)}$. These functions are also expressed
through the ${}_3F_2$ series and look more complicated than the polynomials
themselves. However, they possess a well--defined limit (convergence with all
the derivatives) which is again expressed through the Gauss hypergeometric
function. This completes the computation of the asymptotics.

The question is: how did we find these convenient meromorphic functions?
The answer lies in the definition of the kernel $K^{(N)}$ as
$L^{(N)}(1+L^{(N)})^{-1}$. It is not hard to see that the kernel $L^{(N)}$
belongs to the class of (discrete) integrable operators, see \cite{B3}. This
implies that the kernel $K^{(N)}$ can be expressed through a solution of a
(discrete) Riemann--Hilbert problem (RHP, for short), see \cite{B3, Proposition 
4.3}. It is
the solution of this Riemann--Hilbert problem that yields the needed
meromorphic functions.

The problem of finding this solution explicitly
requires additional efforts.  The key fact here is that the jump matrix of
this RHP can be reduced to a constant jump matrix by conjugation. It is a
very general idea of the inverse scattering method that in such a situation
the solution of the RHP must satisfy a difference (differential, in the case
of continuous RHP) equation. Finding this equation and solving it in
meromorphic functions yields the desired solution.

It is worth noting that even though the correct formula for the limit
correlation kernel $K$ can be guessed from just knowing the Askey--Lesky
orthogonal polynomials, the needed convergence of the kernels
$K^{(N)}$ was only possible to achieve through solving the RHP mentioned
above.

Let us also note that computing the limit of the solution of our RHP is not
completely trivial as well. The difficulty here lies in finding,
by making use of numerous known transformation formulas for the
${}_3F_2$ series, a presentation of the solution that would be
convenient for the limit transition .

\subhead (j) The main result \endsubhead
In (f) above we explained how to reduce our problem of harmonic analysis on 
$U(\infty)$ to the problem of computing the correlation functions of the
process $\Cal P$. In this paper we prove that the $n$th correlation function
$\rho_n(x_1,\dots,x_n)$ of $\Cal P$ has the determinantal form
$$
\rho_n(x_1,\dots,x_n)=\det[K(x_i,x_j)]_{i,j=1}^n,\qquad n=1,2,\dots
$$
Here $K(x,y)$ is a kernel on $\Bbb R\setminus \{\pm \frac 12\}$ which can be 
written in the form
$$
K(x,y)=\frac{F_1(x)G_1(y)+F_2(x)G_2(y)}{x-y}\,,\qquad x,y\in \R\setminus
\{\pm\tfrac 12\},
$$
where the functions $F_1,G_1,F_2,G_2$ can be expressed through the Gauss
hypergeometric function ${}_2F_1$. In particular, if $x>\frac 12$ and
$y>\frac 12$ we have
$$
\multline
F_1(x)=-G_2(x)=\frac{\sin(\pi z)\sin(\pi z')}{\pi^2}\\ \times
\left(x-\frac 12\right)^{-(\frac{z+z'}2+w')}\left(x+\frac 12\right)^{\frac{w'-
w}2}
{}_2F_1\left[\matrix
z+w',\,z'+w'\\ z+z'+w+w'\endmatrix\,\Biggl|\,\frac 1{\frac 12 -x}\right]\,,
\endmultline
$$
$$
\multline
G_1(x)=F_2(x)=\frac{\Ga(z+w+1)
\Ga(z+w'+1)\Ga(z'+w+1)\Ga(z'+w'+1)}{\Ga(z+z'+w+w'+1)\Ga(z+z'+w+w'+2)}
\\ \times \left(x-\frac 12\right)^{-(\frac{z+z'}2+w'+1)}\left(x+\frac 
12\right)^{\frac{w'-w}2}\,{}_2F_1\left[\matrix z+w'+1,\,z'+w'+1\\
z+z'+w+w'+2\endmatrix\,\Biggl|\,\frac 1{\frac 12 -x}\right]\,.
\endmultline
$$

A complete statement of the result can be found in Theorem 10.1 below.

\subhead (k) Symmetry of the kernel \endsubhead The correlation kernel $K(x,y)$ 
introduced above satisfies the following symmetry relations
$$
K(x,y)=\cases K(y,x)& \text{  if  } \left(|x|>\tfrac 12,\,
|y|>\tfrac 12\right)
\text{  or  } \left(|x|<\tfrac 12, |y|<\tfrac 12\right),\\
-K(y,x)&  \text{  if  } \left(|x|>\frac 12,\, |y|<\frac 12\right)
\text{  or  } \left(|x|<\tfrac 12, |y|>\tfrac 12\right).
\endcases
$$
Moreover, the kernel is real--valued. This implies that the restrictions
of $K$ to $(-\frac 12,\frac 12)\times (-\frac 12, \frac 12)$ and
$\left(\R\setminus[-\frac 12,\frac 12]\right)\times \left(\R\setminus[-\frac 12, 
\frac 12]\right)$ are Hermitian kernels, while the kernel $K$ on the whole line 
is a $J$--Hermitian\footnote{I.e., Hermitian with respect to the indefinite
inner product defined by the matrix $J=\bmatrix 1&0\\0&-1\endbmatrix$.}
kernel.

We have encountered certain $J$-Hermitian kernels in our work on harmonic
analysis on the infinite symmetric group, see \cite{BO1},
\cite{P.I--V}. At that time we were not aware of the fact
that examples of $J$-Hermitian correlation kernels had appeared
before in works of mathematical physicists on solvable models of
systems with positive and negative charged particles, see \cite{AF},
\cite{CJ1}, \cite{CJ2}, \cite{G}, \cite{F1--3} and references
therein.

As was explained in (f), our system also contains `particles of opposite 
charges'. The property of $J$--symmetry is closely related to this fact, see 
\S5(f),(g) for more details.

\subhead (l) Further development: Painlev\'e VI \endsubhead It is
well known that for a determinantal point process with a correlation
kernel $K$, the probability of having no particles in a region $I$ is
equal to the Fredholm determinant $\det(1-K_I)$, where $K_I$ is the
restriction of $K$ to $I\times I$. It often happens that such
{\it gap probability} can be expressed through a solution of a (second
order nonlinear ordinary differential) Painlev\'e equation. One of the main 
results of \cite{BD} is the following statement.

Let $K_s$ be the
restriction of the kernel $K(x,y)$ of (j) above to 
$(s,+\infty)\times(s,+\infty)$.
Set
$$
\gathered
\nu_1=\frac{z+z'+w+w'}2\,,\quad \nu_3=\frac{z-z'+w-w'}2\,,\quad
\nu_4=\frac{z-z'-w+w'}2\,,\\
\sigma(s)=\left(s^2-\tfrac 14\right)\frac
{d\ln\det(1-K_s)}{ds}-
\nu_1^2\,s+\frac {\nu_3\nu_4}{2}.
\endgathered
$$
Then $\sigma(s)$ satisfies the differential equation
$$
-\sigma'\left(\left(s^2-\tfrac 14\right)\sigma''\right)^2=\left(2\left(s\sigma'-
\sigma\right)\sigma'
-\nu_1^2\nu_3\nu_4\right)^2
-
(\sigma'+\nu_1^2)^2(\sigma'+\nu_3^2)(\sigma'+\nu_4^2).
$$
This differential equation is the so--called $\sigma$-form of the
Painlev\'e VI equation. We refer to \cite{BD, Introduction} for a brief
historical introduction and references on this subject. \cite{BD} also
contains proofs of several important properties of the kernel $K(x,y)$
which we list at the end of \S10 below.

\subhead (m) Connection with previous work \endsubhead
In \cite{BO1}, \cite{BO2}, \cite{B1}, \cite{B2}, \cite{BO4} we worked
out two other problems of harmonic analysis in the situations when
spectral measures live on infinite--dimensional spaces. We will
describe them in more detail and compare them to the problem of the
present paper. 

The problem of harmonic analysis on the group $S(\infty)$
was initially formulated in \cite{KOV}. It consists in decomposing
certain `natural' (generalized regular) unitary representations $T_z$
of the group $S(\infty)\times S(\infty)$, depending on a complex
parameter $z$. In \cite{KOV}, the problem was solved in the case when
the parameter $z$ takes integral values (then the spectral measure has a
finite--dimensional support). The general case presents more
difficulties and we studied it in a cycle of papers \cite{P.I--V},
\cite{BO1--3}, \cite{B1}, \cite{B2}. 
Our main result is that the spectral measure governing the
decomposition of $T_z$ can be described in terms of a determinantal
point process on the real line with one punctured point. The
correlation kernel was explicitly computed, it is expressed through a
confluent hypergeometric function (specifically, through the
W--Whittaker function).  

The second problem deals with decomposition of a family of unitarily
invariant probability measures on the space of all infinite Hermitian
matrices on ergodic components. The measures depend on one complex
parameter and essentially coincide with the measures $\{\mu^{(s)}\}$ 
mentioned in the beginning of (b) above. The problem of decomposition
on ergodic components can be also viewed as a problem of harmonic
analysis on an infinite--dimensional Cartan motion group. The main
result of \cite{BO4} states that the spectral measures in this case
can be interpreted as determinantal point processes on the real line
with a correlation kernel expressed through a confluent
hypergeometric function (this time this is the M--Whittaker
function). 

These two problems and the problem that we deal with in this paper
have a number of similarities. Already the descriptions of the spaces
of irreducible objects (see Thoma \cite{Th} for $S(\infty)$,
Pickrell \cite{Pi} and Olshanski--Vershik \cite{OV} for measures on
Hermitian matrices, and Voiculescu \cite{Vo} for $U(\infty)$) are
quite similar. Furthermore, all three models have some sort of an
approximation procedure using finite--dimensional objects, see
\cite{VK1}, \cite{OV}, \cite{VK2}, \cite{OkOl}. 
The form of the correlation kernels is also essentially the same, with
different special functions involved in different problems. 

It is worth noting that the similarity of theories for the two groups
$S(\infty)$ and $U(\infty)$ seems to be a striking phenomenon. In
addition to mentioned above, it can be traced in the geometric
construction of the `natural' representations and in probabilistic
properties of the corresponding point processes. At present we cannot
completely explain the nature of this parallelism (it looks quite
different from the well--known classical connection between the
representations of the groups $S(n)$ and $U(N)$).  

However, the differences between all these problems should not be 
underestimated. Indeed, the problem of harmonic analysis on
$S(\infty)$ is a problem of asymptotic combinatorics consisting in
controlling the asymptotics of certain explicit probability
distributions on partitions of $n$ as $n\to\infty$. One consequence
of such asymptotic analysis is a simple proof and generalization of
the Baik--Deift--Johansson theorem \cite{BDJ} on longest increasing
subsequences of large random permutations, see \cite{BOO} and
\cite{BO3}. The problem of decomposing measures on Hermitian matrices
on ergodic components is of a different nature. It belongs to Random
Matrix Theory which deals with asymptotics of probability
distributions on large matrices. In fact, for a specific 
value of the parameter, the result of \cite{BO4} reproduces one of the 
basic computations of Random Matrix Theory --- that of the scaling limit
of the Dyson's circular ensemble.

The problem solved in the present paper is more general comparing to
both problems described above. Our model here depends on a larger
number of parameters, it deals with a more complicated group and
representation structure, and the analysis requires substantial
amount of new ideas. Moreover, in appropriate limits this
model degenerates to both models studied earlier. The limits, of
course, are very different. On the level of correlation kernels this
leads to two different degenerations of the Gauss hypergeometric
function to confluent hypergeometric functions. We view the
$U(\infty)$--model as a unifying object for the combinatorial and
random matrix models, and we think that it sheds some light on the
nature of the recently discovered remarkable connections between
different models of these two kinds. 

The model of the present paper can be also viewed as the top of a
hierarchy of (discrete and continuous) probabilistic models leading
to determinantal point processes with `integrable' correlation
kernels. In the language of kernels this looks very much like the
hierarchy of the classical special functions. A description of the
`$S(\infty)$--part' of the hierarchy can be found in \cite{BO3}.
The subject of degenerating the $U(\infty)$--model to simpler models 
(in particular, to the two models discussed above) will be addressed
in a later publication.  

\subhead (n) Organization of the paper \endsubhead
In Section 1 we give a brief introduction to  representation theory and
harmonic analysis of the infinite--dimensional unitary group $U(\infty)$.
Section 2 explains how spectral decompositions of representations of $U(\infty)$ 
can be approximated by those for finite--dimensional groups
$U(N)$. In Section 3 we introduce a remarkable family of characters of 
$U(\infty)$ which we study in this paper. In Section 4 we reformulate the 
problem of harmonic analysis of these characters in the language of random point 
processes. Section 5 is the heart of the paper: there we develop general theory 
of discrete determinantal point processes
which will later enable us to compute the correlation functions of
our concrete processes. In Section 6 we 
show that the point processes introduced in Section 4 are determinantal.
In Section 7 we derive discrete orthogonal polynomials on $\Z$ with the weight
function \tht{0.7}. This allows us to write
out a correlation kernel for approximating point processes associated
with $U(N)$'s. Section 8 is essentially dedicated to representing
this correlation kernel in a form suitable for the limit transition
$N\to\infty$. The main tool here is discrete Riemann--Hilbert
problem. 
Section 9 establishes certain general facts about scaling limits of
point processes associated with restrictions of characters of $U(\infty)$
to $U(N)$. The main result here is that an appropriate scaling limit
yields the spectral measure for the initial character of $U(\infty)$.
In Section 10 we perform such a scaling limit for our concrete family
of characters. Section 11 describes a nice combinatorial degeneration
of 
our characters. In this degeneration the spectral measure loses its
infinite--dimensional support and turns into a Jacobi polynomial
ensemble. Finally, Appendix contains proofs of transformation formulas
for the hypergeometric series ${}_3F_2$ which were used in the computations.

\subhead (o) Acknowledgment \endsubhead
At different stages of our work we have been receiving a lot of
inspiration from conversations with Sergei Kerov and Yurii Neretin
who generously shared their ideas with us. We are extremely grateful
to them. 

We would also like to thank Peter Forrester and Tom Koornwinder for
referring us to the papers by R.~Askey and P.~Lesky, and to Peter
Lesky for sending us his recent preprint \cite{Les2}. 

This research was partially conducted during the period one of the
authors (A.B.) served as a Clay Mathematics Institute Long--Term
Prize Fellow.

\head 1. Characters of the group $U(\infty)$
\endhead

\subhead (a) Extreme characters \endsubhead
Let $U(N)$ be the group of unitary matrices of order $N$.
For any $N\ge2$ we identify $U(N-1)$ with the subgroup in $U(N)$
fixing the $N$th basis vector, and we set
$$
U(\infty)=\varinjlim U(N).
$$

One can view $U(\infty)$ as a group of matrices $U=[U_{ij}]_{i,j=1}^\infty$
such that there are finitely matrix elements $U_{ij}$ not equal to
$\delta_{ij}$, and $U^*=U^{-1}$.

A {\it character\/} of $U(\infty)$ is a
function $\chi:U(\infty)\to\C$ which is constant on conjugacy
classes, positive definite, and normalized at the unity
($\chi(e)=1$). We also assume that $\chi$ is continuous on each
subgroup $U(N)\subset U(\infty)$. The characters form a convex set.
The extreme points of this convex set are called the {\it extreme
characters.\/}

A fundamental result of the representation theory
of the group $U(\infty)$ is a complete description of extreme
characters. To state it we need some notation.

Let $\R^\infty$ denote the product of countably many
copies of $\R$, and set
$$
\R^{4\infty+2}=\R^\infty\times\R^\infty\times\R^\infty\times\R^\infty
\times\R\times\R.
$$
Let $\Om\subset\R^{4\infty+2}$ be the subset of sextuples
$$
\om=(\al^+,\be^+;\al^-,\be^-;\de^+,\de^-)
$$
such that
$$
\gather
\al^\pm=(\al_1^\pm\ge\al_2^\pm\ge\dots\ge 0)\in\R^\infty,\quad
\be^\pm=(\be_1^\pm\ge\be_2^\pm\ge\dots\ge 0)\in\R^\infty,\\
\sum_{i=1}^\infty(\al_i^\pm+\be_i^\pm)\le\de^\pm, \quad
\be_1^++\be_1^-\le 1.
\endgather
$$
Set
$$
\ga^\pm=\de^\pm-\sum_{i=1}^\infty(\al_i^\pm+\be_i^\pm)
$$
and note that $\ga^+,\ga^-$ are nonnegative.

To any $\om\in\Om$ we assign a function $\chi^{(\om)}$ on $U(\infty)$:
$$
\chi^{(\om)}(U)=
\prod_{u\in\operatorname{Spectrum}(U)}
\left\{e^{\ga^+(u-1)+\ga^-(u^{-1}-1)}
\prod_{i=1}^\infty\frac{1+\be_i^+(u-1)}{1-\al_i^+(u-1)}
\,\frac{1+\be_i^-(u^{-1}-1)}{1-\al_i^-(u^{-1}-1)}\right\}\,.
$$
Here $U$ is a matrix from $U(\infty)$ and $u$ ranges over the set of
its eigenvalues. All but finitely many $u$'s equal 1, so that the
product over $u$ is actually finite. The product over $i$ is
convergent, because the sum of the parameters is finite. Note also
that different $\om$'s correspond to different functions: here the
condition $\be_1^++\be_1^-\le1$ plays the decisive role, see
\cite{Ol3, Remark 1.6}.

\proclaim{Theorem 1.1} The functions $\chi^{(\om)}$, where $\om$ ranges
over $\Om$, are exactly the extreme characters of the group
$U(\infty)$.
\endproclaim

\demo{Proof} The fact that any $\chi^{(\om)}$ is an extreme character is
due to Voiculescu \cite{Vo}. The fact that the extreme characters are
exhausted by the $\chi^{(\om)}$'s can be proved in two ways: by reduction
to an old theorem due to Edrei \cite{Ed} (see \cite{Boy} and \cite{VK2}) and
by Vershik--Kerov's asymptotic method (see \cite{VK2} and
\cite{OkOl}). \qed
\enddemo

The coordinates $\al^\pm_i$, $\be^\pm_i$, and $\ga^\pm$ (or $\de^\pm$)
are called the {\it Voiculescu parameters\/} of the extreme character
$\chi^{(\om)}$. Theorem 1.1 is similar to
Thoma's theorem which describes the extreme characters of the
infinite symmetric group, see \cite{Th, VK1, Wa, KOO}. Another
analogous result is the classification of invariant ergodic
measures on the space of infinite Hermitian matrices (see \cite{OV} and
\cite{Pi2}).

\subhead (b) Spectral measures \endsubhead
Equip $\R^{4\infty+2}$ with the product topology. It induces a
topology on $\Om$. In this topology, $\Om$ is a locally compact
separable space. On the other hand, we equip the set of characters with
the topology of uniform convergence on the subgroups $U(N)\subset
U(\infty)$, $N=1,2,\dots$\,. One can prove that the bijection
$\om\longleftrightarrow\chi^{(\om)}$ is a homeomorphism with respect to
these two topologies, see \cite{Ol3, \S8}. In particular,
$\chi^{(\om)}(U)$ is a continuous function of $\om$ for any fixed $U\in
U(\infty)$.

\proclaim{Theorem 1.2} For any character $\chi$ of the group
$U(\infty)$ there exists a probability measure $P$ on the space $\Om$
such that
$$
\chi(U)=\int_\Om \chi^{(\om)}(U)\,P(d\om), \qquad U\in U(\infty).
$$
Such a measure $P$ is unique. The correspondence $\chi\mapsto P$ is a
bijection between the set of all characters and the set of all
probability measures on $\Om$.
\endproclaim

Here and in what follows, by a measure on $\Om$ we mean a
Borel measure. We call $P$ the {\it spectral measure\/} of $\chi$.

\demo{Proof} See \cite{Ol3, Theorem 9.1}. \qed
\enddemo

Similar results hold for the infinite symmetric group (see
\cite{KOO}) and for invariant measures on infinite Hermitian matrices
(see \cite{BO4}).

\subhead (c) Signatures \endsubhead
Define a {\it signature\/} $\la$ of length $N$ as an ordered sequence
of integers with $N$ members:
$$
\la=(\la_1\ge \la_2\ge\dots\ge \la_N\,|\,\la_i\in\Bbb Z).
$$
Signatures of length $N$ are naturally identified with highest
weights of irreducible representations of the group $U(N)$, see,
e.g., \cite{Zh}. Thus, there is a natural bijection
$\la\longleftrightarrow\chi^\la$ between signatures of length $N$ and
irreducible characters of $U(N)$ (here we use the term ``character''
in its conventional sense). The character $\chi^\la$ can be viewed as
a {\it rational Schur function} (Weyl's character formula)
$$
\chi^\la(u_1,\dots,u_N)
=\frac{\det [u_i^{\la_j+N-j}]_{i,j=1,\dots,N}}
{\det[u_i^{N-j}]_{i,j=1,\dots,N}}\,.
$$
Here the collection $(u_1,\dots,u_n)$ stands for the spectrum of a
matrix in $U(N)$.

We will represent a signature $\la$ as a pair of Young diagrams
$(\la^+,\la^-)$: one consists of positive $\la_i$'s, the other
consists of minus negative $\la_i$'s, zeros can go in either of the
two:
$$
\la=(\la_1^+,\la_2^+,\dots,-\la_2^-,\la_1^-).
$$
Let $d^+=d(\la)$ and $d^-=d(\la^-)$, where the symbol $d(\,\cdot\,)$
denotes the number of diagonal boxes of a Young diagram. Write the
diagrams $\la^+$ and $\la^-$ in Frobenius notation:
$$
\la^\pm=(p_1^\pm,\dots,p_{d^\pm}^\pm\mid q_1^\pm,\dots,q_{d^\pm}^\pm).
$$
We recall that the Frobenius coordinates $p_i$, $q_i$ of a Young
diagram $\nu$ are defined by
$$
p_i=\nu_i-i, \quad q_i=(\nu')_i-i, \qquad i=1,\dots, d(\nu),
$$
where $\nu'$ stands for the transposed diagram. Following
Vershik--Kerov, we introduce the {\it modified Frobenius
coordinates\/} of $\nu$ by
$$
\wt p_i=p_i+\tfrac12, \quad
\wt q_i=q_i+\tfrac12.
$$
Note that $\sum (\wt p_i+\wt q_i)=|\nu|$, where $|\nu|$ denotes the
number of boxes in $\nu$.

We agree that
$$
\wt p_i=\wt q_i=0, \qquad i>d(\nu).
$$

\subhead (d) Approximation of extreme characters \endsubhead
Recall that the dimension of the irreducible representation of $U(N)$
indexed by $\la$ is given by Weyl's dimension formula
$$
\Dim_N\la=\chi^\la(\,\underbrace{1,\dots,1}_N\,)
=\prod_{i\le i<j\le N}\frac{\la_i-i-\la_j+j}{j-i}\,.
$$
Define the {\it normalized\/} irreducible character indexed by $\la$
as follows
$$
\wt\chi^\la=\frac1{\Dim_N\la}\,\chi^\la.
$$
Clearly, $\wt\chi^\la(e)=1$.

Given a sequence $\{f_N\}_{N=1,2,\dots}$ of functions on the groups
$U(N)$, we say that $f_N$'s {\it approximate\/} a function $f$ defined on
the group $U(\infty)$ if, for any fixed $N_0=1,2,\dots$, the restrictions
of the functions $f_N$ (where $N\ge N_0$) to the subgroup $U(N_0)$
uniformly tend, as $N\to\infty$, to the restriction of $f$ to
$U(N_0)$.

\proclaim{Theorem 1.3} Any extreme character $\chi$ of $U(\infty)$
can be approximated by a sequence $\wt\chi^{\,(N)}$ of
normalized irreducible characters of the groups $U(N)$.

In more detail, write $\wt\chi^{\,(N)}=\wt\chi^{\,\la(N)}$, where
$\{\la(N)\}_{N=1,2,\dots}$ is a sequence of signatures, and let $\wt
p_i^{\,\pm}(N)$ and $\wt q_i^{\,\pm}(N)$ stand for the modified
Frobenius coordinates of $(\la(N))^\pm$. Then the functions
$\wt\chi^{\,(N)}$ approximate $\chi$ if and only if the following
conditions hold
$$
\lim_{N\to\infty}\frac{\wt p_i^{\,\pm}(N)}N=\al_i^\pm, \quad
\lim_{N\to\infty}\frac{\wt q_i^{\,\pm}(N)}N=\be_i^\pm, \quad
\lim_{N\to\infty}\frac{|(\la(N))^\pm|}N=\de^\pm,
$$
where $i=1,2,\dots$, and $\al_i^\pm$, $\be_i^\pm$, $\de^\pm$ are the
Voiculescu parameters of the character $\chi$.
\endproclaim

This claim reveals the asymptotic meaning of the Voiculescu
parameters. Note that for any
$\om=(\al^+,\be^+;\al^-,\be^-;\de^+,\de^-)\in\Om$, there exists a
sequence of signatures satisfying the above conditions, hence any
extreme character indeed admits an approximation.

\demo{Proof} This result is due to Vershik and Kerov, see their
announcement \cite{VK2}. A detailed proof is contained in \cite{OkOl}.
\qed
\enddemo

For analogous results, see \cite{VK1, OV}.

\head 2. Approximation of spectral measures  \endhead

\subhead (a) The graph $\GT$ \endsubhead
For two signatures $\nu$ and $\la$, of length $N-1$ and $N$,
respectively, write $\nu\prec\la$ if
$$
\la_1\ge\nu_1\ge\la_2\ge\nu_2\ge\dots\ge\nu_{N-1}\ge\la_N\,.
$$
The relation $\nu\prec\la$ appears in the {\it Gelfand--Tsetlin
branching rule\/} for the irreducible characters of the unitary
groups, see, e.g., \cite{Zh}:
$$
\chi^\la(u_1,\dots,u_{N-1},1)=\sum_{\nu:\,\nu\prec\la}\chi^\nu\,.
$$

The {\it Gelfand--Tsetlin graph\/} $\GT$ is a $\Z_+$--graded graph
whose $N$th level $\GT_N$ consists of signatures of length $N$.
Two vertices $\nu\in\GT_{N-1}$ and $\la\in\GT_N$ are joined by an edge if
$\nu\prec\la$. This graph is a counterpart of the Young graph
associated with the symmetric group characters \cite{VK1},
\cite{KOO}.

\subhead (b) Coherent systems of distributions \endsubhead
For $\nu\in\GT_{N-1}$ and $\la\in\GT_N$, set
$$
q(\nu,\la)=\cases
\dfrac{\Dim_{N-1}\nu}{\Dim_N\la}\,, & \nu\prec\la,\\
0, & \nu\nprec\la. \endcases
$$
This is the {\it cotransition probability function\/} of the
Gelfand--Tsetlin graph. It satisfies the relation
$$
\sum_{\nu\in\GT_{N-1}}q(\nu,\la)=1, \qquad \forall\, \la\in\GT_N\,.
$$

Assume that for each $N=1,2,\dots$ we are given a probability measure
$P_N$ on the discrete set $\GT_N$. Then the family
$\{P_N\}_{N=1,2,\dots}$ is called a {\it coherent system\/} if
$$
P_{N-1}(\nu)=\sum_{\la\in\GT_N}q(\nu,\la)P_N(\la),
\qquad \forall\, N=2,3,\dots,\quad \forall\, \nu\in\GT_{N-1} \,.
$$
Note that if $P_N$ is an arbitrary probability measure on $\GT_N$
then this formula defines a probability measure on $\GT_{N-1}$
(indeed, this follows at once from the above relation for
$q(\nu,\la)$). Thus, in a coherent system $\{P_N\}_{N=1,2,\dots}$, the
$N$th term is a refinement of the $(N-1)$st one.

\proclaim{Proposition 2.1} There is a natural bijective correspondence
$\chi\longleftrightarrow \{P_N\}$ between characters of the group
$U(\infty)$ and coherent systems, defined by the relations
$$
\chi\mid_{U(N)}=\sum_{\la\in\GT_N}P_N(\la)\wt\chi^\la,
\qquad N=1,2,\dots\,.
$$
\endproclaim

\demo{Proof} See \cite{Ol3, Proposition 7.4}. \qed
\enddemo

A similar claim holds for the infinite symmetric group
$S(\infty)$, see \cite{VK1}, \cite{KOO}, and for the infinite--dimensional
Cartan motion group, see \cite{OV}. Note that $\{P_N\}$ can be viewed as a
kind of Fourier transform of the corresponding character.

The concept of a coherent system $\{P_N\}$ is important for two
reasons. First, we are unable to calculate directly the ``natural''
nonextreme characters but we dispose of nice closed expressions for their
``Fourier coefficients'' $P_N(\la)$, see the next section. Note that
in the symmetric group case the situation is just the same, see
\cite{KOV}, \cite{BO1-3}. Second, the measures $P_N$ approximate the
spectral measure $P$, see below.

\subhead (c) Approximation $P_N\to P$ \endsubhead
Let $\chi$ be a character of $U(\infty)$ and let $P$ and $\{P_N\}$ be the
corresponding spectral measure and coherent system.

For any $N=1,2,\dots$, we embed the set $\GT_N$ into
$\Om\subset\R^{4\infty+2}$ as follows
$$
\gather
\GT_N\ni\la \longmapsto (a^+,b^+;a^-,b^-;c^+,c^-)\in\R^{4\infty+2},\\
a^\pm_i=\frac{\wt p_i^\pm}N, \quad
b^\pm_i=\frac{\wt q_i^\pm}N, \quad
c^\pm=\frac{|\la^\pm|}N\,,
\endgather
$$
where $i=1,2,\dots$, and $\wt p^\pm_i$, $\wt q^\pm_i$ are the modified
Frobenius coordinates of $\la^\pm$.

Let $\un P_N$ be the pushforward of $P_N$ under this
embedding. Then $\un P_N$ is a probability measure on $\Om$.

\proclaim{Theorem 2.2} As $N\to\infty$, the measures $\un P_N$
weakly tend to the measure $P$. I.e., for any
bounded continuous function $F$ on $\Om$,
$$
\lim_{N\to\infty}\langle F, \un P_N\rangle
=\langle F, P\rangle.
$$
\endproclaim

\demo{Proof} See \cite{Ol3, Theorem 10.2}.
\qed
\enddemo

This result should be compared with \cite{KOO, Proof of Theorem B in \S8} and
\cite{BO4, Theorem 5.3}. Its proof is quite similar to that of
\cite{BO4, Theorem 5.3}.

Theorem 2.2 shows that the spectral measure can be, in principle, computed if
one knows the coherent system $\{P_N\}$.

\head 3. ZW--measures \endhead

The goal of this section is to introduce a family of characters
$\chi$ of the group $U(\infty)$, for which we solve the problem of
harmonic analysis. We describe these characters in terms of the
corresponding coherent systems $\{P_N\}$. For detailed proofs we
refer to \cite{Ol3}.

Let $z,z',w,w'$ be complex parameters.
For any $N=1,2,\dots$ and any $\la\in\GT_N$ set
$$
\gathered
P'_N(\la\mid z,z'w,w')=\Dim_N^2(\la)\\
\times \prod_{i=1}^N
\frac1{\Gamma(z-\la_i+i)\Gamma(z'-\la_i+i)
\Gamma(w+N+1+\la_i-i)\Gamma(w'+N+1+\la_i-i)}\,,
\endgathered
$$
where $\Dim_N\la$ was defined in \S1. Clearly, for any fixed $N$ and $\la$,
$P'_N(\la\mid z,z',w,w')$ is an entire function on $\C^4$.
Set
$$
\Cal D=\{(z,z',w,w')\in\C^4\mid \Re(z+z'+w+w')>-1\}.
$$
This is a domain in $\C^4$.

\proclaim{Proposition 3.1} Fix an arbitrary $N=1,2,\dots$\,.
The series of entire functions
$$
\sum_{\la\in\GT_N} P'_N(\la\mid z,z',w,w')
$$
converges in the domain $\Cal D$, uniformly on compact sets. Its sum
is equal to
$$
S_N(z,z'w,w')=\prod_{i=1}^N
\frac{\Gamma(z+z'+w+w'+i)}
{\Gamma(z+w+i)\Gamma(z+w'+i)\Gamma(z'+w+i)\Gamma(z'+w'+i)\Gamma(i)}
$$
\endproclaim

\demo{Proof} See \cite{Ol3, Proposition 7.5}. \qed
\enddemo

Note that in the special case $N=1$, the set $\GT_1$ is simply $\Z$
and the identity
$$
\sum_{\la\in\GT_1}P'_1(\la\mid\zw)=S_1(\zw)
$$
is equivalent to the well--known Dougall's formula (see
\cite{Er, vol 1, \S1.4}).

Consider the subdomain
$$
\gather
\Cal D_0=\{(z,z',w,w')\in\Cal D\mid
z+w,\,z+w',\,z'+w,\,z'+w'\ne-1,-2,\dots\}\\
=\{(\zw)\in\Cal D\mid S_N(z,z',w,w')\ne0\}.
\endgather
$$
For any $(\zw)\in\Cal D_0$ we set
$$
P_N(\la\mid\zw)=\frac{P'_N(\la\mid\zw)}{S_N(\zw)}\,, \qquad
N=1,2,\dots,\quad \la\in\GT_N.
$$
Then, by Proposition 3.1,
$$
\sum_{\la\in\GT_N}P_N(\la\mid\zw)=1, \qquad (\zw)\in\Cal D_0\,,
$$
uniformly on compact sets in $\Cal D_0$.

\proclaim{Proposition 3.2} Let $(\zw)\in\Cal D_0$. For any
$N=2,3,\dots$, the coherency relation of \S2(b) is satisfied,
$$
P_{N-1}(\nu\mid\zw)=\sum_{\la\in\GT_N}q(\nu,\la)P(\la\mid\zw).
$$
\endproclaim

\demo{Proof} See \cite{Ol3, Proposition 7.7}. \qed
\enddemo

Combining this with Proposition 2.1 we conclude that  $\{P_N(\,\cdot\,\mid 
\zw)\}$, where $N=1,2,\dots$,
is a coherent system provided that $(\zw)\in\Cal D_0$ satisfies the
{\it positivity condition\/}: for any $N=1,2,\dots$, the expression
$P'_N(\la\mid \zw)$ is nonnegative for all $\la\in\GT_N$. (Note that
there always exists $\la$ for which $P'_N(\la\mid \zw)\ne0$, because
the sum over $\la$'s is not 0.)
We proceed to describing a set of quadruples $(\zw)\in\Cal D_0$
satisfying the positivity condition.

Define the subset $\Cal Z\subset\C^2$ as follows:
$$
\gather
\Cal Z=\Zp\sqcup\Zc\sqcup\Zd, \\
\Zp=\{(z,z')\in\C^2\setminus\R^2\mid z'=\bar z\}, \\
\Zc=\{(z,z')\in\R^2\mid \exists m\in\Z, \, m<z,z<m+1\}, \\
\Zd=\underset{m\in\Z}\to{\sqcup}\Zdm, \\
\Zdm=\{(z,z')\in\R^2\mid z=m, \,z'>m-1, \quad
\text{or} \quad z'=m,\, z>m-1\},
\endgather
$$
where ``princ'', ``compl'', and ``degen'' are abbreviations for
``principal'', ``complementary'', and ``degenerate'', respectively.
For an explanation of this terminology, see \cite{Ol3}.

\proclaim{Proposition 3.3} Let $(z,z')\in\C^2$.

{\rm(i)} The expression
$(\Gamma(z-k)\Gamma(z'-k))^{-1}$ is nonnegative for all $k\in\Z$ if
and only if $(z,z')\in\Cal Z$.

{\rm(ii)} If $(z,z')\in\Zp\sqcup\Zc$ then this expression is
strictly positive for all $k\in\Z$.

{\rm(iii)} If $(z,z')\in\Zdm$ then this expression vanishes for
$k=m,m+1,\dots$ and is strictly positive for $k=m-1,m-2,\dots$\,.
\endproclaim

\demo{Proof} See \cite{Ol3, Lemma 7.9}. \qed
\enddemo

\example{Definition 3.4} The {\it set of admissible values\/} of the
parameters $\zw$ is the subset $\Dadm\subset\Cal D$ of quadruples
$(\zw)$ such that both $(z,z')$ and $(w,w')$ belong to
$\Cal Z$. When both $(z,z')$ and $(w,w')$ are in $\Cal Z$,
an extra condition is added: let $k,l$ be such that $(z,z')\in\Zdk$ and
$(w,w')\in\Zdl$; then we require $k+l\ge0$. A quadruple $(\zw)$ will
be called {\it admissible\/} if it belongs to the set $\Dadm$.
\endexample

Note that in this definition we do not assume {\it a priori\/} that
$(\zw)$ belongs to the subdomain $\Cal D_0\subset\Cal D$. However the
conditions imposed on $(\zw)$ imply that $\Dadm\subset\Cal D_0$, see below.

\proclaim{Proposition 3.5} Let $(\zw)\in\Dadm$ and let
$N=1,2,\dots$\,. Then $P'_N(\la\mid\zw)\ge0$ for any $\la\in\GT_N$,
and there exists $\la\in\GT_N$ for which the above inequality is
strict.
\endproclaim

\demo{Proof} The first claim follows from Proposition 3.3 (i). Now we
shall describe the set of those $\la\in\GT_N$ for which
$P'_N(\la\mid\zw)>0$.

When both $(z,z')$ and $(w,w')$ are in $\Zp\sqcup\Zc$ then, by
Proposition 3.3 (ii), this is the whole $\GT_N$.

When $(w,w')\in\Zp\sqcup\Zc$ and $(z,z')\in\Cal Z$, say,
$(z,z')\in\Zdm$, then this set is formed by $\la$'s satisfying the
condition $\la_1\le m$. Indeed, this readily follows from claims (ii)
and (iii) of Proposition 3.3.

Likewise, when $(z,z')\in\Zp\sqcup\Zc$ and $(w,w')\in\Zdm$, then the
condition takes the form $\la_N\ge-m$.

Finally, when both $(z,z')$ and $(w,w')$ are in $\Cal Z$, say
$(z,z')\in\Zdk$ and $(w,w')\in\Zdl$, then the set in question is
described by the conditions $\la_1\le k$, $\la_N\ge -l$. The set is
nonempty provided that $k\ge-l$, which is exactly the extra condition
from Definition 3.4.

Note that if $k=-l$ then this set consists of a single element
$\la=(k,\dots,k)$. \qed
\enddemo

Proposition 3.5 implies that $\Dadm\subset\Cal D_0$. Of course, this
can be checked directly, but the claim is not entirely obvious, for
instance, when both $(z,z')$ and $(w,w')$ are in $\Zc$.

Now we can summarize the above definitions and results in the
following theorem.

\proclaim{Theorem 3.6} For any admissible quadruple $(\zw)$, the family
$\{P_N(\,\cdot\,\mid\zw)\}$, where $N=1,2,\dots$, is a coherent system, so that
it determines a character $\chi_\zw$ of the group $U(\infty)$.
\endproclaim

\demo{Proof} Indeed, let $(\zw)$ be admissible. Since $(\zw)$ is in
$\Cal D_0$, the definition of $P_N(\la\mid \zw)$'s makes sense. By
Proposition 3.5, for any $N$, $P_N(\,\cdot\,\mid \zw)$ is a probability
distribution on $\GT_N$. By Proposition 3.2, the family
$\{P_N(\,\cdot\,\mid \zw)\}_{N=1,2,\dots}$ is a coherent system. By
Proposition 2.1, it defines a character of $U(\infty)$.
\qed
\enddemo

\example{Remark 3.7} The set of characters of the form $\chi_\zw$ is
stable under tensoring with one--dimensional characters
$(\det(\,\cdot\,))^k$, where $k\in\Z$. Indeed, the sets $\Cal D$,
$\Cal D_0$, and $\Dadm$ are invariant under the shift
$$
(\zw)\mapsto(z+k,z'+k,w-k,w'-k),
$$
and we have
$$
P_N(\la+(k,\dots,k)\mid\zw)=P_N(\la\mid z+k,z'+k,w-k,w'-k).
$$
On the other hand, in terms of coherent systems, tensoring with
$(\det(\,\cdot\,))^k$ is equivalent to shifting $\la$ by
$(k,\dots,k)$.
\endexample

\example{Remark 3.8} In the special case when both $(z,z')$ and
$(w,w')$ are in $\Zd$, a detailed study of the distributions
$P_N(\,\cdot\,\mid \zw)$ from a combinatorial point of view was
given by Kerov \cite{Ke}.
\endexample

\example{Remark 3.9} As we see, the structure of the set of all admissible 
parameters is fairly complicated. However, all the major formulas that will be 
obtained below do not feel this structure; they hold for all admissible 
parameters. The explanation of this phenomenon
is rather simple: the quantities in question (like correlation functions)
can usually be defined for the parameters varying in domain which is much larger 
than $\Dadm$, see e.g. Propositions 3.1 and 3.2 above. Thus, the formulas for 
these quantities usually hold on an open subset of $\C^4$ containing $\Dadm$. It 
is only when we require certain quantities to be positive in order to
fit our computations in the framework of probability theory, we need to
restrict ourselves to the smaller set of admissible parameters.
\endexample

\head 4. Two discrete point processes \endhead

In this section we will explain two different ways to associate to the
measure $P_N$ introduced in the previous section a discrete point process.
We also show how the two resulting processes can be obtained one from the
other.

First, we recall the general definition of a random point process.

Let $\X$ be a locally compact separable topological space. A {\it
multiset\/} $X$ in $\X$ is a collection of points with possible
multiplicities and with no ordering
imposed. A {\it locally finite point configuration\/} ({\it
configuration,\/} for short) is a
multiset $X$ such that for any compact set $A\subset\X$ the
intersection $X\cap A$ is finite (with multiplicities counted). This
implies that $X$ itself is either finite or countably infinite.

The set of all configurations in $\X$ is denoted by $\Conf(\X)$. Given a
relatively compact Borel set $A\subset\X$, we introduce a function $\N_A$ on
$\Conf(\X)$ by setting $\N_A(X)=|X\cap A|$. We equip $\Conf(\X)$ with
the Borel structure generated by all functions of this form.

A {\it random point process\/} on $\X$ (point process, for short;
another term is random point field) is a probability Borel measure
$\Cal P$ on the space $\Conf(\X)$.

We do not need the full generality of these definitions in this section. Here
the situation is rather simple: all our processes are discrete (that is, the
space $\x$ is discrete), and the point configurations are finite. However, in
\S9 we will consider a continuous point process with infinitely many
particles, and then we will need the above definitions.

Consider the lattice
$$
\x=\x^{(N)}=\cases \Z, &N\text{ is odd},\\
                   \Z+\frac 12, &N \text{ is even},
		  \endcases
$$
and divide it into two parts
$$
\gathered
\x=\xin\sqcup\xout,\\
\xin=\left\{-\frac{N-1}2,-\frac{N-3}2,\dots,\frac{N-3}2,
\frac{N-1}2\right\},\quad |\xin|=N,\\
\xout=\left\{\dots,-\frac{N+3}2,-
\frac{N+1}2\right\}\sqcup\left\{\frac{N+1}2,\frac{N+3}2,\dots\right\},\quad 
|\xout|=\infty.
\endgathered
$$
Let $\rho_i=\frac{N+1}2-i$, $i=1,\dots, N$. For any $\la\in\GT_{N}$ we
set
$$
\cl(\la)=\{\la_1+\rho_1,\dots,\la_N+\rho_N\}.
$$
Clearly,
$\la\mapsto\cl(\la)$ defines a bijection between $\GT_N$ and the set
of $N$-point multiplicity free configurations on $\X$.

Now we define another correspondence between signatures and point
configurations. Let us represent a signature $\la$ as a pair of Young diagrams
$(\la^+,\la^-)$, see \S1(c).

Finally, we define a point configuration as
$$
X(\la)=\left\{p_i^++\frac{N+1}2\right\}\sqcup
\left\{\frac{N-1}2-q_i^+\right\}\sqcup
\left\{-p_j^--\frac{N+1}2\right\}\sqcup
\left\{-\frac{N-1}2+q_j^-\right\},
\tag 4.1
$$
where $i=1,\dots,d^+$ and $j=1,\dots,d^-$, see \S1(c) for the notation. Note
that if $\la=0$ then the configuration is empty.

{}From the inequalities
$$
\gathered
p_1^+>\dots> p_{d^+}^+\ge0, \quad
q_1^+>\dots> q_{d^+}^+\ge0, \quad
\\
p_1^->\dots> p_{d^-}^-\ge0, \quad
q_1^->\dots> q_{d^-}^-\ge0,\\
d^++d^-\le N
\endgathered
$$
it follows that $X(\la)$ consists of an even number of distinct points (equal
to $2(d^++d^-)$), of which a half lies in $\xout$ while another half
lies in $\xin$. Finite point configurations with this property will
be called {\it balanced.}

Conversely, each balanced, multiplicity free configuration on $\X$ is of the
form $X(\la)$ for one and only one signature $\la\in\GT_N$. Thus, the map
$\la\mapsto\cl(\la)$ defines a bijection between $\GT_N$ and the set
of finite balanced configurations on $\X$ with no multiplicities.

Define an involution on the set $\Conf(\X)$ of multiplicity free point
configurations on $\x$ by
$$
X\mapsto X^\triangle=X\,\triangle\,\xin
=(X\cap\xout)\cup(\xin\setminus X).
$$
Since $|\xin|=N$,
this involution defines a bijection between
$N$-point configurations and finite balanced configurations.

\proclaim{Proposition 4.1} In the above notation,
$X(\la)=\cl(\la)^\triangle$ for any signature $\la\in\GT_N$.
\endproclaim

{}For instance, let $N=7$ and $\la=(4,2,2,0,-1,-2,-2)$. Then
$$\cl(\la)=\{7,4,3,0,-2,-4,-5\}$$ and, since $\xin=\{3,2,1,0,-1,-2,-3\}$, we 
have
$$
X(\la)=\cl(\la)^\triangle=\{7,4,2,1,-1,-3,-4,-5\}.
$$
On the other hand, $\la^+=(4,2,2)=(3,0\,|\, 2,1)$,
$\la^-=(2,2,1)=(1,0\,|\, 2,0)$, and (4.1) gives the same $X(\la)$.

Another example: for the zero signature $\underline0$ we have
$\cl(\underline0)=\xin$ and $X(\underline0)=\varnothing$.

\demo{Proof of Proposition 4.1} Geometric constructions described below are
illustrated by a picture at the end of the text. Consider a plane with
Cartesian coordinates $(x,y)$ and  put the lattice $\X=\X^{(N)}$ on the
vertical axis $x=0$, so that each  $a\in\X$ is identified with the point
$(0,a)$ of the plane. Draw a square grid in the plane, formed by the
horizontal lines  $y=a+\frac12$, where $a$ ranges over $\X$, and by the
vertical lines $x=b$, where $b$ ranges over $\Z$. We represent $\la$ by an
infinite polygonal line $\Bbb L$ on the grid, as follows.

Denote by $A_0,\dots,A_N$ the horizontal lines defined by
$y=\frac{N}2$, $y=\frac{N}2-1$, $\dots$, $y=-\frac{N}2$,
respectively. We remark that these lines belong to the grid: indeed,
$\X$ coincides with $\Z$ shifted by $\frac{N-1}2$, so that the
points $\frac{N}2,\frac{N}2-1,\dots,-\frac{N}2$ belong to the shift
of $\X$ by $\frac12$.

The polygonal line $\Bbb L$ first goes along
$A_0$, from right to left, starting at $x=+\infty$, up to the point
with the coordinate $x=\la_1$. Then it changes the direction and goes
downwards until it meets the next horizontal line $A_1$. Then it goes along
$A_1$, again from right to left, up to the point with the coordinate
$x=\la_2$, etc. Finally, after reaching the lowest line $A_N$ at the point
with the coordinate $x=\la_N$, it goes only to the left, along this line.

Further, we define a bijective correspondence $a\leftrightarrow s$
between the points $a\in\X$ and the sides $s$ of $\Bbb L$, as
follows. Given $a$, we draw the line $x+y=a$, it intersects $\Bbb L$ at
the midpoint of a side, which is, by definition, $s$. Let us call $a$
a ``$v$-point'' or a ``$h$-point'' according to whether the
corresponding side $s$ is vertical or horizontal. Thus, the whole set
$\X$ is partitioned into ``$v$-points'' and ``$h$-points''.

The ``$v$-points'' of $\X$ are exactly those of
the configuration $\cl(\la)$. Consequently, the collection
$\bigl(\cl(\la)\cap\xout\bigr)\sqcup\bigl(\xin\setminus\cl(\la)\bigr)$
is formed by the ``$v$-points'' from $\xout$ and the ``$h$-points''
from $\xin$.

On the other hand, the correspondence $a\leftrightarrow s$ makes it
possible to interpret the same collection of points in terms of the
Frobenius coordinates of the diagrams $\la^+$ and $\la^-$. Indeed, the diagram
$\la^+$ can be identified with the figure bounded by the horizontal line
$A_0$, the vertical line $x=0$, and by $\Bbb L$. Then the line $x+y=\frac{N}2$
coincides with the  diagonal of $\la^+$. Above this line, there are $d^+$
vertical sides of $\Bbb L$, say, $s_1,\dots,s_{d^+}$, which lie in the rows of
$\la^+$ with numbers $1,\dots, d^+$. The corresponding Frobenius coordinates
are $p^+_i=\la_i-i$, where $i=1,\dots,d^+$. It easily follows that the
midpoint of the side $s_i$ lies on the line $x+y=\frac{N+1}2+p^+_i$, i.e.,
$s_i$ corresponds to $\frac{N+1}2+p^+_i$. In this way we get the first
component $\left\{p_i^++\frac{N+1}2\right\}\subset X(\la)$, see (4.1).  The
remaining three components are interpreted similarly. \qed

\enddemo

Fix any admissible quadruple $(z,z',w,w')$ of parameters and consider
the corresponding probability measure $P_N$ on $\GT_N$, see \S3. Taking the 
pushforwards of the
measure $P_N$ under the maps $\la\mapsto\cl(\la)$ and $\la\mapsto
X(\la)$ we get two point processes on the lattice $\X$, which we
denote by $\tP^{(N)}$ and $\Cal P^{(N)}$, respectively. We are mainly
interested in the process $\Cal P^{(N)}$, which is defined by $\la\mapsto
X(\la)$; the process $\tP^{(N)}$ defined by $\la\mapsto\cl(\la)$ will play an
auxiliary role. Proposition 4.1 implies that
$\tP^{(N)}(X)=\Cal P^{(N)}(X^\triangle)$ for any finite configuration $X$.

\head 5. Determinantal point processes. General theory
\endhead

\subhead (a) Correlation measures \endsubhead
Let $\Cal P$ be a point process on $\X$ (see the definition in the beginning
of \S4), and let $A$ denote an arbitrary  relatively compact Borel subset of
$\X$. Then $\N_A$ is a random variable with  values in $\{0,1,2,\dots\}$. We
assume that for any $A$ as above, $\N_A$ has finite moments of all orders.

Let $n$ ranges over $\{1,2,\dots\}$. The {\it $n$th correlation measure} of
$\Cal P$,  denoted as $\rho_n$, is a Borel measure on $\X^n$, uniquely
defined by
$$
\rho_n(A^n)=\E[\N_A(\N_A-1)\dots(\N_A-n+1)],
$$
where the symbol $\E$ means expectation with respect to the probability space
$(\Conf(\X),\Cal P)$.

Equivalently, for any bounded compactly supported Borel function $F$ on $\X^n$,
$$
\langle F, \rho_n\rangle=
\int_{X\in\Conf(X)}\left(\sum\Sb x_1,\dots,x_n\in X\\
\text{pairwise distinct}\endSb F(x_1,\dots,x_n)\right) \Cal P(dX),
$$
where the summation is taken over all ordered $n$--tuples of pairwise distinct
points  taken from the (random) configuration $X$ (here a multiple point is
viewed as a  collection of different elements).

The measure $\rho_n$ takes finite values on the compact subset of $\X^n$.
The measure $\rho_n$ is symmetric with respect to the permutations of the 
arguments.

Under mild assumptions about the growth of $\rho_n(A^n)$ as $n\to\infty$ (here
$A$ is an arbitrary compact set), the collection of the correlation measures
$\rho_1,\rho_2,\dots$ defines the initial process $\Cal P$ uniquely. See
\cite{Len} and \cite{So, (1.6)}.

When there is a ``natural'' reference measure $\mu$ on $\X$ such
that, for any $n$, $\rho_n$
is absolutely continuous with respect to the product measure $\mu^{\otimes n}$,
the density of $\rho_n$ is  called the {\it $n$th correlation function}. For
instance, this always holds if the  space $\X$ is discrete: then as $\mu$ one
takes the counting measure on $\X$.  The correlation functions are denoted as
$\rho_n(x_1,\dots,x_n)$.

If the space $\x$ is discrete and the process is multiplicity free then
$\rho_n(x_1,\dots,x_n)$ is the probability that the random point configuration
contains the points $x_1,\dots,x_n$ (here $x_i$'s are pairwise distinct,
otherwise $\rho_n(x_1,\dots,x_n)=0$).

For a general discrete process,
$\rho_n(x_1,\dots,x_n)$ is equal to the sum of weights of the point
configurations with certain combinatorial prefactors computed as follows: if
$x$ has multiplicity $k$ in the multiset $(x_1,\dots,x_n)$ and has multiplicity
$m$ in the point configuration in question, then this produces the prefactor
$m(m-1)\cdots(m-k+1)$ (such prefactor is computed for every element of the set
$\{x_1,\dots,x_n\}$). Note that this prefactor vanishes unless $m\ge k$ for
every $x\in \{x_1,\dots,x_n\}$.

\subhead (b) Determinantal processes \endsubhead
A point process is called {\it determinantal} if there exists a function
$K(x,y)$ on  $\X\times\X$ such that, for an appropriate reference measure
$\mu$, the  correlation functions are given by the determinantal formula
$$
\rho_n(x_1,\dots,x_n)=\det[K(x_i,x_j)]_{i,j=1}^n\,, \qquad n=1,2,\dots\,.
$$

The function $K$ is called the {\it correlation kernel} of the process. It is
not  unique: replacing $K(x,y)$ by $f(x)K(x,y)f(y)^{-1}$,  where $f$ is an
arbitrary  nonzero function on $\X$, leaves the above expression for the
correlation  functions intact.

If the reference measure is multiplied by a positive function $f$ then the
correlation kernel should be appropriately transformed. For instance, one can
multiply it by $(f(x)f(y))^{-1/2}$.

It is often useful to view $K(x,y)$ as the kernel of an integral operator acting
in the Hilbert space $L^2(\X,\mu)$. We will denote this operator by the same
symbol $K$.

Assume that a function $K(x,y)$ is Hermitian symmetric (i.e.,
$K(x,y)=\overline{K(y,x)})$ and locally of trace class (i.e., its
restriction to
any compact set $A\subset\X$ defines a trace class operator in $L^2(A,\mu)$,
where $\mu$ is a fixed reference measure). Then $K(x,y)$ is the correlation
kernel of a determinantal point process if and only if the operator $K$ in
$L^2(\X,\mu)$ satisfies the condition $0\le K\le 1$, see \cite{So}.
However, there are important examples of correlation kernels which are not
Hermitian symmetric, see below.

If $\X$ is a discrete countably infinite space then any multiset with
finite multiplicities is a configuration. As $\mu$ we will always take the
counting measure. A correlation kernel is simply an infinite matrix with the
rows and columns labeled by the points of $\X$. For any determinantal
process
on $\X$ the random configuration is multiplicity free with probability 1.
Indeed, if any two arguments of the $n$th correlation function coincide then
the defining determinant above vanishes.

\subhead (c) The complementation principle \endsubhead
Assume that $\X$ is discrete and fix a subset $Z\subseteq\X$. For a
subset $X$
in $\X$ let $X\triangle Z$ denote its symmetric difference with $Z$, i.e.,
$X\triangle Z=(X\cap \bar Z)\cup(Z\setminus X)$, where $\bar Z=\X\setminus Z$.
The map $X\mapsto X\triangle Z$, which we will denote by the symbol $\triangle$,
is an involution on multiplicity free configurations. If the process $\Cal P$
lives on the multiplicity free configurations, we can define its image $\Cal
P^\triangle$ under $\triangle$.

Assume further that $\Cal P$ is determinantal and let $K$ be its correlation
kernel. Then the process $\Cal P^\triangle$ is also determinantal. Its
correlation kernel $K^\triangle$ can be obtained from $K$ as follows:
$$
K^\triangle(x,y)=\cases K(x,y), & x\in\bar Z, \\
\de_{xy}-K(x,y), & x\in Z,
\endcases \tag 5.1
$$
where $\de_{xy}$ is the Kronecker symbol. See \cite{BOO, \S{A.3}}.

Note that one could equally well  use the formula
$$
K^\triangle(x,y)=\cases K(x,y), & y\in\bar Z, \\
\de_{xy}-K(x,y), & y\in Z,
\endcases
$$
obtained from \tht{5.1} by multiplying the kernel by the function
$\ep(x)\ep(y)$, where $\ep(\,\cdot\,)$ is equal to 1 on $\bar Z$ and
to $-1$ on $Z$.  This operation does not affect the correlation
functions, see \S5(b).

We call the passage from the process $\Cal P$ to the process $\Cal P^\triangle$,
together with formula \tht{5.1} the {\it complementation principle}. The idea
was  borrowed from unpublished work notes by Sergei Kerov connected with an
early  version of \cite{BOO}.

Note that Proposition 4.1 can now be restated as follows:
$$
\tP^{(N)}=(\Cal P^{(N)})^\triangle,
$$
where the role of the set $Z$ is played by $\xin$.

\subhead (d) Discrete polynomial ensembles \endsubhead
Here we assume that $\X$ is a finite or countably infinite subset of $\R$
without limit points.

Assume that we are given a nonnegative function $f(x)$ on $\X$. Fix a natural
number $N$. We consider $f$ as a weight function: denoting by $\mu$ the counting
measure on $\X$ we assign to $f$ the measure $f\mu$ on $\X$.

We impose on $f$ two basic assumptions:

(*) $f$ has finite moments at least up to order $2N-2$, i.e.,
$$
\sum_{x\in\X}x^{2N-2}f(x)<\infty.
$$

(**) $f$ does not vanish at least at $N$ distinct points.

Under these assumptions the functions $1,x,\dots, x^{N-1}$ on $\X$ are linearly
independent and lie in the Hilbert space $L^2(\X,f\mu)$. Let
$p_0=1,\,p_1,\,\dots,\,p_{N-1}$ be the monic polynomials obtained by
orthogonalizing  the system $(1,x,\dots, x^{N-1})$ in $L^2(\X,f\mu)$.

We set
$$
h_n=(p_n,p_n)_{L^2(\X,f\mu)}=\sum_{x\in\X}p_n^2(x)f(x),
\qquad n=0,\dots,N-1,
$$
and consider the {\it Christoffel--Darboux kernel}
$$
\sum_{n=0}^{N-1}\frac{p_n(x)p_n(y)}{h_n}\,,
\qquad x,y\in\X.
$$
This kernel defines an orthogonal projection operator in $L^2(\X,f\mu)$; its
range is the $N$--dimensional subspace spanned by $1,x,\dots, x^{N-1}$.

Consider an isometric embedding $L^2(\X,f\mu)\to\ell^2(\X)$ which is defined
as multiplication by $\sqrt{f(\,\cdot\,)}$. Under this isomorphism the
Christoffel--Darboux kernel turns into another kernel which we will call the
{\it normalized} Christoffel--Darboux kernel and denote as $K^\cd$:
$$
K^\cd(x,y)=\sqrt{f(x)f(y)}\cdot\sum_{n=0}^{N-1}\frac{p_n(x)p_n(y)}{h_n}\,,
\qquad x,y\in\X.
\tag 5.2
$$
This kernel defines a projection operator in $\ell^2(\X)$ of rank $N$.

Let $\Conf_N(\X)$ denote the set of $N$--point multiplicity free configurations
(subsets) in $\X$.  For $X\in\Conf_N(\X)$ we set
$$
V^2(X)=\prod_{1\le i<j\le N} (x_i-x_j)^2,
$$
where $x_1,\dots,x_N$ are the points of $X$ written in any order.

Under the assumptions (*) and (**) we have
$$
0<\sum_{X\in\Conf_N(\X)}\left(\prod_{x\in X}f(x)\cdot V^2(X)\right)<\infty.
$$
Therefore, we can form a point process on $\X$ which lives on $\Conf_N(\X)$ and
for which the probability of a configuration $X$ is given  by
$$
\operatorname{Prob}(X)=\const\cdot\prod_{x\in X}f(x)\cdot V^2(X),
\qquad X\in\Conf_N(\X),
\tag 5.3
$$
where $\const$ is the normalizing constant.
This process is called the $N$--point {\it polynomial ensemble} with the
weight  function $f$.

\proclaim{Proposition 5.1} Let $\X$ and $f$ be as above, and $f$
satisfy the assumptions \tht{*}, \tht{**}. Then the $N$--point
polynomial ensemble with
the weight  function $f$ is a determinantal point process whose correlation
kernel is the  normalized Christoffel--Darboux kernel \tht{5.2}.
\endproclaim

\demo{Proof} A standard argument from the Random Matrix Theory, see, e.g.,
\cite{Me, \S5.2}. \qed
\enddemo

\example{Remark 5.2} Under a stronger than (*) condition
$$
\sum_{x\in\X}|x|^{2N-1} f(x)<\infty,
$$
there exists a monic polynomial $p_N$ of degree $N$, orthogonal to
$1,x,\dots,x^{N-1}$ in $L^2(\X,f\mu)$. Then the Christoffel--Darboux
kernel can be written as
$$
\frac1{h_{N-1}}\,\frac{p_N(x)p_{N-1}(y)-p_{N-1}p_N(y)}{x-y}\,.
$$
The value at the diagonal $x=y$ is determined via the L'Hospital rule.

According to this, the normalized Christoffel--Darboux kernel can be written
in the form
$$
K^\cd(x,y)=\frac{\sqrt{f(x)f(y)}}{h_{N-1}}\,
\frac{p_N(x)p_{N-1}(y)-p_{N-1}p_N(y)}{x-y}\,.
\tag 5.4
$$
\endexample

\subhead (e) L--ensembles \endsubhead
Let $\X$ be an arbitrary discrete space (finite or countably
infinite). We are dealing with the Hilbert space
$\ell^2(\X)=L^2(\X,\mu)$, where, as usual, $\mu$
denotes the counting measure on $\X$. Let $\Conff(\X)$ denote the set of all
finite, multiplicity free configurations in $\X$ (i.e., simply finite
subsets).

Let $L$ be an operator in $\ell^2(\X)$ and $L(x,y)$ be its matrix
($x,y\in\X$). For $X\in\Conff(\X)$ we denote by $L_X(x,y)$ the
submatrix of $L(x,y)$ of order $|X|$ whose rows and columns are
indexed by the points $x\in X$. The determinants $\det L_X$ are
exactly the diagonal minors of the matrix $L(x,y)$.

We impose on $L$ the following two conditions:

(*) $L$ is of trace class.

(**) All finite diagonal minors $\det L_X$ are nonnegative.

Under these assumptions we have
$$
\sum_{X\in\Conff(\X)}\det L_X=\det(1+L)<\infty.
$$
We agree that $\det L_\varnothing=1$. Hence, the sum above is always
strictly positive.

Now we form a point process on $\X$ living on the finite
multiplicity free configurations $X\in\Conff(\X)$ with the probabilities
given by
$$
\operatorname{Prob}(X)=(\det (1+L))^{-1} \det L_X,
\qquad X\in\Conff(\X).
\tag 5.5
$$

It is convenient to have a name for the processes obtained in this way;
let us call them the {\it L--ensembles}.

\proclaim{Proposition 5.3} Let $L$ satisfy the conditions (*) and (**) above.
Then  the associated L--ensemble is a determinantal process with the
correlation  kernel $K=L(1+L)^{-1}$.
\endproclaim

\demo{Proof} See \cite{DVJ, Exercise 5.4.7}, \cite{BO2, Proposition 2.1}, 
\cite{BOO, Appendix}. \qed
\enddemo

The condition (*) can be slightly relaxed, see \cite{BOO}. The condition (**)
holds, for instance, when $L$ is Hermitian nonnegative. However, this
is by no means necessary, see \S5(f) below.

The relation between $L$ and $K$ can also be written in the form
$$
1-K=(1+L)^{-1}.
$$

\example{Remark 5.4} Assume that $K$ is a finite--dimensional orthogonal
projection  operator in $\ell^2(\X)$ (for instance, $K(x,y)=K^\cd(x,y)$ as in
\S5(d)). One  can prove that there exists a determinantal point process $\Cal
P$ for which $K$  serves as the correlation kernel. $\Cal P$ is not an
$L$--ensemble, because $1- K$ is not invertible (except $K=0$). However, $\Cal
P$ can be approximated by  certain $L$--ensembles. To see this, replace $K$ by
$K_\ep=\ep K$, where  $0<\ep<1$. The matrices $L_\ep=(1-K_\ep)^{-1}-1$ satisfy
both (*) and (**). The  process $\Cal P$ arises in the limit of the
L--ensembles associated with the  matrices $L_\ep$ as $\ep\nearrow 1$. One can
check that the probabilities $$
\lim_{\ep\nearrow1} \det(1+L_\ep)^{-1}\det (L_\ep)_X\,,
\qquad X\in\Conff(\X),
$$
are correctly defined.
\endexample

For a special class of matrices $L$ there exists a complex analytic problem
the solution of which yields the resolvent matrix $K$.

We will follow the exposition of \cite{B3}.

Let $\x$ be a discrete locally finite subset of $\C$. We call an operator
$L$ acting in $\ell^2(\x)$ {\it integrable} if its matrix has the form
$$
L(x,x')=\cases \dfrac {\sum_{j=1}^M f_j(x)g_j(x')}{x-x'},&\quad x\ne x',\\
0,&\quad x=x',
\endcases
\tag 5.6
$$
for some functions $f_j,\ g_j$ on $\x$, $j=1,\dots,M$, satisfying the relation
$$
\sum_{j=1}^M f_j(x)g_j(x)=0,\quad x\in\x.
\tag 5.7
$$
We will assume that $f_j,g_j\in\ell^2(\x)$ for all $j$.

Set
$$
f=(f_1,\dots,f_M)^t,\quad g=(g_1,\dots,g_M)^t.
$$
Then \tht{5.7} can be rewritten as $g^t(x)f(x)=0$.
We will also assume that the operator
$$
(Th)(x)=\sum_{x'\in\x,\, x'\ne x}\frac{h(x')}{x-x'}
\tag 5.8
$$
is a bounded operator in $\ell^2(\x)$. For example, this holds for $\x=\Z+c$
for any $c\in\C$.\footnote{Indeed, then $T$ is a Toeplitz operator with the
symbol $\sum\limits_{n\ne 0}\dfrac{u^n}n\in L^{\infty}(S^1)$.}  Under these
assumptions, it is easy to see that $L$ is a bounded operator in $\ell^2(\x)$.

Now we introduce the complex analytic object.

Let $w$ be a map from $\x$ to $\operatorname{Mat}(k,\C)$, $k$ is a fixed
integer.

We say that a matrix function $m:\C\setminus\x\to
\operatorname{Mat}(k,\C)$  with simple poles at the points $x\in\x$ is a
solution of the {\it discrete Riemann--Hilbert problem}\footnote{DRHP, for
short} $(\x,w)$ if the following conditions are satisfied
$$ \align
&\bullet\quad  m(\zeta)\text{ is analytic in }\Bbb C\setminus \x,\\
&\bullet\quad
\r m(\ze)=\lim_{\zeta\to x}\left(m(\zeta)w(x)\right),\quad x\in\x,\\
&\bullet\quad m(\zeta)\to I \text{ as } \zeta\to \infty.
\endalign
$$
Here $I$ is the $k\times k$ identity matrix. The matrix $w(x)$ is called the
{\it jump matrix}.

If the set $\x$ is infinite, the last condition must be made more precise.
Indeed, a function with poles accumulating at infinity cannot have asymptotics
at infinity. One way to make this condition precise is to require
the uniform asymptotics on a sequence of expanding contours, for example, on a
sequence of circles $|\ze|=a_k$, $a_k\to+\infty$.

In order to guarantee the uniqueness of solutions
of the DRHPs considered below, we always assume that there exists a
sequence of expanding contours such that the distance from these contours to
the set $\x$ is bounded from zero, and we will require a solution $m(\zeta)$
to uniformly converge to $I$ on these contours.

The setting of the DRHP above is very similar to the pure soliton case in
the inverse scattering method, see \cite{BC}, \cite{BDT}, \cite{NMPZ,
Ch. III}.

\proclaim{Proposition 5.5 \cite{B3, Proposition 4.3}}  Let $L$ be an
integrable operator as described above such that the operator $(1+L)$ is
invertible, and $m(\zeta)$ be a solution of the DRHP $(\x, w)$ with
$$
w(x)=-f(x)g(x)^t\in \operatorname{Mat}(M,\C).
$$
Then the matrix $K=L(1+L)^{-1}$ has the form
$$
K(x,x')=
\cases\dfrac {G^t(x')F(x)}{x-x'},\quad &x\ne x',\\
 G^t(x)\lim\limits_{\zeta\to x}\left(m'(\zeta)
 \,f(x)\right), &x=x',
\endcases
$$
where $m'(\ze)=\dfrac{dm(\ze)}{d\ze}$, and
$$
F(x)=\lim_{\zeta\to x}\left(m(\zeta)\, f(x)\right), \qquad
G(x)=\lim_{\zeta\to x}\left((m^{t}(\zeta))^{-1}\, g(x)\right).
$$
\endproclaim
\demo{Comments} 1) The continuous analog of this result was originally proved
in \cite{IIKS}, see also \cite{De} and \cite{KBI}.

2) It can be proved that the solution of the DRHP stated in Proposition 5.5
exists and is unique, see \cite{B3, (4.9)} for the existence and \cite{B3, Lemma 
4.7} for the uniqueness.

3) The requirement of matrix $L$'s vanishing on the diagonal can be
substantially weakened, see \cite{B3, Remark 4.2}. A statement similar to
Proposition 5.5 can be proved if the diagonal elements of $L$ are bounded from
$-1$.

4) Proposition 5.5 holds {\it without} the assumptions \tht{*}, \tht{**}
stated in the beginning of this subsection.

5) If the operator $L$ is bounded, has the form \tht{5.6}, but the functions
$f_j$ and $g_j$ are not in $\ell^2(\x)$, then it may happen that the operator
$K=L/(1+L)$ is well--defined while the corresponding DRHP fails to have a
solution.

\enddemo

\subhead (f) Special matrices $L$ \endsubhead
Let $\X$ be a discrete space with a fixed splitting into the union of two
disjoint subsets,
$$
\X=\X_{I}\sqcup \X_{II}.
$$
The splitting induces an orthogonal decomposition of $\ell^2(\X)$,
$$
\ell^2(\X)=\ell^2(\X_{I})\oplus\ell^2(\X_{II}).
$$
According to this decomposition we will write operators in $\ell^2(\X)$ (or
matrices of the format $\X\times\X$) in the block form. For instance,
$$
L=\bmatrix L_{I,I}& L_{I,II} \\ L_{II,I} & L_{II,II}\endbmatrix,
$$
where $L_{I,I}$ acts from $\ell^2(\X_I)$ to $\ell^2(\X_I)$,
$L_{I,II}$ acts from
$\ell^2(\X_{II})$ to $\ell^2(\X_I)$, etc.

We are interested in the matrices $L$ of the following special form:
$$
L=\bmatrix 0 & A \\ -A^* & 0 \\ \endbmatrix,
\tag 5.9
$$
where $A$ is an operator from $\ell^2(\X_{II})$ to $\ell^2(\X_I)$ and
$A^*$ is
the adjoint operator.

For such $L$, the condition (**) of \S{5(e)} is satisfied, while the
condition
(*) is equivalent to saying that $A$ is of trace class. It can be shown that
the construction of \S{5(e)} holds even if $A$ is a Hilbert--Schmidt operator,
see \cite{BOO, Appendix}.

Note that the matrices of the form \tht{5.9} are not Hermitian symmetric but
J--symmetric. I.e., the corresponding operator is Hermitian with
respect to the indefinite inner product on the space $\ell^2(\X)$
defined by the matrix
$J=\bmatrix 1 & 0\\ 0 &-1\endbmatrix$. It follows that the matrices
$K=L(1+L)^{- 1}$ are J--symmetric, too. This provides a class of
determinantal processes
whose correlation kernels are not Hermitian symmetric.

Now let us look at an even more special situation. Assume that $\x$ is a
locally finite subset of $\C$ such that the operator $T$ defined by \tht{5.8}
is bounded.

Let $h_{I}(\,\cdot\,)$, $h_{II}(\,\cdot\,)$ be two functions defined on
$\x_I$ and $\x_{II}$, respectively. We assume that $h_{I}\in\ell^2(\x_I)$,
$h_{II}\in\ell^2(\x_{II})$. (The functions $h_I$, $h_{II}$ should not
be confused with the constants $h_n$ attached to orthogonal
polynomials, see \S5(d).) 

Set
$$
L=\bmatrix 0 & A\\ -A^* & 0\endbmatrix,
\qquad \text{where}\quad A(x,y)=\frac{h_I(x)h_{II}(y)}{x-y}\,.
\tag 5.10
$$

The matrix $A$ is well defined, because $x$ and $y$ range over disjoint
subsets  $\X_I$ and $\X_{II}$ of $\x$.

As is explained in \cite{B3, \S6}, such $L$ is an integrable operator in
the sense of \S5(e) with $M=2$. Let us assume that the functions $h_I$ and
$h_{II}$ are real--valued. Then we have $L^*=-L$, and $-1$ cannot
belong to the
spectrum of L, that is, $(1+L)$ is invertible. Thus, the DRHP of Proposition
5.5 has a unique solution.

Let us introduce a special notation for this solution $m(\ze)$.
We define four meromorphic functions $R_{I}$, $S_{I}$,
$R_{II}$,  $S_{II}$ by the relation
$$
m=\bmatrix m_{11}& m_{12}\\ m_{21}& m_{22}\endbmatrix=\bmatrix R_{I}&
-S_{II}\\
-S_{I}& R_{II}\endbmatrix.
$$

Then the DRHP of Proposition 5.5 for our special $L$ given by \tht{5.10} can
be restated as follows, see \cite{B3, \S6}:

\noindent$\bullet$  matrix elements $m_{11}=R_I$ and $m_{21}=-S_I$ are
holomorphic in  $\C\setminus\x_{II}$;

\noindent$\bullet$  matrix elements $m_{12}=-S_{II}$ and $m_{22}=R_{II}$ are
holomorphic in  $\C\setminus\x_{I}$;

\noindent$\bullet$  $R_I$ and $S_I$ have simple poles at the points of
$\x_{II}$, and for $x\in\x_{II}$
$$
\gathered
\r R_I(\ze)=h_{II}^2(x)S_{II}(x),\\
\r S_{I}(\ze)=h_{II}^2(x)R_{II}(x);
\endgathered
$$
\noindent$\bullet$  $R_{II}$ and $S_{II}$ have simple poles at the points of
$\x_{I}$, and for $x\in\x_{I}$
$$
\gathered
\r R_{II}(\ze)=h_{I}^2(x)S_{I}(x),\\
\r S_{II}(\ze)=h_{I}^2(x)R_{I}(x);
\endgathered
$$
\noindent$\bullet$ $R_{I},R_{II}\to 1$, $S_{I},S_{II}\to 0$ as
$\ze\to\infty$.
\smallskip

As before, the last condition is understood as the uniform convergence on
a sequence of expanding contours such that the distance from these contours to
the set $\x$ is bounded from zero.

It can be proved that these conditions imply the relations
$$
\gather
R_{I}(\ze)=1-\sum_{y\in\x_{II}}\frac
{h_{II}^2(y)S_{II}(y)}{y-\ze}\,,\quad S_{I}(\ze)=-\sum_{y\in\x_{II}}\frac
{h_{II}^2(y)R_{II}(y)}{y-\ze} \,,
\tag 5.11\\
R_{II}(\ze)=1-\sum_{y\in\x_{I}}\frac
{h_{I}^2(y)S_{I}(y)}{y-\ze}\,,\quad S_{II}(\ze)=-\sum_{y\in\x_{I}}\frac
{h_{I}^2(y)R_{I}(y)}{y-\ze}\,.
\tag 5.12
\endgather
$$
The inverse implication also holds if we know that the functions
$R_{I},S_I,R_{II},S_{II}$ have the needed asymptotics as $\ze\to\infty$.

The next statement is a direct corollary of Proposition 5.5.

\proclaim{Proposition 5.6 \cite{B3, Proposition 6.1}} Let
$$
m=\bmatrix R_{I}& -S_{II}\\
-S_{I}& R_{II}\endbmatrix.
$$
be a solution of the DRHP stated above, where $h_I\in\ell^2(\x_I)$ and
$h_{II}\in\ell^2(\x_{II})$ are real--valued. Then $(1+L)$ is invertible, and
the matrix of the operator $K=L/(1+L)$, with respect to the splitting
$\x=\x_I\sqcup\x_{II}$, has the form
$$
\gathered
K_{I,I}(x,y)={h_I(x)h_I(y)}
\,\frac {R_I(x)S_I(y)-S_I(x)R_I(y)}{x-y}\,,\\
K_{I,II}(x,y)={h_I(x)h_{II}(y)}
\,\frac {R_I(x)R_{II}(y)-S_I(x)S_{II}(y)}{x-y}\,,\\
K_{II,I}(x,y)=h_{II}(x)h_I(y)
\,\frac {R_{II}(x)R_I(y)-S_{II}(x)S_I(y)}{x-y}\,,\\
K_{II,II}(x,y)=h_{II}(x)h_{II}(y)
\,\frac {R_{II}(x)S_{II}(y)-S_{II}(x)R_{II}(y)}{x-y}\,,
\endgathered
$$
where the
indeterminacy on the diagonal $x=y$ is resolved by the L'Hospital rule:
$$
\gathered
K_{I,I}(x,x)=h_I^2(x)
\left((R_I)'(x)S_I(x)-(S_I)'(x)R_I(x)\right),
\\
K_{II,II}(x,x)=h_{II}^2(x)
\left((R_{II})'(x)S_{II}(x)-(S_{II})'(x)R_{II}(x)\right).
\endgathered
$$
\endproclaim

\subhead (g) Connection between discrete polynomial ensembles and
L--en\-sembles \endsubhead

Here we adopt the following assumptions:

$\bullet$ $\X=\X_I\sqcup\X_{II}$ is a finite or countably infinite subset of
$\R$ without limit points.

$\bullet$ The set $\X_{II}$ is finite, $|\X_{II}|=N$.

$\bullet$  $h_I$ is a nonnegative function on $\X_I$ such that
$$
\sum_{x\in\X_I}\frac{h^2_I(x)}{1+x^2}<\infty.
\tag 5.13
$$

$\bullet$ $h_{II}$ is a strictly positive function on $\X_{II}$.

To these data we associate a function $f$ on $\X$ as follows:
$$
f(x)=\cases \dfrac{h^2_I(x)}{\prod\limits_{y\in\X_{II}}(x-y)^2}\,,
& x\in\X_I\,, \\
\dfrac1{h^2_{II}(x)\prod\limits\Sb y\in\X_{II}\\ y\ne x\endSb(x-y)^2}\,,
& x\in\X_{II}\,.
\endcases
\tag 5.14
$$

We note that $f$ is nonnegative on $\X$, strictly positive on $\X_{II}$, and
its $(2N-2)$nd moment is finite,
$$
\sum_{x\in\X}x^{2N-2}f(x)<\infty.
$$
Conversely, given $f$ with such properties, we can define the functions
$h_I$ and $h_{II}$ by inverting \tht{5.14}, and then the condition
\tht{5.13} will be satisfied.

Since the function $f$ satisfies the two basic assumptions for a
weight function
stated in \S5(d) (the moment of order $2N-2$ is finite, and $f$ is strictly
positive on an $N$--point subset), we can attach to it a discrete polynomial
ensemble.

On the other hand, let us define a matrix $L$ using \tht{5.10}.
By virtue of \tht{5.13}, all the columns of $A$ are
vectors from $\ell^2(\X_I)$. Since the total number of columns in $A$ is finite,
the trace class condition for $A$ holds for trivial reasons. According to
\S{5(f)}, such a matrix $L$ defines a determinantal point process.

\proclaim{Proposition 5.7} Under the above assumptions, the orthogonal
polynomial ensemble with the weight function $f$ and the L-ensemble
associated with the  matrix \tht{5.10} are connected by the involution
$\triangle$ corresponding to  $Z=\X_{II}$.
\endproclaim

\demo{Proof} We will prove that for any balanced configuration $X$,
the probability of $X$ in the the L--ensemble is equal to the probability of
$X^\triangle$ in the orthogonal polynomial ensemble.  We have to compare two
expressions, \tht{5.3} and \tht{5.5}, which both involve a normalizing
constant.  Since we know that we are dealing with probability measures, we may
ignore constant factors.

Let $X$ be a finite balanced configuration with no multiplicities. Write it as
$A\sqcup B$, where $A=X\cap\x_{I}=\{a_1,\dots,a_d\}$,
$B=X\cap\x_{II}=\{b_1,\dots,b_d\}$. In this notation, the probability of $X$ in
the  L--ensemble is equal, up to a constant factor, to
$$
\det L_{A\sqcup B}=\prod_{i=1}^d h^2_I(a_i)h^2_{II}(b_i)\cdot
{\det}^2 \left[\frac1{a_i-b_j}\right]_{1\le i,j\le d}\,.
\tag 5.15
$$

For arbitrary finite configurations
$C=\{c_1,\dots,c_m\}$ and $D=\{d_1,\dots,d_n\}$ we will
abbreviate
$$
V^2(C)=\prod_{1\le i<j\le m}(c_i-c_j)^2, \qquad
V^2(C;D)=\prod_{i=1}^m\prod_{j=1}^n(c_i-d_j)^2.
$$

By the well--known formula for Cauchy's determinant,
$$
{\det}^2\left[\frac1{a_i-b_j}\right]=\frac{V^2(A)V^2(B)}{V^2(A;B)}\,,
$$
the expression \tht{5.15} is equal to
$$
\prod_{a\in A}h^2_I(a)\cdot\prod_{b\in B}h^2_{II}(b)\cdot
\frac{V^2(A)V^2(B)}{V^2(A;B)}\,.
\tag 5.16
$$

On the other hand, $X^\triangle=A\sqcup\bar B$, where
$\bar B=\X_{II}\setminus B$. The probability of $X^\triangle$ in the orthogonal
polynomial ensemble is equal, up to a constant factor, to
$$
\prod_{x\in X^\triangle}f(x)\cdot V^2(X^\triangle)=
\prod_{a\in A}f(a)\cdot
\prod_{\bar b\in\bar B}f(\bar b)\cdot
V^2(A\sqcup\bar B).
\tag 5.17
$$
Let us transform this expression. We have
$$
\prod_{\bar b\in\bar B}f(\bar b)=
\const\cdot\prod_{b\in B}\frac1{f(b)}\,, \quad
\const=\prod_{x\in\X_{II}}f(x).
$$
Then
$$
\gather
V^2(A\sqcup\bar B)=V^2(A)V^2(\bar B)V^2(A;\bar B)\\
=\frac{V^2(A)V^2(B)}{V^2(A;B)}\cdot
\frac{V^2(A;\bar B)V^2(A;B)}{V^4(B)V^2(B;\bar B)}\cdot
V^2(B)V^2(\bar B)V^2(B;\bar B)\\
=\const\cdot\frac{V^2(A)V^2(B)}{V^2(A;B)}\cdot
\frac{\prod\limits_{a\in A}\,\prod\limits_{y\in\X_{II}}(a-y)^2}
{\prod\limits_{b\in B}\,\prod\limits_{y\in\X_{II}\setminus\{b\}}(b-y)^2}\,,
\endgather
$$
where
$$
\const=V^2(B)V^2(\bar B)V^2(B;\bar B)=V^2(\X_{II}).
$$
It follows that \tht{5.17} is equal, up to a constant factor, to
$$
\prod_{a\in A}\left(f(a)\prod_{y\in\X_{II}}(a-y)^2\right)\cdot
\prod_{b\in B}\left(f(b)\prod_{y\in\X_{II}\setminus\{b\}}(b-y)^2\right)^{-
1}\cdot
\frac{V^2(A)V^2(B)}{V^2(A;B)}\,.
$$
By virtue of the connection between $f$ and $\{h_I,h_{II}\}$, see
\tht{5.14}, this is equal to  \tht{5.16}.
\qed
\enddemo

\subhead (h) Connection between two correlation kernels \endsubhead
We keep the assumptions  of \S{5(g)}. In particular, $f$
is related to $h_I$ and $h_{II}$ by \tht{5.14}.

Recall that if we assume that $h_I\in\ell^2(\x_I)$ and the operator $T$
(see \tht{5.8}) is bounded then the DRHP of \S5(f) has a unique solution which
defines the meromorphic functions $R_I,S_I,R_{II},S_{II}$. Also,
$h_I\in\ell^2(\x_I)$ implies that the function $f$ defined through \tht{5.14}
has a finite $(2N)$th moment, and, hence, we can define monic, orthogonal with
respect to the weight $f$ polynomials $\{p_0,p_1,\dots\}$ at least up to the
$N$th one, see \S5(d). Note that the condition $h_I\in\ell^2(\x_I)$
is stronger than \tht{5.13}.

\proclaim{Proposition 5.8} Under the assumptions of \S5(g), if
$h_I\in\ell^2(\x_I)$ and the operator $T$ is bounded, then we have
$$
\gathered
R_{I}(\ze)=\frac{p_N(\ze)}{\prod\limits_{y\in\x_{II}}(\ze-y)}\,, \qquad
S_{I}(\ze)=\frac{p_{N-1}(\ze)}{h_{N-1}\prod\limits_{y\in\x_{II}}(\ze-y)}\,.
\endgathered
\tag 5.18
$$
\endproclaim
\demo{Proof} Denote the right--hand sides of \tht{5.18} by $\wt R_I$ and $\wt
S_I$, respectively, and define
$$
\wt R_{II}(\ze)=1-\sum_{y\in\x_{I}}\frac
{h_{I}^2(y)\wt S_{I}(y)}{y-\ze}\,,\qquad \wt
S_{II}(\ze)=-\sum_{y\in\x_{I}}\frac {h_{I}^2(y)\wt R_{I}(y)}{y-\ze}\,,
$$
cf. \tht{5.12}. We will show that the matrix
$$
\wt m=\bmatrix \wt R_{I}& -\wt S_{II}\\
-\wt S_{I}& \wt R_{II}\endbmatrix
$$
solves the DRHP of \S5(f). By uniqueness of the solution we will conclude that
$m=\wt m$.

The condition $h_I\in\ell^2(\x_I)$ guarantees that the formulas above define
meromorphic functions $\wt R_{II}$, $\wt S_{II}$ with needed asymptotics and
location of poles. Thus, we  only need to check the relations involving
residues at the poles. The equalities $$
\gathered
\r \wt R_{II}(\ze)=h_{I}^2(x)\wt S_{I}(x),\\
\r \wt S_{II}(\ze)=h_{I}^2(x)\wt R_{I}(x),
\endgathered
$$
are obviously satisfied.

The relation $\r \wt S_{I}(\ze)=h_{II}^2(x)\wt R_{II}(x)$ is equivalent to the
equality
$$
-\frac{p_{N-1}(x)}{h_{N-1}\prod\limits_{y\in\x_{II},\, y\ne x}(x-y)}=
-h_{II}^2(x)\left(1-\sum_{y\in\x_{I}}\frac
{h_{I}^2(y)\wt S_{I}(y)}{y-x}\right)\, \tag 5.19
$$
which can be rewritten as ($x\in\x_{II}$)
$$
-{p_{N-1}(x)f(x)\prod\limits_{t\in\x_{II},\, t\ne x}(x-t)}=-h_{N-
1}+\sum_{y\in\x_{I}}
{p_{N-1}(y)f(y)\prod_{t\in\x_{II},\, t\ne x}(y-t)}.
$$
But this is the relation
$$
\left\langle p_{N-1}(y),\prod_{t\in\x_{II},\, t\ne x}(y-t)\right\rangle =h_{N-1}
$$
which directly follows from the definition of the orthogonal polynomials.

The relation $\r \wt R_{I}(\ze)=h_{II}^2(x)\wt S_{II}(x)$
is equivalent to the equality
$$
\frac{p_{N}(x)}{\prod\limits_{y\in\x_{II},\, y\ne x}(x-y)}=
-h_{II}^2(x)\sum_{y\in\x_{I}}\frac {h_{I}^2(y)\wt R_{I}(y)}{y-x}\,.
\tag 5.20
$$
which is just the orthogonality relation
$$
\left\langle p_{N}(y),\prod_{t\in\x_{II},\, t\ne x}(y-t)\right\rangle =0.\qed
$$
\enddemo

\proclaim{Corollary 5.9}
Under the assumptions of \S5(g), if $h_I\in\ell^2(\x)$ and the operator $T$ is
bounded, then for any $x\in \x_{II}$ we have
$$
R_{II}(x)=\frac{p_{N-1}(x)}{h_{N-1}\,h_{II}^2(x)\prod\limits_{y\in\x_{II},\,
y\ne x}(x-y)}\,,
\quad
S_{II}(x)=\frac{p_{N}(x)}{h_{II}^2(x)\prod\limits_{y\in\x_{II},\,
y\ne x}(x-y)}\,. \tag 5.21
$$
\endproclaim
\demo{Proof} Follows from the relations
$$
\r  R_{I}(\ze)=h_{II}^2(x) S_{II}(x),\qquad
\r  S_{I}(\ze)=h_{II}^2(x) R_{II}(x),
$$
and \tht{5.18}. \qed
\enddemo

For the next statement we drop the assumption $h_I\in\ell^2(\x_I)$ and use the
weaker assumption \tht{5.13} instead.

\proclaim{Theorem 5.10} Under the assumptions of {\rm\S5(g)}, let $L$
be given by \tht{5.10}, $K=L(1+L)^{-1}$, and  $K^\cd$ be the $N$th normalized
Christoffel--Darboux kernel for the weight function $f$, as  defined in
{\rm\S5(d)}. Introduce the following function on $\X$ taking values in
$\{\pm1\}$:
$$
\ep(x)=\cases \sgn\left(\prod\limits_{y\in\X_{II}}(x-y)\right), & x\in\X_I\\
\sgn \left(\prod\limits_{y\in\X_{II}\setminus\{x\}}(x-y)\right),
& x\in\X_{II}.
\endcases
$$

Then we have
$$
K(x,y)=\ep(x)(K^\cd)^\triangle(x,y)\ep(y),
\tag 5.22
$$
where the operation $(\,\cdot\,)^\triangle$ is defined by \tht{5.1} with
$Z=\X_{II}$, $\bar Z=\X_I$.
\endproclaim

Before proceeding to the proof let us make a couple of comments.

\demo{Comments} 1) By Proposition 5.3 the kernel $K$ describes the correlation
functions of  the L--ensemble in question. On the other hand, the same
correlation functions  are also expressed in terms of the kernel
$(K^\cd)^\triangle$, see Proposition 5.7 and \S5(c). This does not mean that
both kernels must coincide, because the correlation kernel of a  determinantal
point process is not defined uniquely, see \S5(b). And indeed, we see that the
kernels turn out to be conjugated by a nontrivial diagonal matrix. Note that
conjugating by a diagonal matrix is the only possible ``generic''
transformation, because this is the only operation on the ``generic'' matrix
which preserves all diagonal minors. In our situation, both kernels are real
and  possess a symmetry property (J--symmetry), so that it is not surprising
that the  diagonal entries of this diagonal matrix are equal to $\pm1$.
However, the exact  values of these $\pm1$'s, as given in the theorem, are not
evident.

2) The theorem makes it possible to calculate the kernel $K=L(1+L)^{-1}$
provided that we know the orthogonal polynomials with the weight function $f$.
Both kernels, $K$ and $(K^\cd)^\triangle$, are suitable to describing the
correlation functions of the L--ensembles. However, from the
computations that follow in subsequent sections we will see that, for the
particular L--ensemble we are interested in, the former kernel survives in a
scaling limit procedure (see \S10) while the latter kernel does not.
\enddemo

\demo{Proof of Theorem 5.10} Let us assume that the stronger condition
$h_I\in\ell^2(\x_I)$ is satisfied. Then the normalized Christoffel--Darboux
kernel can be written in the form \tht{5.4},  which makes it possible to
express $(K^\cd)^\triangle$ as follows:
$$
\aligned
&(K^\cd)^\triangle_{I,I}(x,y)=\ep(x)\ep(y)\,\frac 
{\sqrt{f(x)f(y)}}{h_{N-1}}
\,\frac {p_N(x)p_{N-1}(y)-p_{N-1}(x)p_N(y)}{x-y}\,,\\
&(K^\cd)^\triangle_{I,II}(x,y)=\ep(x)\ep(y)\frac 
{\sqrt{f(x)f(y)}}{h_{N-1}}\,\frac
{p_N(x) p_{N-1}(y)-p_{N-1}(x)p_{N}(y)}{x-y}\,,\\
&(K^\cd)^\triangle_{II,I}(x,y)=\ep(x)\ep(y)\frac 
{\sqrt{f(x)f(y)}}{h_{N-1}}\,\frac
{p_{N-1}(x) p_{N}(y)-p_N(x) p_{N-1}(y)}{x-y}\,,\\
&(K^\cd)^\triangle_{II,II}(x,y)\\
&\qquad \qquad=\delta(x-y)+\ep(x)\ep(y)\,\frac
{\sqrt{f(x)f(y)}}{h_{N-1}} \frac {p_{N-1}(x)p_{N}(y)-p_N(x)p_{N-1}(y)}{x-y}\,.
\endaligned
$$
Here $f$ is the function defined in \tht{5.14}.

Assuming that the operator $T$ of \tht{5.8} is bounded, we see that
\tht{5.18}, \tht{5.21}, and Proposition 5.6 directly imply the claim of
the theorem everywhere except for the diagonal set $(x,x)\in \x\times\x$. But
on this diagonal set it is immediately  seen that both sides of the equality
$$
K(x,x)=(K^\cd)^\triangle(x,x)
$$
represent the probability that the corresponding L--ensemble has a particle at
the point $x$, see Comment 1 above.

Thus, we have verified \tht{5.22} assuming that
$h_I\in\ell^2(\x_I)$ and the boundedness of $T$. Now we will show how to get
rid of these extra conditions.

Let $H_I$ denote the set of all nonnegative functions $h_I(x)$
satisfying \tht{5.13}. We equip $H_I$ with the weakest topology for
which all the evaluations $h_I\mapsto h_I(x)$ ($x\in\x_I$) and the
sum in \tht{5.13} are continuous.  Fix $h_{II}$ and let $h_I$ vary
over $H_I$. We claim that for any fixed $x,y\in\x$, both values
$K^{\operatorname{CD}}(x,y)$ and $K(x,y)$ depend continuously on $h_I\in H_I$.

By the definition of the topology in $H_I$, the moments
$$
c_j=\sum_{x\in\x}x^j f(x), \qquad j=0,\dots,2N-2,
$$
depend continuously on $h_I$. Let $G=[c_{i+j}]_{i,j=0}^{N-1}$ be the
Gram matrix of the vectors $1,x,\dots,x^{N-1}\in L^2(\x,f\mu)$. The
Christoffel--Darboux kernel can be expressed through the moments as
follows:
$$
\sum_{i,j=0}^{N-1}(G^{-1})_{ij}x^iy^j.
$$
(This can be derived from the classical determinantal expression for
the orthogonal polynomials, see \cite{Er, Vol. 2, 10.3(4)}, or by making use of
a simple direct argument, see, e.g., \cite{B4}. )
Hence
$$
K^{\operatorname{CD}}(x,y)=\sqrt{f(x)f(y)}\sum_{i,j=0}^{N-1}(G^{-1})_{ij}x^iy^j.
$$
Clearly, this expression is a continuous function of $h_I$.

On the other hand, again by the definition of the topology, the
columns of the matrix $A(x,y)$ (see \tht{5.10}), viewed as vectors in $\ell^2(\x_I)$,
are continuous in $h_I$. It follows that the operators $L$ and
$K=L(1+L)^{-1}$ are continuous in $h_I$ with respect to uniform
operator topology. Hence $K(x,y)$ is continuous, too.

Thus, we have proved that both sides of the required equality
\tht{5.22} are continuous in $h_I\in H_I$.

Finally, let $H^0_I\subset H_I$ be the subset of those functions which
have finite support. Note that \tht{5.22} holds for any $h_I\in H_I^0$.
Indeed, \tht{5.22} does not change if we replace $\x_I$ by
$\operatorname{supp}h_I$, and for a finite $\x_I$ the extra conditions are
obviously satisfied. Since $H^0_I$ is dense in $H_I$, we conclude that
\tht{5.22} holds for any $h_I\in H_I$. \qed
\enddemo

\head 6. $\Cal P^{(N)}$ and $\tP^{(N)}$ as determinantal point processes
\endhead

Recall that in \S4 we defined two discrete point processes
$\tP^{(N)}$ and $\Cal P^{(N)}$. These processes live on the lattice
$\x^{(N)}$ and depend on four parameters $(\zw)\in\Dadm$. As was
mentioned in \S5(c), $\tP^{(N)}$ and $\Cal P^{(N)}$ are related by the
complementation principle: $\tP^{(N)}=(\Cal P^{(N)})^\triangle$, where
the special set $Z$ is equal to $\xin$. In this section we will show
that $\tP^{(N)}$ is a discrete polynomial ensemble (as defined in
\S5(d)) and $\Cal P^{(N)}$ is an L--ensemble (as defined in \S5(e)).

Consider the following weight function on the lattice $\x^{(N)}$:
$$
f(x)=\frac 1{\Gamma\left(z-x+\frac{N+1}2\right)\Gamma\left(z'-
x+\frac{N+1}2\right)\Gamma
\left(w+x+\frac{N+1}2\right)\Gamma\left(w'+x+\frac{N+1}2\right)}\,.
\tag 6.1
$$
Here we assume that $(\zw)\in\Dadm$. The expression \tht{6.1} comes
from the expression for $P'_N(\la\mid\zw)$, see \S3. Namely, in the
notation of \S4, we have
$$
P'_N(\la\mid\zw)=\Dim_N^2(\la)\cdot\prod_{x\in\Cal L(\la)}f(x).
$$

\proclaim{Proposition 6.1} Let $(\zw)\in\Dadm$. The function $f(x)$
is nonnegative on $\x^{(N)}$ and satisfies the assumptions (*) and
(**) of \S5(d).
\endproclaim

\demo{Proof} The nonnegativity of $f(x)$ follows from Proposition 3.5
and the definition of $\Dadm$, see Definition 3.4.

The condition (*) of \S5(d) says that
$$
\sum_{x\in\x^{(N)}}x^{2N-2}f(x)<\infty.
$$
This follows from the estimate
$$
f(x)\le \const\cdot(1+|x|)^{-(z+z'+w+w'+2N)}
$$
and the fact that $z+z'+w+w'>-1$ for $(\zw)\in\Dadm$. As for the
estimate above, it readily follows from the asymptotics of the gamma
function, see \cite{Ol3, \tht{7.6}}.

Finally, the condition (**) of \S5(d) says that $f(x)$ does not
vanish at least on $N$ distinct points. This follows from the fact
that $P'_N(\la\mid\zw)$ does not vanish identically, see Proposition
3.5. \qed
\enddemo

Note that if $(z,z')\in\Zd$ then $f(x)$ vanishes for positive large $x$. 
Similarly, if $(w,w')\in\Zd$ then $f(x)$ vanishes for negative $x$
such that $|x|$ is large enough.

\proclaim{Corollary 6.2} Let $(\zw)\in\Dadm$. Then $\tP^{(N)}$ is a
discrete polynomial ensemble with the weight function $f(x)$ given by
\tht{6.1}. That is, for any $N$-point configuration
$X=\{x_1,\dots,x_N\}$
$$
\tP^{(N)}(X)=\const\cdot\prod_{i=1}^Nf(x_i)\cdot\prod_{1\le i<j\le
N}(x_i-x_j)^2.
$$
\endproclaim

\demo{Proof} Indeed, by Proposition 6.1, the assumptions of \S5(d)
are satisfied. The formula above is a direct corollary of the definition of
$P_N(\la)$, see \S3, and the formula for $\operatorname{Dim}_N\la$, see
\S1(d).\qed
\enddemo

Now let turn to the process $\Cal P^{(N)}$. We will apply the
formalism of \S5(g) for $\x=\x^{(N)}$, $\x_{I}=\xout$,
$\x_{II}=\xin$. Recall that \S5(g) relies on four assumptions. The first three 
of them hold
for any quadruple $(\zw)\in\Dadm$. As for the fourth assumption, it
is equivalent to the strict positivity of the weight
function $f(x)$ on $\xin$, see \S5(g). It may happen that for some admissible
quadruples $(\zw)$ this requirement is violated:
$f(x)$ vanishes at certain points of $\xin$. Specifically, this
happens whenever $(z,z')$ or $(w,w')$ belongs to $\Zdm$ with $m<0$.
For this reason we have to impose an additional restriction on the
parameters.

\example{Definition 6.3} Let $\Dadm'$ denote the subset of $\Dadm$
formed by the quadruples $(\zw)\in\Dadm$ such that neither $(z,z')$ nor
$(w,w')$ belongs to $\Zdm$ with $m<0$.
\endexample

In the rest of the section we assume that $(\zw)\in\Dadm'$, so that
$f(x)>0$ for any $x\in\xin$. Note that, in terms of signatures $\la$,
this condition means that $P_N(\la\mid\zw)$ does not vanish at
$\la=(0,\dots,0)$.

\example{Remark 6.4} Recall that in Remark 3.7 we introduced a
natural shift on the set $\Dadm$,
$$
(\zw)\mapsto(z+k,z'+k,w-k,w'-k), \quad k\in\Z.
$$
This shift of the parameters is equivalent to the shift of all
configurations of $\tP^{(N)}$ by $k$. The definition of $\Dadm$
implies that any quadruple from $\Dadm\setminus\Dadm'$ can be moved
into $\Dadm'$ by an appropriate shift. Thus, for the study
of the process $\tP^{(N)}$, the restriction of the admissible
quadruples to $\Dadm'$ is not essential. This argument does not work
for the process $\Cal P^{(N)}$, because the shift above
does not preserve the splitting $\x=\xout\sqcup\xin$. However, as
will be shown later (see the proof of Theorem 10.1), the effect of the shift is 
negligible in the limit transition as $N\to\infty$.
\endexample

Let us define functions $\psiin^{(N)}$ and $\psiout^{(N)}$ on $\xin$
and $\xout$, respectively, by the formulas
$$
\psiin^{(N)}(x)=\Gamma\left[\matrix -x+z+\frac{N+1}2,-
x+z'+\frac{N+1}2,x+w+\frac{N+1}2,x+w'+\frac{N+1}2\\ -x+\frac{N+1}2,-
x+\frac{N+1}2,x+\frac{N+1}2,x+\frac{N+1}2\endmatrix\right]\,,
\tag 6.2
$$
$$
\psiout^{(N)}(x)=\frac{\left(\left(x-\frac{N-1}2\right)_N\right)^2}
{\Gamma(-x+z+\frac{N+1}2)\Gamma(-x+z'+\frac{N+1}2)
\Gamma(x+w+\frac{N+1}2)\Gamma(x+w'+\frac{N+1}2)}
\,,
\tag 6.3
$$
where we use the notation
$$
(a)_k=\frac{\Ga(a+k)}{\Ga(a)}=a(a+1)\cdots(a+k-1),\quad
\Gamma\left[\matrix a,\,b,\,\dots\\ c,\,d,\,\dots\endmatrix\right]=
\frac{\Gamma(a)\Gamma(b)\dots}{\Gamma(c)\Gamma(d)\dots}\,.
$$

Note that
$$
\gathered
\psiin^{(N)}(x)=\frac{1}{\left(\Gamma\left(-x+\frac{N+1}2\right)
\Ga\left(x+\frac{N+1}2\right)\right)^2 f(x)}\,,\\ \psiout^{(N)}(x)=\cases
\left(\dfrac{\Ga\left(x+\frac {N+1}2\right)}{\Ga\left(x-\frac
{N-1}2\right)}\right)^2 f(x),&x\ge\frac {N+1}2\,.\\
\left(\dfrac{\Ga\left(-x+\frac {N+1}2\right)}{\Ga\left(-x-\frac
{N-1}2\right)}\right)^2 f(x),&x\le-\frac {N+1}2\,.\endcases
\endgathered
$$

The function $\psiout^{(N)}$ is nonnegative and the function $\psiin^{(N)}$ is
strictly positive. Note also that both functions are invariant with
respect to the substitution
$$
(z,z')\longleftrightarrow (w,w'),\quad x\longleftrightarrow -x.
$$

Introduce a matrix $L^{(N)}$ of format $\x\times\x$ which in the block form
corresponding to the splitting $\x=\xout\sqcup\xin$ is given by
$$
L^{(N)}_{\xout\sqcup\xin}=\left[\matrix 0&\A^{(N)}\\
-(\A^{(N)})^*&0\endmatrix\right]\,,
$$
cf. \S5(f), where $\A^{(N)}$ is a matrix of format $\xout\times\xin$,
$$
\A^{(N)}(a,b)=\frac{\sqrt{\psiout^{(N)}(a)\psiin^{(N)}(b)}}{a-b}\,, \qquad
a\in\xout, \quad b\in\xin.
$$

\proclaim{Proposition 6.5} Let $(\zw)\in\Dadm'$. The process $\Cal
P^{(N)}$ introduced in \S4 is an L--ensemble with the matrix
$L^{(N)}$ introduced above. That is, for any finite configuration $X$
$$
\Cal P^{(N)}(X)=\frac{\det L_X^{(N)}}{\det(1+L^{(N)})}\,,
$$
where $L_X^{(N)}$ denotes the submatrix of $L^{(N)}$ of finite format $X\times
X$.  \endproclaim

\demo{Proof} This is a direct corollary of Propositions 5.7 and 6.1,
where for Proposition 5.7 we take
$$
\x_{I}=\xout, \quad \x_{II}=\xin,\quad h_I^2=\psiout^{(N)}, \quad
h_{II}^2=\psiin^{(N)}.
$$
The relations between the functions $\psiout^{(N)}$, $\psiin^{(N)}$, and $f$
above exactly coincide with \tht{5.14}. \qed
\enddemo

\head 7. The correlation kernel of the process $\tP^{(N)}$
\endhead

The goal of this section is to compute the normalized
Christoffel--Darboux kernel, see \S5(d), associated with the weight
function $f$ on the lattice $\x^{(N)}$ given by \tht{6.1}. According
to Proposition 5.1, this kernel can be taken as a correlation kernel
for the process $\tP^{(N)}$. We will denote this kernel by $\tK$.

We will show below that the orthogonal polynomials with the weight
$f$ can be expressed through the hypergeometric function of type
$(3,2)$. This is an analytic function in one complex variable $u$ defined
inside the unit circle by its Taylor series
$$
{}_3F_2\left[\matrix a,\,b,\,c\,\\
e,\,f\endmatrix\,\Bigr|\,u\right]=\sum_{k=0}^\infty
\frac{(a)_k(b)_k(c)_k}{k!(e)_k(f)_k}\, u^k.
$$
Here $a,b,c,e,f$ are complex parameters, $e,f\notin\{0,-1,-2,\dots\}$.
We will only need the value of this function at the point $u=1$. Then the
series above converges if $\Re (e+f-a-b-c)>0$. Moreover, the function
$$
\frac{1}{\Ga(e)\Ga(f)\Ga(e+f-a-b-c)}\,{}_3F_2\left[\matrix a,\,b,\,c\,\\
e,\,f\endmatrix\,\Bigr|\,1\right]
$$
can be analytically continued to an entire function in 5 complex variables
$a,b,c,e,f$. \footnote{We could not
find a proof of this fact in the literature. We give our own proof
in the end of the Appendix.} Note also that if
one of the parameters
$a,b,c$ is a nonpositive integer, say, $a\in\{0,-1,-2,\dots\}$, then
the series above has only finitely many nonzero terms. It is easy to
see that in such a case
$$
\frac{1}{\Ga(e)\Ga(f)}\,{}_3F_2\left[\matrix a,\,b,\,c\,\\
e,\,f\endmatrix\,\Bigr|\,1\right]
$$
is an entire function in $b,c,e,f$.

Let us return to the process $\tP^{(N)}$. Set $\s=z+z'+w+w'$.

\proclaim{Theorem 7.1} For any $(z,z',w,w')\in \Dadm$, the
normalized Christoffel--Darboux kernel $\tK$ is given by
$$
\tK(x,y)=\frac
1{h_{N-1}}\,\frac{\p_N(x)\p_{N-1}(y)-\p_{N-1}(x)\p_N(y)}
{x-y}\,\sqrt{f(x)f(y)}\,,
\tag 7.1
$$
where $x,y\in\X^{(N)}$,
$$
\gathered
\p_N(x)=\frac{\Gamma(x+w'+\frac{N+1}2)}{\Gamma(x+w'-\frac{N-1}2)}\,
{}_3F_2\left[\matrix -N,\,z+w',\,z'+w'\\ \s,\,x+w'-\frac{N-
1}2\endmatrix\Biggl|\,1\,\right],\\
\p_{N-1}(x)=\frac{\Gamma(x+w'+\frac{N+1}2)}{\Gamma(x+w'-\frac{N-1}2+1)}\,
{}_3F_2\left[\matrix -N+1,\,z+w'+1,\,z'+w'+1\\ \s+2,\,x+w'-\frac{N-
1}2+1\endmatrix\Biggl|\,1\,\right],\\
h_{N-1}=\Gamma\left[\matrix N,\,\s+1,\,\s+2\\
\s+N+1,\,z+w+1,\,z+w'+1,\,z'+w+1,\,z'+w'+1\endmatrix\right],
\endgathered
\tag 7.2
$$
and $f(x)$ is given by \tht{6.1}.

Equivalently,
$$
\tK(x,y)=\frac
1{h_{N-1}}\,\frac{\wt \p_N(x)\p_{N-1}(y)-\p_{N-1}(x)\wt \p_N(y)}
{x-y}\,\sqrt{f(x)f(y)}\,,
\tag 7.3
$$
where
$$
\wt \p_N(x)=\frac{\Gamma(x+w'+\frac{N+1}2)}{\Gamma(x+w'-\frac{N-1}2)}\,
{}_3F_2\left[\matrix -N,\,z+w',\,z'+w'\\ \s+1,\,x+w'-\frac{N-
1}2\endmatrix\Biggl|\,1\,\right].
\tag 7.4
$$
\endproclaim

\demo{Comment} 1) If $\s=0$ then the formula for $\p_N$ above does not make
sense because it involves a hypergeometric function with a zero lower
index. The formula \tht{7.3} gives an explicit continuation of the
right--hand side of \tht{7.1} to the set $\s=0$.

2) For any $x\in\C$, the values $\p_{N-1}(x)$, $\wt \p_N(x)$, as well as
the constant $h_{N-1}$, are analytic functions in $(z,z',w,w')\in\Cal D$. The
value $\p_{N}(x)$ is an analytic function in $(z,z',w,w')\in\Cal
D\setminus\{\s=0\}$. If $(z,z',w,w')\in \Cal D_0$ then $h_{N-1}\ne 0$.
\enddemo

\demo{Proof of Theorem 7.1} It is convenient to set
$$
t=x+\tfrac{N-1}2, \quad
u=z+N-1, \quad u'=z'+N-1, \quad
v=w, \quad  v'=w'.
$$
Then $t$ ranges over the lattice $\Z$ and the function $f(x)$ turns
into the function
$$
g(t)=\frac1{\Gamma(u+1-t)\Gamma(u'+1-t)\Gamma(v+1+t)\Gamma(v'+1+t)}\,,
\quad t\in \Z\,.
$$

We consider $g(t)$ as a weight function on the lattice $\Z$. Note that
$$
g(t)\le \const\cdot (1+|t|)^{-u-u'-v-v'-2},\quad |t|\to\infty.
$$
Indeed, this follows from the estimate of the function $f$ given in
the beginning of the proof of Proposition 6.1.

We aim to study monic orthogonal polynomials $p_0=1,p_1,p_2,\dots$
corresponding to the weight function $g(t)$. In general, there
are only finitely many such polynomials, because, for any
(nonintegral) values of the parameters $u,u',v,v'$, the weight
function $g$ has only finitely many moments. This is a major
difference between our polynomial ensemble and polynomial ensembles
which are usually considered in the literature, cf. \cite{NW},
\cite{J}.

The number of existing polynomials with the weight function $g(t)$
depends on the number of finite moments of $g(t)$, i.e., on
$u+u'+v+v'$. Specifically, the $m$th polynomial exists if $g(t)$ has
finite moments up to the order $(2m-1)$ (this follows, e.g., from the
 determinantal formula expressing orthogonal polynomials
through the moments, see \cite{Er, Vol. 2, 10.3(4)}). This condition is 
satisfied if
$u+u'+v+v'>2m-2$. Let
$$
(p_m,p_m)=\sum_{t\in\Z} p^2_m(t)g(t)
$$
denote the square of the norm of the $m$th polynomial. This quantity
is well defined if $g(t)$ has finite $(2m)$th moment, which holds if
a slightly stronger condition is satisfied: $u+u'+v+v'>2m-1$.
(Note that it may happen that $p_m$ exists but $(p_m,p_m)=\infty$.)

\proclaim{Proposition 7.2} Set $\frak S=u+u'+v+v'$ and let
$m=0,1,\dots$.

If $\frak S>2m-2$ then
$$
p_m(t)=\frac{\Gamma(v'+1+t)}{\Gamma(v'+1+t-m)}\,
{}_3F_2\left[\matrix -m,\,u+v'+1-
m,\,u'+v'+1-m\\\goth S+2-2m,\,v'+1+t-m\endmatrix\biggl|\,1\,\right]\, .
$$

If $\frak S>2m-1$ then
$$
\gathered
(p_m,p_m)\equiv\sum_{t\in\Z}p_m^2(t)g(t)= \\
\Gamma\left[\matrix m+1,\, \goth S+1-2m,\, \goth S+2-2m\\
\goth S-m+2,\, u+v+1-m,\, u+v'+1-m,\, u'+v+1-m,\,
u'+v'+1-m\endmatrix\right]\,.
\endgathered
$$
\endproclaim

We will give the proof of Proposition 7.2 at the end
of this section. Now we proceed with the proof of Theorem 7.1.

Note that the condition $\Si=z+z'+w+w'>-1$ entering the definition of
$\Dadm$ is equivalent to the condition $\frak S=u+u'+v+v'>2N-3$. This
implies the existence of the $N$th Christoffel--Darboux kernel
associated with the weight function $g(t)$,
$$
\sum_{m=0}^{N-1}\frac{p_m(t_1)\,p_m(t_2)}{(p_m,p_m)}\,,
\qquad t_1,t_2\in\Z.
$$
Hence, the $N$th Christoffel--Darboux kernel associated with the
weight function $f$ has the form
$$
S^{(N)}(x,y):=\sum_{m=0}^{N-1}
\frac{p_m(x-\tfrac{N-1}2)\,p_m(y-\tfrac{N-1}2)}
{(p_m,p_m)}\,, \qquad x,y\in\x^{(N)}.
$$
By the general definition of the normalized Christoffel--Darboux
kernel (see \tht{5.2}),
$$
\tK(x,y)=S^{(N)}(x,y)\sqrt{f(x)f(y)}.
$$
On the other hand, let $T_1^{(N)}(x,y)$ denote the right--hand side of
\tht{7.1} with the term $\sqrt{f(x)f(y)}$ removed. Similarly, let
$T_2^{(N)}(x,y)$ denote the right--hand side of \tht{7.3} with the
term $\sqrt{f(x)f(y)}$ removed. It suffices to prove that
$$
S^{(N)}(x,y)=T_1^{(N)}(x,y)=T_2^{(N)}(x,y).
$$

\proclaim{Lemma 7.3} For any $x,y\in\x^{(N)}$, the kernel $S^{(N)}(x,y)$
viewed as a function in $(\zw)$, can be extended to a
holomorphic function on the domain $\Cal D_0$.
\endproclaim

Recall that the domain $\Cal D_0$ was introduced in \S3.

\demo{Proof} Indeed, the claim of the lemma holds for
$p_m(x-\tfrac{N-1}2)$, $p_m(y-\tfrac{N-1}2)$, and $(p_m,p_m)^{-1}$,
where $m=0,\dots,N-1$. This can be verified using explicit formulas of
Proposition 7.2. \qed
\enddemo

Let us prove first that $S^{(N)}(x,y)=T_1^{(N)}(x,y)$ under the
additional restriction $\Si>0$.
Then $\frak S>2N-2$, so that the $N$th polynomial $p_N$ exists. Hence,
the kernel $S^{(N)}$ can be written in the form (see Remark 5.2)
$$
\multline
S^{(N)}(x,y)\\
=\frac1{(p_{N-1},p_{N-1})}\,
\frac{p_N(x-\tfrac{N-1}2)p_{N-1}(y-\tfrac{N-1}2)-
p_{N-1}(x-\tfrac{N-1}2)p_N(y-\tfrac{N-1}2)}{x-y}\,.
\endmultline
$$
By explicit formulas of Proposition 7.2,
$$
p_N(x-\tfrac{N-1}2)=\p_N(x), \quad
p_{N-1}(x-\tfrac{N-1}2)=\p_{N-1}(x), \quad
(p_{N-1},p_{N-1})=h_{N-1}.
$$
This implies the desired equality when $\Si>0$.

Next, observe that the expression $T_1^{(N)}(x,y)$ is well defined
when $(\zw)$ ranges over the domain $\Cal D_0\setminus\{\Si=0\}$ and is a
holomorphic function on this domain. Together with Lemma 7.3 this
proves the equality $S^{(N)}(x,y)=T_1^{(N)}(x,y)$ provided that
$(\zw)\in \Cal D_0\setminus\{\Si=0\}$. As a consequence, we get
\tht{7.1} under the restriction $\Si\ne0$.

Finally, the series representation for ${}_3F_2$ easily implies the
identity
$$
{}_3F_2\left[\matrix a,\,b,\,c\,\\
e,\,f\endmatrix\,\Bigr|\,1\right]={}_3F_2\left[\matrix a,\,b,\,c\,\\
e+1,\,f\endmatrix\,\Bigr|\,1\right]+\frac{abc}{e(e+1)f}\,{}_3F_2\left[\matrix
a+1,\,b+1,\,c+1\,\\ e+2,\,f+1\endmatrix\,\Bigr|\,1\right].
$$
Applying this identity to $\p_N$ we see that $\p_N(x)=\wt
\p_N(x)+\const\p_{N-1}(x)$. It follows that $T_1^{(N)}(x,y)=
T_2^{(N)}(x,y)$ provided that $\p_N(x)$ makes sense. Thus, we see that
the singularity in $T_1^{(N)}(x,y)$ on the hyperplane $\Si=0$ is removable. 
Specifically, this
singularity is explicitly removed by means of the equality
$T_1^{(N)}(x,y)=T_2^{(N)}(x,y)$ (equivalently, by means of
\tht{7.3}). This completes the proof of Theorem 7.1. \qed
\enddemo

\demo{Proof of Proposition 7.2} We apply the general formalism
explained in \cite{NSU, Chapter 2}. Consider the following difference equation
on the lattice $\Z$:
$$
\sigma(t)\Delta\nabla y(t)+\tau(t)\Delta\, y(t) +\ga\, y(t)=0,
\tag 7.5$'$
$$
where
$$
\gathered
\sigma(t)=-(t+v)(t+v'),\\
\tau(t)=\goth S\, t+vv'-uu',\\
\nabla y(t)=y(t)-y(t-1),\quad \Delta y(t)=y(t+1)-y(t).
\endgathered
\tag 7.5$''$
$$
According to \cite{NSU} this equation is of hypergeometric type, that
is, $\sigma(t)$ is a polynomial of degree 2, $\tau(t)$ is a
polynomial of degree 1, and $\ga$ is a constant. The crucial fact is
that this equation can be rewritten in the selfadjoint form with the
weight function $g(t)$:
$$
(\Delta\circ\sigma g\circ\nabla)\,y+\ga\,g y=0,
$$
which easily follows from the relation
$$
\frac{g(t)}{g(t+1)}=\frac{\sigma(t+1)}{\sigma(t)+\tau(t)}\,.
$$
(Note that in \cite{NSU} the weight function is denoted by $\rho$
and the spectral parameter $\ga$ is denoted by $\la$.)

We will seek $p_m$ as a monic polynomial of degree $m$, which
satisfies the difference equation above with an appropriate value of
$\ga$. According to \cite{NSU}, this value must be equal to
$$
\ga=\ga_m=-m\tau'-\frac{m(m-1)}2\sigma'',
$$
where the derivatives $\tau'$ and $\sigma''$ are constants, because
$\tau$ has degree 1 and $\sigma$ has degree 2. Further, the polynomial in
question exists and is unique provided that
$$
\mu_k:=\ga_m+k\tau'+\frac{k(k-1)}2\sigma''\ne0,
\qquad k=0,\dots,m-1.
$$

In our case $\tau'=\frak S$, $\sigma''=-2$, so that
$$
\ga_m=m(m-1-\frak S)
$$
and the nonvanishing condition $\mu_k\ne0$ turns into
$$
\frak S \ne m+k-1, \qquad k=0,\dots,m-1,
$$
which is certainly satisfied if $\frak S>2m-2$.

Our additional condition $\s>0$ is equivalent to $\frak S>2N-2$. It
follows that the desired polynomial solutions $p_m$ certainly exist for $m\le
N$.

Following the notation of \cite{NSU} set
$$
\gather
A_{mm}=m!\,\prod_{k=0}^{m-1}(\tau'+\tfrac{m+k-1}2\sigma''), \\
g_m(t)=g(t)\,\prod_{l=1}^m \sigma(t+l),\\
S_m=\sum_{t\in\Z} g_m(t).
\endgather
$$
Then, according to \cite{NSU},
$$
\gather
p_m(x)=\frac{(-1)^m\, m!}{A_{mm}} \,
\sum_{k=0}^m\frac{(-m)_kg_m(x-m+k)}{k!\, g(x)}\,,\\
(p_m,p_m)=\frac{(-1)^m (m!)^2S_m}{A_{mm}}\,.
\endgather
$$

Applying this recipe in our concrete case we get
$$
\gather
A_{mm}=m!\,\frac{\Ga(\frak S -m+2)}{\Ga(\frak S-2m+2)}\,, \\
g_m(t)=\frac{(-1)^m}{\Ga(u+1-m-t)\Ga(u'+1-m-t)
\Ga(v+1+t)\Ga(v'+1+t)}
\endgather
$$
and, by Dougall's formula \cite{Er, Vol.1, 1.4(1)},
$$
\multline
S_m=\sum_{t\in\Z}g_m(t)\\
=(-1)^m\,\Ga\,\bmatrix \frak S+1-2m \\
u+v+1-m,\, u'+v'+1-m,\, u'+v+1-m\,, u+v'+1-m \endbmatrix\,.
\endmultline
$$
This implies the expression for $(p_m,p_m)$ given in Proposition 7.2.

Further, after simple transformations we get
$$
\aligned
p_m(t)=&\,\Ga\bmatrix \frak S -2m+2, \, v+1+t, \, v'+1+t\\
\frak S-m+2, \, v+1-m+t, \,v'+1-m+t \endbmatrix
\\ &\times \sum_{k=0}^m \frac{(-m)_k(t-u)_k(t-u')_k}
{(t+v+1-m)_k(t+v'+1-m)_k k!}\\
=&\,\Ga\bmatrix \frak S -2m+2, \, v+1+t, \, v'+1+t\\
\frak S-m+2, \, v+1-m+t, \,v'+1-m+t \endbmatrix
\\ &\times{}_3F_2\left[\matrix -m,\, t-u,\, t-u'\\
t+v+1-m,\, t+v'+1-m \endmatrix \biggl|\,1\,\right]\,.
\endaligned
$$
Applying the transformation
$$
{}_3F_2\left[\matrix a,\, b,\, c\\
e,\, f \endmatrix \biggl|\,1\,\right]=\Ga\bmatrix e, \, e+f-a-b-c\\
e+f-b-c, \, e-a \endbmatrix {}_3F_2\left[\matrix a,\, f-b,\, f-c\\
e+f-b-c,\, f \endmatrix \biggl|\,1\,\right]\,,
\tag 7.6
$$
see Appendix, to the last expression we
arrive at the desired formula for $p_m(t)$. \qed

\enddemo

The orthogonal polynomials $\{p_m\}$ were discovered by
R.~Askey \cite{As}. Later they were independently found by Peter A.~Lesky,
see \cite{Les1}, \cite{Les2}.  We are grateful to Tom Koornwinder for informing
us about Lesky's work, and to Peter A. Lesky for sending us his preprint
\cite{Les2}.

\head 8. The correlation kernel of the process $\Cal P^{(N)}$
\endhead

In this section we study the process $\Cal P^\n$ on the lattice
$\x^\n$ introduced in \S4. In \S6 we showed that $\Cal P^\n$ is an
L--ensemble, hence, it is a determinantal point process. Denote by
$K^\n$ the $\x^\n\times\x^\n$ matrix defined by
$K^\n=L^\n(1+L^\n)^{-1}$, where $L^\n$ was defined in \S6. By
Proposition 5.3, $K^\n$ is a correlation kernel for $\Cal P^\n$. Our
goal in this section is to provide analytic expressions for $K^\n$
suitable for a future limit transition as $N\to\infty$.

Our analysis is based on the general results of \S5(f)--(h), where we
take $\x=\x^\n$, $\x_I=\xout$, $\x_{II}=\xin$. To compute the kernel
$K^\n$ we use the method of \S5(f). Proposition 5.6 expresses the
kernel in terms of four functions $R_I$, $R_{II}$, $S_I$, $S_{II}$,
which solve a discrete Riemann--Hilbert problem. In our concrete
situation we will redenote these functions by $\Rout^\n$, $\Rin^\n$,
$\Sout^\n$, $\Sin^\n$.

In order to apply this method we need to impose an additional
restriction on the parameters. Specifically, we require that
$(\zw)\in\Dadm'\cap\{\Si>0\}$, where $\Dadm'$ was defined in
Definition 6.2 and $\Si=z+z'+w+w'$. Later we show that the final
expression for $K^\n$ remains valid without the restriction $\Si>0$.

Thanks to Proposition 5.8, the functions $\Rout^\n$ and $\Sout^\n$ are
immediately expressed through the orthogonal polynomials $\p_N$,
$\p_{N-1}$ evaluated in \S7. The remaining two functions $\Rin^\n$ and
$\Sin^\n$ are uniquely determined by their relations to $\Rout^\n$
and $\Sout^\n$, see \tht{5.12} and \tht{8.3} below. We provide
certain explicit expressions for $\Rin^\n$, $\Sin^\n$ and check that
they satisfy the needed relations. About the derivation of these
expressions, see Remark 8.9 below.

To remove the restriction $\Si>0$ we use Theorem 5.10 and Lemma 7.3.
These results show that the kernel $K^\n$ divided by a simple
factor, admits a holomorphic continuation (as a function of the
parameters) to a certain domain $\Cal D'_0\supset\Dadm'$.

Now we proceed to the realization of of the plan described above.

We assume that $(\zw)\in\Dadm'$. As was mentioned in \S6, this
guarantees that the weight function $f(x)$ does not vanish on $\xin$.
Next, we temporarily assume that $\Si>0$. This ensures that the
function $h_I$ belongs to $\ell^2(\xin)$ as required in \S5(h).

Following Proposition 5.8 and using Theorem 7.1, we define 2 meromorphic 
functions
$\Rout^{(N)}$ and $\Sout^{(N)}$ on the complex plane with poles in $\xin$ as
follows:  $$
\multline
\Rout^{(N)}(x)=\frac{\p_N(x)}{(x-\frac{N-1}2)\cdots(x+\frac{N-1}2)}\\=
\Gamma\bmatrix
x-\frac{N-1}2,\,x+w'+\frac{N+1}2\\x+\frac{N+1}2,\,x+w'-\frac{N-1}2\endbmatrix
{}_3F_2\left[\matrix -N,\,z+w',\,z'+w'\\\Sigma,\,x+w'-\frac{N-
1}2\endmatrix\biggl|\,1\,\right]\,,
\endmultline
\tag 8.1
$$
and
$$
\multline
\Sout^{(N)}(x)=\frac{\p_{N-1}(x)}
{h_{N-1}(x-\frac{N-1}2)\cdots(x+\frac{N-1}2)}\\=
\frac 1{h_{N-1}}\,\Gamma\bmatrix
x-\frac{N-1}2,\,x+w'+\frac{N+1}2\\x+\frac{N+1}2,\,x+w'-\frac{N-3}2\endbmatrix
{}_3F_2\left[\matrix -N+1,\,z+w'+1,\,z'+w'+1\\
\Sigma+2,\,x+w'-\frac{N-3}2\endmatrix\biggl|\,1\,\right]\,,
\endmultline
\tag 8.2
$$
where $h_{N-1}$ was defined in \tht{7.2}. Observe that
$$
\Rout^{(N)}(x)=1+O(x^{-1}),\quad \Sout^{(N)}=O(x^{-1}),\qquad x\to\infty.
$$

Note that the right--hand side of \tht{8.1} makes sense for any
$x\in\C\setminus\xin$ and $(z,z',w,w')\in\Cal D\setminus\{\s=0\}$, and
the right--hand side of \tht{8.2} makes sense for any $x\in\C\setminus\xin$
and $(z,z',w,w')\in\Cal D_0$.

Using \tht{5.12} as a prompt, we set
$$
\Sin^{(N)}(x)=-\sum_{y\in\xout}\frac
{\psiout^{(N)}(y)\Rout^{(N)}(y)}{y-x}\,,\quad
\Rin^{(N)}(x)=1-\sum_{y\in\xout}\frac
{\psiout^{(N)}(y)\Sout^{(N)}(y)}{y-x}\,,
\tag 8.3
$$
where the functions $\psiout^\n$ and $\psiin^\n$ were introduced in
\S6.

Note that $\psiout^{(N)}(x)\le\const\, (1+|x|)^{-\s}$ (this follows
from the estimate of the weight function $f$, see \S6). Since, by
assumption, $\Si>0$, it follows that the series \tht{8.3} converge
and define meromorphic functions with poles in $\xout$.

\proclaim{Proposition 8.1} We have
$$
\aligned
\Rin^{(N)}(x)=&-\frac {\sin \pi z}{\pi}\,\Gamma\left[\matrix z'-
z,\,z+w+1,\,z+w'+1\\ \s+1\endmatrix\right]\\ &\times
\Gamma\left[\matrix
x+\frac{N+1}2,\,-x+\frac{N+1}2,\,N+1+\s\\-
x+z'+\frac{N+1}2,\,x+w+\frac{N+1}2,N+1+z+w'
\endmatrix\right]\\
&\times
{}_3F_2\left[\matrix  z+w'+1,\,
-z'-w,\,-x+z+\frac{N+1}2\\z-z'+1,\,N+1+z+w'\endmatrix\,\Biggl|\,1\right]\\ &
-\{ \text{ similar expression with  } z \text{  and  } z' \text{
interchanged }\},
\endaligned
\tag 8.4
$$
$$
\aligned
\Sin^{(N)}(x)=&-\frac {\sin \pi z}{\pi}\,\Gamma\left[\matrix z'-
z,\,\s\\ z'+w,\,z'+w'\endmatrix\right]\\ &\times
\Gamma\left[\matrix
x+\frac{N+1}2,\,-x+\frac{N+1}2,\,N+1\\-
x+z'+\frac{N+1}2,\,x+w+\frac{N+1}2,N+1+z+w'
\endmatrix\right]\\
&\times
{}_3F_2\left[\matrix  -z'-w+1,\,
z+w',\,-x+z+\frac{N+1}2\\z-z'+1,\,N+1+z+w'\endmatrix\,\Biggl|\,1\right]\\ &
-\{ \text{ similar expression with  } z \text{  and  } z' \text{
interchanged }\}.
\endaligned
\tag 8.5
$$
\endproclaim
\demo{Singularities}
Using the structure of singularities of a general ${}_3F_2$ function, see the
beginning of \S7, it is easy to verify that the formulas above make sense for
$x\in\C\setminus \xout$ and
$(z,z',w,w')\in\Cal D_0\setminus\left(\{\s=0\}\cup\{z-z'\in\Z\}\right).$
$\{\s=0\}$ is indeed a singular set for $\Sin^{(N)}$. The singularities
$\{z-z'\in\Z\}$, however, are removable, as long as
we are in $\Cal D_0\setminus\{\s=0\}$. For $\s>0$ this follows from the
definitions \tht{8.3}.
\enddemo

\demo{Proof} Since both parts of \tht{8.4} and \tht{8.5} are analytic in the
parameters, it is enough to give a proof for nonintegral values of
$z,z',w,w'$, such that $z-z'\notin \Z$, $\s\notin \Z$.

We start with rewriting the expressions above in a somewhat more
suitable form.

Introduce a meromorphic function
$$
\multline
F(x)=\frac 1{\sin(\pi(z'-z))\sin(\pi(z+z'+w+w'))}\\ \times
\left(\dfrac{\sin(\pi(z+w))\sin(\pi(z+w'))\sin(\pi
z')}
{\sin(\pi(-x+z'+\frac{N+1}2))}-\dfrac{\sin(\pi(z'+w))\sin(\pi(z'+w'))\sin(\pi
z)}{\sin(\pi(-x+z+\frac{N+1}2))}\right)\,.
\endmultline
$$
It is not difficult to see that if $x\in\x$, that is, if
$x-\frac{N-1}2\in \Z$, then $F(x)=(-1)^{{x-\frac{N-1}2}}$.

Let us also introduce a more detailed notation $h(N-1,z,z',w,w')$ for the
constant $h_{N-1}$ defined in \tht{7.2}.

\proclaim{Lemma 8.2} The right--hand side of \tht{8.4} can be
written in the form
$$
\gathered
\Gamma\bmatrix
x+\frac{N+1}2,\,x-z-\frac{N-1}2\\x-\frac{N-1}2,\,x-
z+\frac{N+1}2\endbmatrix{}_3F_2\left[\matrix N,\,-z-w',\,-z-w\\-\Sigma,\,x-
z+\frac{N+ 1}2\endmatrix\biggl|\,1\,\right]\\
+\Gamma\bmatrix x+\frac{N+1}2,\,-x+\frac{N+1}2\\
-x+z+\frac{N+1}2,\,-x+z'+\frac{N+1}2,\,x+w+\frac{N+1}2,\,
x+w'-\frac{N-3}2\endbmatrix\\ \times \frac{F(x)}{h(N-1,z,z',w,w')}\
{}_3F_2\left[\matrix -N+1,\,z+w'+1,\,z'+w'+1\\ \Sigma+2,\,x+w'-\frac{N-
3}2\endmatrix\biggl|\,1\,\right],
\endgathered
\tag 8.6
$$
and the right--hand side of \tht{8.5} can be
written in the form
$$
\gathered
h\left(N,z-\frac 12,z'-\frac 12,w-\frac 12,w'-\frac 12\right)
\Gamma\bmatrix
x+\frac{N+1}2,\,x-z-\frac{N-1}2\\x-\frac{N-1}2,\,x-z+\frac{N+3}2\endbmatrix\\
\times {}_3F_2\left[\matrix N+1,\,-z-w'+1,\,-z-w+1\\-\Sigma+2,\,x-z+\frac{N+
3}2\endmatrix\biggl|\,1\,\right]\\ + \Gamma\bmatrix
x+\frac{N+1}2,\,-x+\frac{N+1}2\\
-x+z+\frac{N+1}2,\,-x+z'+\frac{N+1}2,\,x+w+\frac{N+1}2,\,
x+w'-\frac{N-1}2\endbmatrix\\ \times F(x)\  {}_3F_2\left[\matrix
-N,\,z+w',\,z'+w'\\ \Sigma,\,x+w'-\frac{N-1}2\endmatrix\biggl|\,1\,\right]\,.
\endgathered
\tag 8.7
$$
\endproclaim
\demo{Singularities} It looks like \tht{8.6} may have more
singular points that the right--hand side of \tht{8.4}. For
example, the first summand in \tht{8.6} has poles at the points, where
$$
x-u-\tfrac{N-1}2\in\{0,-1,\dots\}, \quad u=z\text{  or  } z',
$$
and same is true about the second
summand. The fact that these poles cancel out is not obvious.
Similar cancellations happen in \tht{8.7}.
\enddemo

The proof of the lemma can be found in the Appendix. It is rather tedious and
is based on the known transformation formulas for ${}_3F_2$ with the unit
argument.

Our next step is to prove the following statement.

\proclaim{Lemma 8.3} The only singularities of the right--hand of \tht{8.4}
regarded as a function in $x\in \C$, are simple poles at the points of
$\xout$. The residue of the right--hand side of \tht{8.4} at any point
$x\in\xout$ equals $\psiout^{(N)}(x)\Sout^{(N)}(x)$. Similarly, the only
singularities of the right--hand of \tht{8.5} regarded as a function in $x\in
\C$, are simple poles at the points of $\xout$. The residue of the right--hand
side of \tht{8.5} at any point $x\in\xout$ equals
$\psiout^{(N)}(x)\Rout^{(N)}(x)$.
\endproclaim

\demo{Proof}
The location of poles part follows from the general structure of singularities
of ${}_3F_2$, see the beginning of \S7.

To evaluate the residue of the right--hand side of \tht{8.4} we will use the
formula \tht{8.6}. We easily see that first summand of \tht{8.6} takes finite
values on $\xout$. As for
the second summand of \tht{8.6}, the ${}_3F_2$ is a terminating series and it
has no singularities in $\xout$. Furthermore, $F(x)=(-1)^{{x-\frac{N-1}2}}$
for $x\in \x$, in particular, for $x\in \xout$. It remains to examine the
residue of $\Ga(x+\frac{N+1}2)\Ga(-x+\frac{N+1}2)$. 

Assume that $x\in\xout$ and $x>\frac{N-1}2$. Then
$$
\ru\Ga\left(-u+\tfrac{N+1}2\right)=
\frac{(-1)^{x-\frac{N-1}2}}{\Gamma(x-\frac{N-1}2)}.
$$
Thus, the residue of \tht{8.6} at $x$ is equal to
$$
\gathered
\Gamma\bmatrix x+\frac{N+1}2\\
x-\frac{N-1}2,\,
-x+z+\frac{N+1}2,\,-x+z'+\frac{N+1}2,\,x+w+\frac{N+1}2,\,
x+w'-\frac{N-3}2\endbmatrix\\ \times \frac{1}{h(N-1,z,z',w,w')}\
{}_3F_2\left[\matrix -N+1,\,z+w'+1,\,z'+w'+1\\ \Sigma+2,\,x+w'-\frac{N-
3}2\endmatrix\biggl|\,1\,\right].
\endgathered
$$
By a direct comparison we see that this is equal to
$\psiout^{(N)}(x)\Sout^{(N)}(x)$.

Similarly, if $x\in\xout$ and $x<-\frac{N-1}2$, we have
$$
\ru\Ga\left(u+\tfrac{N+1}2\right)=
\frac{(-1)^{-x-\frac{N+1}2}}{\Gamma(-x-\frac{N-1}2)},
$$
and the residue of \tht{8.6} at $x$ equals
$$
\gathered
\Gamma\bmatrix -x+\frac{N+1}2\\
-x-\frac{N-1}2,\,
-x+z+\frac{N+1}2,\,-x+z'+\frac{N+1}2,\,x+w+\frac{N+1}2,\,
x+w'-\frac{N-3}2\endbmatrix\\ \times \frac{(-1)^N}{h(N-1,z,z',w,w')}\
{}_3F_2\left[\matrix -N+1,\,z+w'+1,\,z'+w'+1\\ \Sigma+2,\,x+w'-\frac{N-
3}2\endmatrix\biggl|\,1\,\right],
\endgathered
$$
which is again equal to $\psiout^{(N)}(x)\Sout^{(N)}(x)$.

Thus, we have proved the first statement of the lemma. The proof of the second
one is very similar. \qed
\enddemo

Lemma 8.3 shows that the right--hand sides of \tht{8.3} and of \tht{8.4},
\tht{8.5} have the same singularity structure. Clearly, the sums in the
right--hand sides of \tht{8.3} decay as $x\to\infty$ and $x$ keeps finite
distance from the points of $\xout$. We aim to prove that the right--hand
sides of \tht{8.4} and \tht{8.5} have a similar property. Using Lemma 8.2, we
may consider the expressions \tht{8.6} and \tht{8.7}.

Consider the formula \tht{8.6}. We first assume that $\Re x\ge0$. Let
us examine the first summand. The 
gamma--factors form a rational function which tends to 1 uniformly in
$\operatorname{arg}(x)$ as $|x|\to\infty$.

In order to handle the asymptotics of the ${}_3F_2$ factor, let us apply the
formula \tht{7.6} to it with
$$
a=-z-w,\quad b=-z-w',\quad c=N,\quad e=x-z+\frac{N+1}2,\quad f=-\s.
$$
We get
$$
\gathered
{}_3F_2\left[\matrix N,\,-z-w',\,-z-w\\
-\Sigma,\,x-z+\frac{N+1}2\endmatrix\biggl|\,1\,\right]=
\Ga\bmatrix x-z+\frac{N+1}2,\, x-z'-\frac{N-1}2\\
x-z-z'-w-\frac{N-1}2,\,x+w+\frac{N+1}2\endbmatrix\\ \times
{}_3F_2\left[\matrix -z-w,\,-z'-w,\,-\s-N\\
-\Sigma,\,x-z-z'-w-\frac{N-1}2\endmatrix\biggl|\,1\,\right].
\endgathered
$$
Using the asymptotic relation
$$
\frac{\Gamma(x+\alpha)}{\Ga(x+\beta)}=x^{\al-\be}(1+O(|x|^{-1})),
\tag 8.8
$$
which is uniform in $\operatorname{arg}(x)\in
(-\pi+\varepsilon,\pi-\varepsilon)$ for any $\varepsilon>0$, we see that the
gamma-factors tend to 1 as $|x|\to \infty$.

As for the ${}_3F_2$, we see that
the sum of the lower parameters minus the sum of the upper parameters is equal
to $x+w+\frac{N+1}2$. Then if $\Re x\ge 0$ and $\Re w+\frac{N+1}2>0$, the
defining series for this ${}_3F_2$ converges uniformly in $x$, provided that
the distance from $x$ to the lattice $\Z+z+z'+w+\frac{N-1}2$ is bounded from
zero. One proof of this fact can be obtained by applying the general estimate
(see \cite{Er, Vol. 1, 1.9(8)})
$$
c_1\alpha^{\be-\ga}<
\frac{\Ga(\alpha+\beta)}{\Ga(\alpha+\gamma)}<c_2\al^{\be-\ga},\qquad
\alpha>0,\quad \Re\be,\,\Re\ga>\operatorname{const}>0,
$$
$c_1,c_2>0$ are some fixed constants, to the terms of the series.

Thus, we can pass to the limit $|x|\to\infty$ term-wise, and since all terms
of the series starting from the second one converge to zero uniformly in
$\operatorname{arg}(x)$, we conclude, that the first summand of \tht{8.6}, as
$|x|\to\infty$, converges to 1 uniformly in $x$ such that $\Re x\ge 0$ and the
distance from $x$ to the lattice $\Z+z+z'+w+\frac{N-1}2$ is bounded from zero.

Let us proceed to the second summand of \tht{8.6}. Now the ${}_3F_2$ part is a
rational function which tends to 1 as $|x|\to\infty$ uniformly in
$\operatorname{arg}(x)$. The remaining part --- the product of gamma-factors
and $F(x)$ --- can be written in the form
$$
\gathered
\const_1\,\Gamma\bmatrix x+\frac{N+1}2,\,-x+\frac{N+1}2,\,x-z-\frac{N-1}2\\
-x+z'+\frac{N+1}2,\,x+w+\frac{N+1}2,\,
x+w'-\frac{N-3}2\endbmatrix \\+ \const_2
\,\Gamma\bmatrix x+\frac{N+1}2,\,-x+\frac{N+1}2,\,x-z'-\frac{N-1}2\\
-x+z+\frac{N+1}2,\,x+w+\frac{N+1}2,\,
x+w'-\frac{N-3}2\endbmatrix,
\endgathered
\tag 8.9
$$
where for $F(x)$ we used the formulas ($u=z$ or $z'$)
$$
\frac 1{\sin(\pi(-x+u+\frac{N+1}2))}=\frac
{\Ga(-x+u+\frac{N+1}2)\Ga(x-u-\frac{N-1}2)}\pi\,.
$$
Further, using the relations ($u=z$ or $z'$)
$$
\frac{\Ga\left(-x+\frac{N+1}2\right)}{\Ga\left(-x+u+\frac{N+1}2\right)}=\frac
{\sin(\pi(-x+u+\frac{N+1}2))\,\Ga(x-u-\frac{N-1}2)}{\sin(\pi(-
x+\frac{N+1}2))\,\Ga(x-\frac{N-1}2)}\,,
$$
observing that the ratio of sines is bounded as long as $x$ is bounded from
the lattice $\x=\Z+\frac{N-1}2$, and employing \tht{8.8}, we see that
the absolute value of \tht{8.9} is bounded by a constant times $|x|^{-\s-1}$,
and the bound is uniform in $x$ such that $\operatorname{arg}(x)$ is bounded
from $\pm\pi$ and $x$ is bounded from $\x$. Since $\s+1>0$, the second term of
\tht{8.6} tends to zero.

We conclude that, under the condition $\Re w+\frac{N+1}2>0$, the expression
\tht{8.6} converges to 1
as $|x|\to \infty$ uniformly in $x$ such that $\Re x\ge 0$ and the distance
from $x$ to the lattices $\Z+\frac{N-1}2$ and $\Z+z+z'+w+\frac{N-1}2$ is
bounded from zero.

To extend this estimate to the domain $\Re x\le 0$ we use the following
\proclaim{Lemma 8.4} The expression \tht{8.6} is invariant with respect to the
following change of variable and parameters:
$$
x\mapsto -x, \qquad (z,z')\longleftrightarrow (w',w).
$$
\endproclaim

We give a proof of this lemma in the Appendix. Similarly to the proof of Lemma
8.2, it is rather technical and is based on certain transformation formulas
for ${}_3F_2$.

Lemma 8.4 immediately implies that, under the conditions $\Re
w+\frac{N+1}2>0$, $\Re z'+\frac{N+1}2>0$, the expression \tht{8.6} converges
to 1 as $|x|\to \infty$ uniformly in $x$ such that the distance from $x$ to the
lattices $\Z+\frac{N-1}2$, $\Z+z+z'+w+\frac{N-1}2$, $\Z+z'+w+w'+\frac{N-1}2$
is bounded from zero.

This statement together with Lemma 8.3 shows that the right--hand side of the
second formula of \tht{8.3} and the right--hand side of \tht{8.4} have the
same singularities and asymptotics at infinity. Thus, they must be equal.
This proves \tht{8.4} under the additional conditions $\Re
w+\frac{N+1}2>0$, $\Re z'+\frac{N+1}2>0$. Since both sides of \tht{8.4} depend
on the parameters $z,z',w,w'$ analytically, the additional conditions may
be removed. The proof of \tht{8.4} is complete.

The proof of \tht{8.5} is very similar. It is based on the following technical
statement, which will also be addressed in the Appendix.

\proclaim{Lemma 8.5} The expression \tht{8.7} is skew-invariant with respect to

$$
x\mapsto -x, \qquad (z,z')\longleftrightarrow (w',w). \qed
$$
\endproclaim

The proof of Proposition 8.1 is now complete.\qed
\enddemo

\example{Remark 8.6} Having proved the formulas \tht{8.6}, \tht{8.7}, it is
easy to verify \tht{5.21} directly. Indeed, if $x\in\xin$
then the first summands of \tht{8.6} and \tht{8.7} vanish thanks to
$\Ga(x+\frac{N+1}2)/\Ga(x-\frac{N-1}2)$, while in the second ones we have
$F(x)=(-1)^{x-\frac{N-1}2}$. A direct comparison of the resulting expressions
with \tht{8.1} and \tht{8.2} yields \tht{5.21}.
\endexample

\proclaim{Theorem 8.7} Let $(z,z',w,w')\in\Dadm'$ {\rm(}see
Definition \tht{6.3}{\rm)}, $\Cal P^\n$ be the
corresponding point process defined in \S4, and
$K^\n=L^\n(1+L^\n)^{-1}$, where $L^\n$ is as in \S6.

Then the kernel $K^\n(x,y)$, represented in the block form
corresponding to the splitting $\x=\xout\sqcup\xin$, is equal to
$$
\gathered
K^{(N)}_{\out,\out}(x,y)=\sqrt{\psiout^{(N)}(x)\psiout^{(N)}(y)}
\,\frac {\Rout^{(N)}(x)\Sout^{(N)}(y)-\Sout^{(N)}(x)\Rout^{(N)}(y)}{x-y}\,,\\
K^{(N)}_{\out,\inr}(x,y)=\sqrt{\psiout^{(N)}(x)\psiin^{(N)}(y)}
\,\frac {\Rout^{(N)}(x)\Rin^{(N)}(y)-\Sout^{(N)}(x)\Sin^{(N)}(y)}{x-y}\,,\\
K^{(N)}_{\inr,\out}(x,y)=\sqrt{\psiin^{(N)}(x)\psiout^{(N)}(y)}
\,\frac {\Rin^{(N)}(x)\Rout^{(N)}(y)-\Sin^{(N)}(x)\Sout^{(N)}(y)}{x-y}\,,\\
K^{(N)}_{\inr,\inr}(x,y)=\sqrt{\psiin^{(N)}(x)\psiin^{(N)}(y)}
\,\frac {\Rin^{(N)}(x)\Sin^{(N)}(y)-\Sin^{(N)}(x)\Rin^{(N)}(y)}{x-y}\,,
\endgathered
$$
where the functions $\Rout^{(N)}$, $\Sout^{(N)}$, $\Rin^{(N)}$, $\Sin^{(N)}$
are given by \tht{8.1}, \tht{8.2}, \tht{8.4}, \tht{8.5}, the
functions $\psiout^\n$, $\psiin^\n$ are given by \tht{6.2},
\tht{6.3}, and the indeterminacy on the diagonal $x=y$ is resolved by
the L'Hospital rule:
$$
\gathered
K^{(N)}_{\out,\out}(x,x)=\psiout^{(N)}(x)
\left((\Rout^{(N)})'(x)\Sout^{(N)}(x)-(\Sout^{(N)})'(x)\Rout^{(N)}(x)\right),
\\
K^{(N)}_{\inr,\inr}(x,x)=\psiin^{(N)}(x)
\left((\Rin^{(N)})'(x)\Sin^{(N)}(x)-(\Sin^{(N)})'(x)\Rin^{(N)}(x)\right).
\endgathered
$$
\endproclaim

\demo{Singularities} We know that the functions $\Rout^{(N)}$ and $\Sin^{(N)}$
are singular when $\s=0$, see the beginning of the section. However, as will
be clear from the proof, the value $K^{(N)}(x,y)$ is a well--defined
continuous function on the whole $\Dadm$ including the set $\{\s=0\}$, for any
$x,y\in\x^\n$.
\enddemo

\demo{Proof} Under the additional restriction $\s>0$ the statement of
the theorem follows from the above results. Indeed, the functions
$\Rout^\n$, $\Sout^\n$, $\Rin^\n$, $\Sin^\n$ defined by \tht{8.1},
\tht{8.2}, \tht{8.3} solve the discrete Riemann--Hilbert problem
associated with $L^\n$, see \tht{5.18} and \tht{5.12}. Proposition
5.6 explains how the kernel $K^\n$ is written in terms of these
functions. Proposition 8.1 provides explicit expressions \tht{8.4},
\tht{8.5} for $\Rin^\n$, $\Sin^\n$, and this completes the proof when
$\Si>0$.

Now we want to get rid of this restriction.

Introduce a function $\Psi^{(N)}:\x^{(N)}\to\C$ by
$$
\Psi^{(N)}(x)=\cases \psiout^{(N)}(x),&x\in\xout,\\
\dfrac 1{\psiin^{(N)}(x)}\,,&x\in\xin.
\endcases
$$
We have
$$
\Psi^{(N)}(x)=f(x)\cdot\cases \left((x-\tfrac{N-1}2)_N\right)^2,&x\in\xout,\\
\left(\Ga(-x+\tfrac{N+1}2)\Ga(x+\tfrac{N+1}2)\right)^2,&x\in\xin,\endcases
\tag8.10
$$
where $f(x)$ was defined in \tht{6.1}.
Note that $\Psi^{(N)}(x)$ is entire in $(z,z',w,w')$ for any $x\in
\x$.

Set
$$
\Cal D'_0=\{(\zw)\in\Cal D_0\mid z,z',w,w'\ne -1,-2,\dots\}.
$$
Recall that $\Cal D_0$ contains $\Dadm$ and note that
$\Cal D'_0\cap\Dadm=\Dadm'$. In particular, $\Cal D'_0$ is a domain
in $\C^4$ containing $\Dadm'$.

\proclaim{Lemma 8.8} Let $(\zw)\in\Dadm'$. The kernel $K^\n(x,y)$ can
be written in the form
$$
K^\n(x,y)=\sqrt{\Psi^\n(x)\Psi^\n(y)}\Kcirc^\n(x,y), \tag8.11
$$
where $\Kcirc^\n(x,y)$ can be extended {\rm (}as a function in $(\zw)${\rm )} to 
a
holomorphic function on the domain $\Cal D'_0$.
\endproclaim

\demo{Proof of the lemma} Recall (see \S7) that
$$
\tK(x,y)=S^\n(x,y)\sqrt{f(x)f(y)},
$$
where $S^\n$ is the (ordinary, not normalized) Christoffel--Darboux
kernel, which admits a holomorphic continuation to the domain $\Cal
D_0$ (see Lemma 7.3). Theorem 5.10 provides a connection between the
kernels $K^\n$ and $\tK$. Specifically, \tht{5.22} reads
$$
K^\n(x,y)=\ep(x)(\tK)^\triangle(x,y)\ep(y),
$$
where the symbol $(\,\cdot\,)^\triangle$ is explained in \S5(c).
Since $\ep(\,\cdot\,)$ does not depend on $(\zw)$, it suffices to
prove that the kernel $(\tK)^\triangle(x,y)$ can be represented as
the product of $\sqrt{\Psi^\n(x)\Psi^\n(y)}$ and a holomorphic
function on $\Cal D'_0$.

By \tht{8.10}, $\Psi^\n(x)$ differs from $f(x)$ by a positive factor
not depending on $(\zw)$. Hence, it is enough to prove that
$(\tK)^\triangle(x,y)$ can be written as the product of
$\sqrt{f(x)f(y)}$ and a holomorphic function on $\Cal D'_0$.

By \tht{5.1} we have
$$
(\tK)^\triangle(x,y)=\cases f(x)\left(\frac1{f(x)}-S^\n(x,x)\right),
& x=y\in\xin\,, \\
\pm\sqrt{f(x)f(y)}S^\n(x,y), & \text{otherwise},
\endcases
$$
where the sign `$\pm$' does not depend on $(\zw)$. {}From the
definition of $f(x)$ it follows that for $x\in\xin$, $\frac1{f(x)}$ admits a
holomorphic extension (as a function in the parameters) to $\Cal
D'_0$. Since $S^\n(x,y)$ is holomorphic on $\Cal D_0$, the lemma
follows. \qed
\enddemo

Now we complete the proof of Theorem 8.7. Since we have already
proved its claim when $\Si>0$, we see that under this restriction
$$
\gathered
\Kcirc^{(N)}_{\out,\out}(x,y)=
\frac {\Rout^{(N)}(x)\Sout^{(N)}(y)-\Sout^{(N)}(x)\Rout^{(N)}(y)}{x-y}\,,\\
\Kcirc^{(N)}_{\out,\inr}(x,y)=\Psi^{(N)}(y)
\,\frac {\Rout^{(N)}(x)\Rin^{(N)}(y)-\Sout^{(N)}(x)\Sin^{(N)}(y)}{x-y}\,,\\
\Kcirc^{(N)}_{\inr,\out}(x,y)=\Psi^{(N)}(x)
\,\frac {\Rin^{(N)}(x)\Rout^{(N)}(y)-\Sin^{(N)}(x)\Sout^{(N)}(y)}{x-y}\,,\\
\Kcirc^{(N)}_{\inr,\inr}(x,y)=\Psi^{(N)}(x)\Psi^{(N)}(y)
\,\frac {\Rin^{(N)}(x)\Sin^{(N)}(y)-\Sin^{(N)}(x)\Rin^{(N)}(y)}{x-y}\,,
\endgathered
$$

Since the functions $\Rout^\n$, $\Sout^\n$, $\Rin^\n$, $\Sin^\n$, and
$\Psi^\n$ are all holomorphic on $\Cal D'_0\setminus\{\Si=0\}$, these
formulas hold on $\Cal D'_0$, and the singularity at $\Si=0$ is
removable. As $\Cal D'_0\supset\Dadm'$, the proof of Theorem 8.7 is
complete. \qed
\enddemo

\example{Remark 8.9} The reader might have noticed that the proof of
Proposition 8.1 that we gave above is a verification that our
formulas give the
right answer rather than a derivation of these formulas. Unfortunately, at
this point we cannot suggest any derivation procedure for which we would
be able to prove that it produces the formulas we want. However, we can
explain how we obtained the answer.

Recall that in our treatment of the orthogonal polynomials in \S7 the crucial
role was played by the difference equation \tht{7.5}. Corollary 5.9 shows that
if we know a difference equation for the orthogonal polynomials then we can
write down difference equations for the values of $\Rin^{(N)}$ and
$\Sin^{(N)}$ on the lattice $\xin$. At this point we make the assumption that
the {\it meromorphic functions\/} $\Rin^{(N)}(\ze)$ and $\Sin^{(N)}(\ze)$
satisfy the same difference equations. {\it A priori\/}, it is not
clear at all
why this should be the case. However, the general philosophy of the
Riemann--Hilbert problem suggests that $\Rin^{(N)}(\ze)$ and $\Sin^{(N)}(\ze)$
should satisfy some difference equations. So we proceed and find meromorphic
solutions of the difference equation which we got from the lattice.
This can be
done using general methods of solving difference equations with polynomial
coefficients, see \cite{MT, Chapter XV}. We also want our
solutions to be holomorphic in $\C\setminus \xout$ and to have a fixed
asymptotics at infinity, and this leads (through heavy computations!) to the
formulas of Lemma 8.2. Unfortunately, these formulas are not suitable for the
limit transition $N\to\infty$. So we had to play around with the transformation
formulas for ${}_3F_2$ to get the convenient formulas of Proposition 8.1.
\endexample

\head 9. The spectral measures and continuous point processes
\endhead

Define a {\it continuous phase space}
$$
\x=\R\setminus\{\pm\tfrac 12\}
$$
and divide it into two parts
$$
\gather
\x=\xin\sqcup\xout,\\
\xin=(-\tfrac 12,\tfrac 12),\quad\xout=
(-\infty,-\tfrac 12)\sqcup(\tfrac 12,+\infty).
\endgather
$$

As in \S4, by $\Conf(\x)$ we denote the space of configurations in
$\x$. We define a map $\iota:\Om\to\Conf(\x)$ by
$$
\multline
\om=(\al^+,\be^+;\al^-,\be^-;\delta^+,\delta^-) \\
\overset{\iota}\to{\longmapsto}
\Cal C(\om)=\{\al_i^++\tfrac{1}2\}\sqcup\{\tfrac{1}2-
\be_i^+\}\sqcup\{-\al_j^--\tfrac{1}2\}\sqcup\{-
\tfrac{1}2+\be_j^-\},
\endmultline \tag9.1
$$
where we omit possible 0's in $\al^+,\be^+,\al^-,\be^-$, and also
omit possible 1's in $\be^+$ or $\be^-$.

\proclaim{Proposition 9.1} The map \tht{9.1} is a well--defined Borel
map.
\endproclaim

\demo{Proof} We have to prove that, for any compact set $A\subset\x$,
the function $\N_A(\iota(\om))$ takes finite values and is a Borel
function, see the definitions in the beginning of \S4. Without loss
of generality, we may assume that $A$ is a closed interval, entirely
contained in $(-\tfrac12,\tfrac12)$, $(\tfrac12,+\infty)$ or
$(-\infty,-\tfrac12)$.

Assume that $A=[x,y]\subset(-\tfrac12,\tfrac12)$. Then
$$
\gather
\N(\iota(\om))=\operatorname{Card}
\{i\mid\tfrac12-y\le\be^+_i\le\tfrac12-x\}+
\operatorname{Card}
\{j\mid\tfrac12+x\le\be^-_j\le\tfrac12+y\}\\
=\sum_{i=1}^\infty{\bold1}_{[1/2-y,1/2-x]}(\be^+_i)+
\sum_{j=1}^\infty{\bold1}_{[1/2+x,1/2+y]}(\be^-_j),
\endgather
$$
where ${\bold1}_{[\dots]}$ stands for the indicator of an interval.
Since $\sum\be^+_i\le\de^+$, $\sum\be^-_j\le\de^-$, and
$\tfrac12+x>0$, $\tfrac12-y>0$, the sums of
the indicator functions above are actually finite. Clearly, these are
Borel functions in $\om$.

When $A$ is contained in $(\tfrac12,+\infty)$ or
$(-\infty,\tfrac12)$, the argument is the same. \qed
\enddemo

Let $\chi$ be a character of the group $U(\infty)$ and let $P$ be its
spectral measure. By Proposition 9.1, the pushforward $\iota(P)$ of
$P$ is a well--defined probability measure on the space $\Conf(\x)$. We view
it as a point process and denote it by $\Cal P$. Our purpose  is to
describe the spectral measure $P$ (for concrete characters $\chi$) in terms
of the process $\Cal P$. The map $\iota$ glues some points $\om$
together, and it is natural to ask whether we lose any information
about $P$ by passing to $\Cal P$. We discuss this
issue at the end of the section.

Let $\{P_N\}$ be the coherent system corresponding to $\chi$.
Recall that in \S2(c) we have associated with each $P_N$ a probability
measure $\un P_N$ on $\Om$, which is the pushforward of $P_N$ under
an embedding of $\GT_N$ into $\Om$. According to Theorem 2.2, the
measures $\un P_N$ weakly converge to $P$ as $N\to\infty$. Let us
form the probability measures
$$
\un{\Cal P}^{(N)}=\iota(\un P_N),
$$
which are point processes on $\x$. Theorem 2.2 suggests the idea that
the process $\Cal P$ should be a limit of the processes
$\un{\Cal P}^{(N)}$ as $N\to\infty$. For instance, the correlation
measures of $\Cal P$ should be limits of the respective correlation
measures of $\un{\Cal P}^{(N)}$'s. Then this would give us a
possibility to find the correlation measures of $\Cal P$ through the
limit transition in the correlation measures of the processes
$\un{\Cal P}^{(N)}$ as $N\to\infty$. The goal of this section is to
provide a general justification of such a limit transition.

First of all, let us give a slightly different (but equivalent)
definition of the processes $\un{\Cal P}^{(N)}$. Recall that in \S4,
we have associated to each $P_N$ a point process $\Cal P^\n$ on the
lattice $\x^\n$. Consider the map
$$
\x^\n \longrightarrow \un{\x}^\n=\tfrac1N\,\x^\n\subset \x,
\qquad x\mapsto\dfrac xN\, . \tag9.2
$$
Then the process $\un{\Cal P}^{(N)}$ can be identified with the
pushforward of the process $\Cal P^\n$ under this map of the phase
spaces. Hence, denoting by $\rho^\n_k$ and $\un\rho^\n_k$ the $k$th
correlation measure of $\Cal P^\n$ and $\un{\Cal P}^{(N)}$,
respectively, we see that $\un\rho^\n_k$ is the pushforward of
$\rho^\n_k$ under the map \tht{9.2}.

Let us assume that $\Cal P^\n$ is a determinantal process with
a kernel $K^\n(x,y)$ on $\x^\n\times\x^\n$ (which actually holds in
our concrete situation). Then $\un{\Cal P}^{(N)}$ is a determinantal
process, too. It is convenient to take as the reference measure on the
lattice $\un\x^\n$ not the counting measure but the measure
$\un\mu^\n$ such that $\un\mu^\n(\{x\})=\tfrac1N$ for any
$x\in\un\x^\n$. The reason is that, as $N\to\infty$, the measures
$\un\mu^\n$ approach the Lebesgue measure on $\x$. Taking account of
the factor $\tfrac1N$ we see that the kernel
$$
\un K^\n(x,y)=N\cdot K^\n(Nx,Ny), \qquad x,y\in\x^\n.\tag9.3
$$
is a correlation kernel for $\un{\Cal P}^{(N)}$.

We need one more piece of notation: given $x\in\x$, let $x_N$ be the
node of the lattice $\x^\n$ which is closest to $Nx$ (any of two if
$Nx$ fits exactly at the middle between two nodes).

The main result of this section is as follows.

\proclaim{Theorem 9.2} Let $\chi$ be a character of $U(\infty)$ and
let $P$, $\{P_N\}$, $\x$, $\Cal P$, $\Cal P^\n$ be as
above. Assume that each $\Cal P^\n$ is a determinantal process
$\x^\n$ with a correlation kernel $K^\n(x,y)$ whose restriction both
to $\xin^\n\times\xin^\n$ and $\xout^\n\times\xout^\n$ is Hermitian
symmetric. Further, assume that
$$
\lim_{N\to\infty}N\cdot K^\n(x_N,y_N)=K(x,y), \qquad
x,y\in\x,
$$
uniformly on compact subsets of $\x\times\x$, where
$K(x,y)$ is a continuous function on $\x\times\x$.

Then $\Cal P$ is a determinantal point process and
$K(x,y)$ is its correlation kernel.
\endproclaim

The proof will be given after a preparation work. First, we review
a few necessary definitions and facts from \cite{Ol3}.

A {\it path\/} in the Gelfand--Tsetlin graph $\GT$
is an infinite sequence $t=(t_1,t_2,\dots)$ such that $t_N\in\GT_N$ and
$t_N\prec t_{N+1}$ for any $N=1,2,\dots$. The set of the paths will
be denoted by $\Cal T$.

Consider the natural embedding $\Cal T\subset\prod\limits_N \GT_N$.
We equip $\prod\limits_N \GT_N$ with the product topology (the sets
$\GT_N$ are viewed as discrete spaces). The set $\Cal T$ is
closed in this product space. We equip $\Cal T$ with the induced
topology. Then $\Cal T$ turns into a totally disconnected topological space.
Let $\tau=(\tau_1,\dots,\tau_N)$ be an arbitrary {\it finite\/} path in the
graph $\GT$, i.e., $\tau_1\in\GT_1, \dots, \tau_N\in\GT_N$ and
$\tau_1\prec\dots\prec\tau_N$. The cylinder sets of the form
$$
C_\tau=\{t\in\Cal T\mid t_1=\tau_1,\dots,t_N=\tau_N\},
$$
form a base of topology in $\Cal T$.

Consider an arbitrary signature $\la\in\GT_N$. The set of finite paths
$\tau=(\tau_1\prec\dots\prec\tau_N)$ ending at $\la$ has cardinality
equal to $\Dim_N\la=\chi^\la(e)$. The cylinder sets $C_\tau$
corresponding to these finite paths $\tau$ are pairwise disjoint, and
their union coincides with the set of infinite paths $t$ passing
through $\la$.

A {\it central measure\/} is any probability Borel measure on $\Cal
T$ such that the mass of any cylinder set $C_\tau$ depends only of
its endpoint $\la$. These definitions are inspired by \cite{VK1}.

There exists a natural bijective correspondence $M\longleftrightarrow\{P_N\}$
between central measures $M$ and coherent systems $\{P_N\}$, defined
by the relations
$$
\Dim_N\la\cdot M(C_\tau)=P_N(\la),
$$
where $N=1,2,\dots$, $\la\in\GT_N$, and $\tau$ is an arbitrary finite
path ending at $\la$. In other words, for any $N$, the
pushforward of $M$ under the natural projection
$$
\prod_{N=1}^\infty\GT_N\supset \Cal T\to \GT_N \tag9.4
$$
coincides with $P_N$.

By Proposition 2.1 we get a bijective correspondence between central
measures and characters of $U(\infty)$. This correspondence is an
isomorphism of convex sets. So, extreme central measures exactly
correspond to extreme characters.

On the other hand, by virtue of Theorem 1.2, we get a bijection
$M\longleftrightarrow P$ between central measures on $\Cal T$ and
probability measures on $\Om$. In more detail, the correspondence
$M\to P$ has the form
$$
M\to\{P_N\}\to\chi\to P.
$$

Given a path $t$, let $\wt p_i^{\,\pm}(N,t)$
and $\wt q_i^{\,\pm}(N,t)$ denote the modified Frobenius coordinates
of the Young diagram $(t_N)^\pm$. We say that $t$ is a {\it
regular path\/} if there exist limits
$$
\lim_{N\to\infty}\frac{\wt p_i^{\,\pm}(N,t)}N=\al_i^\pm, \quad
\lim_{N\to\infty}\frac{\wt q_i^{\,\pm}(N,t)}N=\be_i^\pm, \quad
\lim_{N\to\infty}\frac{|(t_N))^\pm|}N=\de^\pm,
$$
where $i=1,2,\dots$\,. Then the limit values are the coordinates of a
point $\om\in\Om$, and we
say that $\om$ is the {\it end\/} the regular path $t$ or that $t$
{\it ends\/} at $\om$.

Let $\Treg\subset\Cal T$ be the subset of regular paths. This is a
Borel set. Let $\Treg\to\Om$ be the projection assigning to a regular
path its end. This is a Borel map.

\proclaim{Theorem 9.3} Let $M$ be a central measure on $\Cal T$. Then
$M$ is supported by the Borel set $\Treg$ and hence can be viewed as
a probability measure on $\Treg$. The pushforward of $M$ under
the projection $\Treg\to\Om$ coincides with the spectral measure $P$
that appears in the correspondence $M\to P$ defined above.
\endproclaim

\demo{Proof} See \cite{Ol3, Theorem 10.7}. \qed
\enddemo

\proclaim{Corollary 9.4} Both the measure $P$ and all the measures
$P_N$ can be represented as pushforwards of the measure $M$, with
respect to the maps $\Treg\to\Om$ and $\Treg\to\GT_N$, respectively,
where the latter map is given by restricting \tht{9.4} to
$\Treg\subset\Cal T$.
\endproclaim

\demo{Proof} The claim concerning $P$ is exactly Theorem 9.3. The
claim concerning $P_N$ follows from the discussion above. \qed
\enddemo

For any compact set $A\subset\x$, let $\N_{A,N}$ denote the random
variable $\N_A$ associated with the process $\un{\Cal P}^{(N)}$.
Here it is convenient to consider as the phase space of
$\un{\Cal P}^{(N)}$ not the lattice $\un\x^\n$ but the ambient
continuous space $\x$. Recall that by $\un\rho_k^{(N)}$ we have
denoted the $k$th correlation measure of $\un{\Cal P}^{(N)}$. The
first step towards Theorem 9.2 is the following

\proclaim{Proposition 9.5} Assume that for any compact set
$A\subset\x$ there exist uniform on $N$ bounds for the moments of
$\N_{A,N}$,
$$
\E[(\N_{A,N})^l]\le C_l, \qquad l=1,2,\dots,
$$
where the symbol $\E$ means expectation.

Then for any $k=1,2,\dots$, the $k$th correlation measure $\rho_k$ of
the process $\Cal P$ exists. Moreover, for any continuous compactly
supported function $F$ on $\x^k=\x\times\dots\times\x$,
$$
\lim_{N\to\infty}\langle F,\un\rho^{(N)}_k\rangle=
\langle F,\rho_k\rangle.
$$
\endproclaim

\demo{Proof} We argue as in the proof of \cite{BO4, Lemma
6.2}. First, using Corollary 9.4, we put all the processes on one and
the same probability space, $(\Treg, M)$.

Given a regular path $t=(t_N)_{N=1,2,\dots}$, we assign to it point
configurations
$$
\Cal C(t)\in\Conf(\x), \quad \Cal C_N(t)\in\Conf(\x), \qquad N=1,2,\dots,
$$
as follows.

Let
$$
\om(t)=(\al^+(t), \be^+(t);\al^-(t),\be^-(t);\de^+(t),\de^-(t))\in\Om
$$
be the end of $t$. Then we set
$$
\Cal C(t)=\{\al_i^+(t)+\tfrac{1}2\}
\sqcup\{\tfrac{1}2-\be_i^+(t)\}
\sqcup\{-\al_j^-(t)-\tfrac{1}2\}
\sqcup\{-\tfrac{1}2+\be_j^-(t)\},
$$
cf. \tht{9.1}. Here, as in \tht{9.1}, we omit possible 0's in
$\al^+(t)$, $\be^+(t)$, $\al^-(t)$, $\be^-(t)$, and also possible 1's
in $\be^+(t)$ or $\be^-(t)$. Equivalently,
$$
\Cal C(t)=\Cal C(\om(t))=\iota(\om(t)).
$$

Likewise, for any $N=1,2,\dots$, let $\wt p^{\,\pm}_i(N,t)$ and
$\wt q^{\,\pm}_i(N,t)$ denote the modified Frobenius coordinates of
$(t_N)^\pm$, and let
$$
a^\pm_i(N,t)=\frac{\wt p^{\,\pm}_i(N,t)}N\,, \quad
b^\pm_i(N,t)=\frac{\wt q^{\,\pm}_i(N,t)}N\,, \qquad
i=1,2,\dots\,.
$$
We set
$$
\Cal C_N(t)=\{a_i^+(N,t)+\tfrac{1}2\}
\sqcup\{\tfrac{1}2-b_i^+(N,t)\}
\sqcup\{-a_j^-(N,t)-\tfrac{1}2\}
\sqcup\{-\tfrac{1}2+b_j^-(N,t)\}.
$$
Equivalently, $\Cal C_N(t)$ is the image of $t_N$
under the composite map $\GT_N\to\Om\to\Conf(\x)$.

We view $\Cal C(t)$ and $\Cal C_N(t)$ (for any $N=1,2,\dots$) as
random configurations defined on the common probability space
$(\Treg,M)$. By Corollary 9.4, these are
exactly the same as the random configurations corresponding to the
point processes $\Cal P$ and $\un{\Cal P}^{(N)}$, respectively.

{}From now on all the random variables will be referred to the
probability space $(\Treg,M)$. Fix a continuous compactly supported
function $F$ on $\X^k$. It will be convenient to assume that $F$ is
nonnegative (this does not mean any loss of generality). Introduce
random variables $f$ and $f_N$ as follows:
$$
f(t)=\sum_{x_1,\dots,x_k\in \Cal C(t)}F(x_1,\dots,x_k), \qquad
f_N(t)=\sum_{x_1,\dots,x_k\in \Cal C_N(t)}F(x_1,\dots,x_k), \tag9.5
$$
summed over ordered $k$-tuples of points with pairwise distinct
labels. Any such sum is actually finite because $F$ is compactly
supported and the point configurations are locally finite.

By the definition of the correlation measures,
$$
\langle F, \rho_k\rangle=\E[f], \qquad
\langle F, \un\rho_k^{(N)}\rangle=\E[f_N],
$$
and the very existence of $\rho_k$ is guaranteed if $\E[f]<\infty$
for any $F$ as above.
So, we have to prove that $\E[f_N]\to\E[f]<\infty$ as $N\to\infty$. By a
general theorem (see \cite{Sh, ch. II, \S6, Theorem 4}), it
suffices to check the following two conditions:

{\it Condition 1.\/} $f_N(t)\to f(t)$ for any $t\in\Treg$.

{\it Condition 2.\/} The random variables $f_N$ are uniformly
integrable, that is,
$$
\sup_N\int_{\{t\mid f_N(t)\ge c\}}f_N(t)M(dt)\to 0, \quad
\text{as $c\to+\infty$.}
$$

Let us check Condition 1. This condition does not depend on $M$, it
is a simple consequence of the regularity property. Indeed, for
$\ep>0$ let $\X_\ep$ be obtained from $\X$ be removing the
$\ep$-neighborhoods of the points $\pm\frac12$,
$$
\X_\ep=\X\setminus((-\tfrac12-\ep,-\tfrac12+\ep)\cup
(\tfrac12-\ep,\tfrac12+\ep)).
$$
Choose $\ep$ so small that the function $F$ be supported by
$\X_\ep^k$. Fix $l$ so large that $\al^\pm_j<\ep$, $\be^\pm_j<\ep$
for all indices $j\ge l$. By the definition of $\X_\ep$, if a point
$$
x=\cases \tfrac12+\al^+_i(t)\\ \tfrac12-\be^+_i(t)\\
-\tfrac12-\al^-_i(t)\\-\tfrac12+\be^-_i(t) \endcases \tag9.6
$$
lies in $\X_\ep$ then $i<l$.

By the definition of regular paths, for any index $i$,
$$
a^\pm_i(N,t)\to \al^\pm_i(t), \quad
b^\pm_i(N,t)\to\be^\pm_i(t), \qquad N\to\infty \tag9.7
$$
Therefore, we have $a^\pm_l(N,t)<\ep$
and $b^\pm_l(N,t)<\ep$ for all $N$ large enough. By monotonicity, the
same inequality holds for the indices $l+1,l+2,\dots$ as well. This
means that if $N$ is large enough and a point
$$
x=\cases \tfrac12+a^+_i(N,t)\\ \tfrac12-b^+_i(N,t)\\
-\tfrac12-a^-_i(N,t)\\-\tfrac12+b^-_i(N,t) \endcases \tag9.8
$$
lies in $\X_\ep$ then $i<l$.

It follows that in the sums \tht{9.5}, only the points \tht{9.6} or
\tht{9.8} with indices $i=1,\dots,l-1$ may really contribute. By
\tht{9.7} and continuity of $F$ we conclude that $f_N(t)\to f(t)$.

Let us check Condition 2. Choose a compact set $A$ such that $F$ is
supported by $A^k$. Denoting by $C$ the supremum norm of $F$, we have
the bound
$$
f_N(t)\le C\cdot N_{A, N}(t)(\N_{A,N}(t)-1)
\dots(\N_{A,N}(t)-k+1)
\le C\cdot(\N_{A,N}(t))^k.
$$
Therefore the random variables $f_N$ are uniformly integrable
provided that this is true for the random variables
$(\N_{A,N})^k$ for any fixed $k$. But the latter fact follows from the
assumption of the proposition and Chebyshev's inequality. \qed
\enddemo

To apply Proposition 9.5 we must check the required uniform bound for
the moments. By the assumption of Theorem 9.2,
each point process $\un{\Cal P}_N$ is a
determinantal process on $\un\x^n$ such that its kernel, restricted
to $\un\x_{\inr}^\n\times\un\x_{\inr}^\n$ and to
$\un\x_{\out}^\n\times\un\x_{\out}^n$, is Hermitian symmetric. Here we set
$\un\x_{\inr}^\n=\un\x^\n\cap\x_{\inr}$,
$\un\x_{\out}^\n=\un\x^\n\cap\x_{\out}$. For a
compact set $A\subset\X$ we denote by $\un K^\n_A$ the restriction of
the kernel to $A\cap\un\x^\n$. If $A$ is entirely contained in $\xin$
or $\xout$ then $\un K^\n_A$ is a finite--dimensional {\it
nonnegative\/} operator.

\proclaim{Proposition 9.6} Assume that for any compact set $A$, which
is entirely contained in $\xin$ or $\xout$, we have a bound of the
form $\tr \un K^\n_A\le \const$, where the constant does not depend on
$N$. Then the assumption of Proposition 9.5 is satisfied.
\endproclaim

\demo{Proof} Let $A\subset\X$ be an arbitrary compact set. Then
$A=\ain\cup\aout$, where $\ain\subset\xin$ and
$\aout\subset\xout$ are compact sets. We have
$\N_{A,N}=\N_{\aout,N}+\N_{\ain,N}$. If we know uniform bounds for the
moments of $\N_{\aout,N}$ and $\N_{\ain,N}$ then, using the
Cauchy-Schwarz--Bunyakovskii inequality, we get such bounds for
$\N_{A,N}$, too. Consequently, we can assume that $A$ is entirely contained
in $\xout$ or $\xin$, so that $\un K^\n_A$ is nonnegative.

Instead of ordinary moments we can deal with factorial moments. Given
$l=1,2,\dots$, the $l$-th factorial moment of $\N_{A,N}$ is equal to
$$
\un\rho_l^{(N)}(A^l)=\int_{A^l}\det [\un K^\n_A(x_i,x_j)]_{1\le i,j\le l}\,
dx_1\dots dx_l=l!\tr(\wedge^l \un K^\n_A).
$$
Since $\un K^\n_A$ is nonnegative, we have
$$
\tr(\wedge^l \un K^\n_A)\le\tr(\otimes^l \un K^\n_A)=(\tr(\un K^\n_A))^l.
$$
This concludes the proof, because we have a uniform bound for the
traces by assumption. \qed
\enddemo

\demo{Proof of Theorem 9.2} We shall approximate the process $\Cal P$
by the processes $\un{\Cal P}^\n$. Recall that the process $\un{\Cal
P}^\n$ is a scaled version of the process $\Cal P^\n$, and their
correlation kernels, $\un K^\n$ and $K^\n$, are related as follows
$$
\un K^\n(x,y)=N\cdot K^\n(Nx,Ny), \qquad x,y\in\un\x^\n=\tfrac1N\x.
$$

Let us check that the assumption of Proposition 9.6 is satisfied.
Without loss of generality we may assume that $A$ is a closed
interval $[a,b]$, contained either in $\xin$ or $\xout$. We have
$$
\tr \un K^\n_A=\frac1N\, \sum_{x\in A\cap\un\x^n}\un K^\n(x,x),
$$
where the factor $\tfrac1N$ comes from the reference measure $\mu^\n$
on the lattice $\x^\n$ (recall that $\mu^\n$ assigns weight $\tfrac1N$ to
each node). By the relation between the two kernels, this can be
rewritten as follows
$$
\tr \un K^\n_A=\frac1N\, \sum_{x\in A\cap\un\x^n}N\cdot K^\n(Nx,Nx).
$$
By the assumption of Theorem 9.2, the right--hand side tends, as
$N\to\infty$, to the integral $\int_a^b K(x,x)dx$, hence
the traces above are uniformly bounded, as required.

Proposition 9.6 makes it possible to apply Proposition 9.5. It
shows that the process $\Cal P$ possesses correlations measures
$\rho_k$, which are weak limits of the measures $\un{\rho}^\n_k$.

On the other hand, we know that for any $k$, the correlation measure
$\un{\rho}^\n_k$ is expressed, by a determinantal formula, through
the correlation kernel $\un K^\n$. Further, by the assumption of Theorem
9.2, this kernel tends to the kernel $K$ as $N\to\infty$.
This implies that the limit measure $\rho_k$ is expressed, by the same
determinantal formula, through the limit kernel $K$. \qed
\enddemo

Now we return to the discussion of the correspondence
$P\mapsto\Cal P$. Recall that it is based on the map
$\iota:\Om\to\Conf(\x)$, sending a point
$\om=(\al^+,\be^+;\al^-,\be^-;\de^+,\de^-)$ to a configuration
$\Cal C(\om)$. This is a continuous analog of the map
$\GT_N\to\Conf(\x^\n)$ sending a signature $\la$ to a point configuration
$X(\la)$, see \tht{4.1}. The correspondence $\la\mapsto X(\la)$ is reversible
while the map $\om\mapsto\Cal C(\om)$ is not. This is caused by three factors
listed below.

(i) The coordinates $\be^+_i$ and $\be^-_j$ become indistinguishable:
given a point $x\in\Cal C(\om)\cap\xin$, there is no way to decide
whether it comes from a coordinate $\tfrac12-x$ of $\be^+$ or from a
coordinate $x+\tfrac12$ of $\be^-$. Note that in the discrete case
such a problem does not arise. Indeed, if $d^+$ and $d^-$ stand for
the numbers of points $x\in X(\la)$ that are on the right and on the
left of $\xin$, respectively, then $X(\la)$ has exactly $d^++d^-$
points in $\xin$ of which the leftmost $d^-$ points come from $\la^-$
while the remaining $d^+$ points come from $\la^+$. But in the
continuous case, such an argument fails, because the total number of
points is, generally speaking, infinite.

(ii) The map $\iota$ ignores the coordinates $\de^\pm$.

(iii) The map $\iota$ ignores possible 1's in $\be^\pm$.

Let us discuss the significance of these factors in succession.

{\it Factor\/} (i). Remark that exactly the same effect
of mixture of the plus and minus $\be$--coordinates arises when an
extreme character $\chi^{(\om)}$ is restricted from the group $U(\infty)$
to the subgroup $SU(\infty)$, see \cite{Ol3, Remark 1.7}. Hence, if one
agrees to view characters that coincide on $SU(\infty)$ as equivalent
ones, then the factor in question becomes not too important.

{\it Factor\/} (ii). We conjecture that in our concrete situation (i.e., for
the characters $\chi_\zw$) the spectral measure $P$ is concentrated
on the subset of $\om$'s with $\ga^\pm=0$ (see the definition of
$\ga^\pm$ in \S1). If this is true then $\de^\pm$ is almost surely
equal to $\sum(\al^\pm_i+\be^\pm_i)$. The conjecture is supported by
the fact that the vanishing of the gamma parameters was proved in
similar situations, see \cite{P.I, Theorem 6.1} and \cite{BO4,
Theorem 7.3}. The method of \cite{BO4} makes it possible to
reduce the conjecture to the following claim concerning the first
correlation measures of the processes $\un{\Cal P}^\n$:
$$
\lim_{\varepsilon\to0}\int_{\frac12-\varepsilon}^{\frac12+\varepsilon}
|x-\tfrac12|\,\un{\rho}^\n_{\,1}(dx)=0, \qquad
\lim_{\varepsilon\to0}\int_{-\frac12-\varepsilon}^{-\frac12+\varepsilon}
|x+\tfrac12|\,\un{\rho}^\n_{\,1}(dx)=0
$$
uniformly on $N$.

{\it Factor\/} (iii). Again, we remark that possible 1's in $\be^\pm$
play no role when characters are restricted to the subgroup
$SU(\infty)$. On the other hand, we can prove that our processes, one
can speculate that in concrete the 1's do not appear almost surely.

The above arguments (though not rigorous) present a justification of
the passage $P\mapsto\Cal P$.

We conclude the section with one more general result which will be
used in \S10.

Observe that the set of characters of $U(\infty)$ is stable under the
operation of pointwise multiplication by $\det(\,\cdot\,)^k$, where
$k\in\Z$.

\proclaim{Proposition 9.7} Let $\chi$ be a character of $U(\infty)$,
$P$ be its spectral measure on $\Om$, and $\Cal P$ be the
corresponding point process on $\x=\R\setminus\{\pm\tfrac12\}$. Then
$\Cal P$ does not change under $\chi\mapsto\chi\det(\,\cdot\,)^k$,
$k\in\Z$.
\endproclaim

\demo{Proof} It suffices to prove this for $k=1$. Assume first that
$\chi$ is extreme, so that $\chi=\chi^{(\om)}$, where
$\om=(\al^+,\be^+;\al^-,\be^-;\de^+,\de^-)\in\Om$, and $P$ reduces to
the Dirac mass at $\{\om\}$. {}From the
explicit expression for $\chi^{(\om)}$, see \tht{1.2}, it follows that
$\chi^{(\om)}\det(\,\cdot)$ is an extreme character, too. Moreover,
the corresponding element $\bar\om\in\Om$ looks as follows: the
parameters $\al^\pm$ and $\de^\pm$ do not change, while
$$
\be^+\mapsto(1-\be^-_1,\be^+_1,\be^+_2,\dots), \qquad
\be^-\mapsto(\be^-_2,\be^-_3,\dots).
$$

On the other hand, from the definition of the projection
$\om\mapsto\Cal C(\om)$, see \tht{9.1}, it follows that the change
$\om\mapsto\bar\om$ does not affect the configuration $\Cal C(\om)$.
This proves the needed claim for extreme $\chi$.

Now let $\chi$ be arbitrary. By the very definition of spectral
measures (see Theorem 1.2), the spectral measure of the character
$\chi\det(\,\cdot\,)$ coincides with the pushforward of the spectral
measure $P$ under the map $\om\mapsto\bar\om$ of $\Om$. We have just
seen that this map does not affect the projection $\om\mapsto\Cal
C(\om)$. Since $\Cal P$ is the image of $P$ under this projection, we
conclude that $\Cal P$ remains unchanged. \qed
\enddemo

\head 10. The correlation kernel of the process $\Cal P$
\endhead

Our goal in this section is to compute the correlation functions of the
process $\Cal P$ associated to the spectral measure for zw--measures,
see the beginning of \S9 for the definitions. Theorem 10.1 below is
the main result of the paper.

In the formulas below we use the {\it Gauss hypergeometric function}.
Recall that this is a function in one complex variable (say, $u$)
defined inside the unit circle by the series
$$
{}_2F_1\left[\matrix
a,\,b\\
c\endmatrix\,\Bigl|\,u\right]=\sum_{k\ge
0}^\infty\frac{(a)_k(b)_k}{k!(c)_k}\,u^k.
$$
Here $a,b,c$ are complex parameters, $c\notin \{0,-1,\dots\}$.

This function can be analytically continued to the domain $u\in\C\setminus
[1,+\infty)$, see, e.g., \cite{Er, Vol. 1, Chapter 2}. We will need the fact
that for any fixed $u$ in $\C\setminus[1,+\infty)$, the expression
$$
\frac 1{\Ga(c)}\,{}_2F_1\left[\matrix
a,\,b\\
c\endmatrix\,\Bigl|\,u\right]
$$
defines an entire function in $(a,b,c)\in \C^3$. This follows, e.g.,
from \cite{Er, 2.1.3(15)}.

\proclaim{Theorem 10.1} Let $(\zw)\in\Dadm$ and $\Cal P$ be the
corresponding point process on $\x=\R\setminus\{\pm\tfrac 12\}$
defined in \S9.

For any $n=1,2,\dots$ and $x_1,\dots,x_n\in\x$, the
$n$th correlation function $\rho_n(x_1,\dots,x_n)$ of the process
$\Cal P$ has the determinantal form
$$
\rho_n(x_1,\dots,x_n)=\det[K(x_i,x_j)]_{i,j=1}^n.
$$
The kernel $K(x,y)$ with respect to the splitting
$\x=\xout\sqcup \xin$ has the following form
$$
\gathered
K_{\out,\out}(x,y)=\sqrt{\psiout(x)\psiout(y)}
\,\frac{\Rout(x)\Sout(y)-\Sout(x)\Rout(y)}{x-y}\,,\\
K_{\out,\inr}(x,y)=\sqrt{\psiout(x)\psiin(y)}
\,\frac{\Rout(x)\Rin(y)-\Sout(x)\Sin(y)}{x-y}\,,\\
K_{\inr,\out}(x,y)=\sqrt{\psiin(x)\psiout(y)}
\,\frac{\Rin(x)\Rout(y)-\Sin(x)\Sout(y)}{x-y}\,,\\
K_{\inr,\inr}(x,y)=\sqrt{\psiin(x)\psiin(y)}
\,\frac{\Rin(x)\Sin(y)-\Sin(x)\Rin(y)}{x-y}\,,
\endgathered
$$
where
$$
\gathered
\psiout(x)=\cases C(z,z')\cdot\left(x-\frac 12\right)^{-z-
z'}\left(x+\frac 12\right)^{-w-w'},& x>\frac 12\,,\\
C(w,w')\cdot\left(-x-\frac 12\right)^{-w-w'}\left(-x+\frac
12\right)^{-z-z'}, &
x<-\frac 12\,,
\endcases
\\
\psiin(x)=\left(\tfrac 12-x\right)^{z+z'}\left(\tfrac
12+x\right)^{w+w'},\quad -\tfrac 12<x<\tfrac 12,\\
C(z,z')=\frac{\sin\pi z\sin\pi z'}{\pi^2}\,,\quad C(w,w')=\frac{\sin\pi
w\sin\pi w'}{\pi^2}\,,
\endgathered
$$
and
$$
\aligned
\Rout(x)&=\left(\frac{x+\frac 12}{x-\frac 12}\right)^{w'}\,
{}_2F_1\left[\matrix
z+w',\,z'+w'\\ z+z'+w+w'\endmatrix\,\Biggl|\,\frac 1{\frac 12 -x}\right]\,,\\
\Sout(x)&=\Gamma\left[\matrix z+w+1,\,
z+w'+1,\,z'+w+1,\,z'+w'+1\\z+z'+w+w'+1,\,z+z'+w+w'+2\endmatrix\right]\\
&\times
\frac 1{x-\frac 12}\,\left(\frac{x+\frac 12}{x-\frac
12}\right)^{w'}\,{}_2F_1\left[\matrix z+w'+1,\,z'+w'+1\\
z+z'+w+w'+2\endmatrix\,\Biggl|\,\frac 1{\frac 12 -x}\right]\,,
\endaligned
$$
$$
\aligned
\Rin(x)=&-\frac {\sin \pi z}{\pi}\,\Gamma\left[\matrix z'-
z,\,z+w+1,\,z+w'+1\\z+w+z'+w'+1\endmatrix\right]\\ &\times
\left(\frac 12 +x\right)^{-w}\left(\frac 12 -x\right)^{-z'}
{}_2F_1\left[\matrix
z+w'+1,\, -z'-w\\z-z'+1\endmatrix\,\Biggl|\,\frac 12 -x\right]\\ &-
\frac {\sin \pi z'}{\pi}\,\Gamma\left[\matrix z-
z',\,z'+w+1,\,z'+w'+1\\z+w+z'+w'+1\endmatrix\right]\\ &\times
\left(\frac 12 +x\right)^{-w}\left(\frac 12 -x\right)^{-z}
{}_2F_1\left[\matrix
z'+w'+1,\, -z-w\\z'-z+1\endmatrix\,\Biggl|\,\frac 12 -x\right]\,,
\endaligned
$$
$$
\aligned
\Sin(x)=&-\frac{\sin\pi z}{\pi}\, \Gamma\left[\matrix z'-z,\,z+z'+w+w'\\
z'+w,\,z'+w'\endmatrix \right]\\ &\times
\left(\frac 12 +x\right)^{-w}\left(\frac 12 -x\right)^{-z'}
{}_2F_1\left[\matrix
z+w',\, -z'-w+1\\z-z'+1\endmatrix\,\Biggl|\,\frac 12 -x\right]\\ &-
\frac{\sin\pi z'}{\pi}\, \Gamma\left[\matrix z-z',\,z+z'+w+w'\\
z+w,\,z+w'\endmatrix \right]\\ &\times
\left(\frac 12 +x\right)^{-w}\left(\frac 12 -x\right)^{-z}
{}_2F_1\left[\matrix
z'+w',\, -z-w+1\\z'-z+1\endmatrix\,\Biggl|\,\frac 12 -x\right].
\endaligned
$$

The indeterminacy on the diagonal $x=y$ is resolved by the L'Hospital rule:
$$
\gathered
K_{\out,\out}(x,x)=\psiout(x)
\left(\Rout'(x)\Sout(x)-\Sout'(x)\Rout(x)\right),
\\
K_{\inr,\inr}(x,x)=\psiin(x)
\left(\Rin'(x)\Sin(x)-\Sin'(x)\Rin(x)\right).
\endgathered
$$
\endproclaim

\demo{Singularities} The formulas for the function $\Rout,\Sout,\Rin,\Sin$
above have no singularities for
$$
(z,z',w,w')\in\Cal D_0\setminus\left(\{\s=0\}\cup \{z-z'\in \Z\}\right).
$$
Moreover, we will prove that the value of the kernel $K(x,y)$ can be extended
to a continuous function on $\Cal \Dadm$ for any fixed $x,y\in \x$.
(Recall that the process $\Cal P$ is defined for $(z,z',w,w')\in\Dadm$.)
\enddemo

\demo{Vanishing of the kernel} Note that if $(z,z')\in\Zd$ (see \S3
for the definition of $\Zd$) then the
function $\psiout$ vanishes on $(\tfrac12,+\infty)$, because
$C(z,z')=0$. This implies that $K(x,y)=0$ whenever $x$ or $y$
is greater than $\tfrac12$. It follows that the configurations of the
process $\Cal P$ do not have points in $(\tfrac12,+\infty)$.

Likewise, if $(w,w')\in\Zd$ then the configurations of the process do
not intersect $(-\infty,-\tfrac12)$.
\enddemo

\demo{Proof} First of all observe that it suffices to prove the
theorem when $(\zw)\in\Dadm'$. Indeed, if
$(\zw)\in\Dadm\setminus\Dadm'$ then $(\zw)$ can be moved to $\Dadm'$
by an appropriate shift of the parameters, which is equivalent to
multiplying the initial character $\chi$ by $\det(\,\cdot\,)^k$ with
a certain $k\in\Z$, see Remarks 6.4 and 3.7. Next, according to
Proposition 9.7, multiplication by $\det(\,\cdot)^k$
does not affect the point process $\Cal P$.

To carry out the desired reduction we have to check that the formulas
for the functions $\rho_n$ given in Theorem 10.1 are also invariant
under any shift of the parameters of the form
$$
(\zw)\mapsto(z+k,z'+k,w-k,w'-k), \quad k\in\Z.
$$
Note that the kernel $K(x,y)$ is not invariant under such a shift.
To see what happens with the kernel we observe that
$$
\gather
\psiout(x)\to\left(\frac{x+\tfrac12}{x-\tfrac12}\right)^{2k}\psiout(x),
\quad
\psiin(x)\to\left(\frac{x+\tfrac12}{-x+\tfrac12}\right)^{-2k}\psiin(x),\\
\Rout(x)\to\left(\frac{x+\tfrac12}{x-\tfrac12}\right)^{-k}\Rout(x),
\quad
\Sout(x)\to\left(\frac{x+\tfrac12}{x-\tfrac12}\right)^{-k}\Sout(x),\\
\Rin(x)\to(-1)^k\,\left(\frac{x+\tfrac12}{-x+\tfrac12}\right)^{k}\Rin(x),
\quad
\Sin(x)\to(-1)^k\,\left(\frac{x+\tfrac12}{-x+\tfrac12}\right)^{k}\Sin(x).
\endgather
$$
It follows that
$$
K(x,y)\to\phi(x)K(x,y)(\phi(y))^{-1},
$$
where
$$
\phi(x)=\cases 1, & x\in\xout, \\ (-1)^k, & x\in\xin. \endcases
$$
But such a transformation of the kernel does not affect the
determinantal formula for the correlation functions.

{}From now on we will assume that $(\zw)\in\Dadm'$, as in Theorem 8.7.
At this moment we impose additional restrictions $\Si\ne0$ and
$z-z'\notin\Z$. We will need these restrictions in Propositions
10.3--10.4 below. After that they will be removed.

For any $x\in\x=\R\setminus\{\pm\frac12\}$, let $x_N$ denote the
point of the lattice $\x^\n=\Z+\tfrac{N-1}2$ which is closest to
$Nx$ (if there are two such points then we choose any of them). By
Theorems 8.7 and 9.2, it suffices to prove that 
$$
\lim_{N\to\infty}N\cdot K^\n(x_N,y_N)=K(x,y), \tag10.1
$$
uniformly on compact sets of $\x\times\x$. Here $K^\n$ is the kernel
of Theorem 8.7.

To do this we will establish the uniform convergence of all six
functions $\psiout^\n$, $\psiin^\n$, $\Rout^\n$, $\Sout^\n$,
$\Rin^\n$, $\Sin^\n$ of Theorem 8.7 to the respective functions of
Theorem 10.1. In order to overcome the difficulty arising from
vanishing of the denominators at $x=y$ we will establish the
convergence of $\Rout^\n$, $\Sout^\n$, $\Rin^\n$, $\Sin^\n$ in a
complex region containing $\x$.

The needed convergence \tht{10.1} follows from Propositions
10.2--10.4 below.

\proclaim{Proposition 10.2} There exist limits
$$
\lim_{N\to\infty}N^{\s}\psiout^{(N)}(x_N)=\psiout(x),
\qquad \lim_{N\to\infty}N^{-\s}\psiin^{(N)}(x_N)=\psiin(x),
$$
as $N\to\infty$, where the functions $\psiin^\n$ and $\psiout^\n$
were defined in \tht{6.2}, \tht{6.3}. The convergence is uniform on
the compact subsets of $\xout=\R\setminus [-\frac 12,\frac 12]$ and
$\xin=(-\frac 12,\frac 12)$, respectively.
\endproclaim

\proclaim{Proposition 10.3} Let $I$ be a compact subset of $\xout$. Set
$$
I_\varepsilon=\{\ze\in\C\,|\, \Re\ze\in I,\,|\Im \ze|<\varepsilon\}.
$$
Then for $\varepsilon>0$ small enough, for any $\ze\in I_\varepsilon$ we have
$$
\lim_{N\to\infty} \Rout^{(N)}(N\ze)=\Rout(\ze),\quad
\lim_{N\to\infty}N^{-\s}\Sout^{(N)}(N\ze)=\Sout(\ze),
\tag 10.2
$$
and the convergence is uniform on $I_\varepsilon$.
 \endproclaim

\proclaim{Proposition 10.4} Let $J$ be a compact subset of $\xin$. Set
$$
J_\varepsilon=\{\ze\in\C\,|\, \Re\ze\in J,\,|\Im \ze|<\varepsilon\}.
$$
Then for $\varepsilon>0$ small enough, for any $\ze\in J_\varepsilon$ we have
$$
\lim_{N\to\infty} \Rin^{(N)}(N\ze)=\Rin(\ze),\quad
\lim_{N\to\infty}N^{\s}\Sin^{(N)}(N\ze)=\Sin(\ze),
\tag 10.3
$$
and the convergence is uniform on $J_\varepsilon$.
\endproclaim

\demo{Proof of Proposition 10.2} Follows from the following uniform estimates.

For $x\in\xin$ we have
$$
\gathered
\frac{\Gamma(-x_N+u+\frac{N+1}2)}{\Gamma(-x_N+\frac{N+1}2)}=
N^u\,\left(-x+\frac 12\right)^u(1+O(N^{-1})),\quad u=z\text{  or  }z', \\
\frac{\Gamma(x_N+v+\frac{N+1}2)}{\Gamma(x_N+\frac{N+1}2)}=
N^v\,\left(x+\frac 12\right)^v(1+O(N^{-1})),\quad v=w\text{ or }w',
\endgathered
$$
as $N\to\infty$.

For $x\in\xout$ and $x>\frac 12$ we use the formulas
$$
\gathered
\frac 1{\Gamma(-x_N+u+\frac{N+1}2)}= \frac {\sin
\left(\pi(-x_N+u+\frac{N+1}2)\right)}{\pi}\,
\Gamma\left(x_N-u-\frac{N-1}2\right)\\=(-1)^{-x_N+\frac{N+1}2}\,\frac {\sin
(\pi u)}{\pi}\,\Gamma\left(x_N-u-\frac{N-1}2\right),\quad u=z\text{  or  }z',
\endgathered
$$
and the asymptotic relations
$$
\gathered
\frac{\Gamma(x_N-u-\frac{N-1}2)}{\Gamma(x_N-\frac{N-1}2)}= N^{-u}\,
\left(x-\frac 12\right)^{-u}(1+O(N^{-1})),\quad u=z\text{  or  }z',\\
\frac{\Gamma(x_N+\frac{N+1}2)}{\Gamma(x_N+v+\frac{N+1}2)}=
N^{-v}\left(x+\frac 12\right)^{-v}(1+O(N^{-1})),\quad v=w\text{  or  }w',
\endgathered
$$
as $N\to\infty$.

For $x\in\xout$ and $x<-\frac 12$ we use the formulas
$$
\gathered
\frac 1{\Gamma(x_N+v+\frac{N+1}2)}=\frac {\sin
\left(\pi(x_N+v+\frac{N+1}2)\right)}{\pi}\,
\Gamma\left(-x_N-v-\frac{N-1}2\right)\\=(-1)^{x_N+\frac{N+1}2}\,\frac {\sin
(\pi v)}{\pi}\,\Gamma\left(-x_N-v-\frac{N-1}2\right),\quad v=w\text{  or  }w',
\endgathered
$$
and the asymptotic relations
$$
\gathered
\frac{\Gamma(-x_N+\frac{N+1}2)}{\Gamma(-x_N+u+\frac{N+1}2)}= N^{-u}\,
\left(-x+\frac12\right)^{-u}(1+O(N^{-1})),\quad u=z\text{  or  }z',\\
\frac{\Gamma(-x_N-v-\frac{N-1}2)}{\Gamma(-x_N-\frac{N-1}2)}= N^{-v}\,
\left(-x-\frac 12\right)^{-v}(1+O(N^{-1})),\quad v=w\text{  or  }w',  
\endgathered $$
as $N\to\infty$. \qed
\enddemo

\demo{Proof of Proposition 10.3} We will employ the formulas \tht{8.1} and
\tht{8.2}. We have
$$
\Rout^{(N)}(N\ze)=
\Gamma\bmatrix
N\ze-\frac{N-1}2,\,N\ze+w'+\frac{N+1}2\\N\ze+\frac{N+1}2,\,N\ze+w'-\frac{N-
1}2\endbmatrix
{}_3F_2\left[\matrix -N,\,z+w',\,z'+w'\\\Sigma,\,N\ze+w'-\frac{N-
1}2\endmatrix\biggl|\,1\,\right].
$$

Assume $\Re \ze<-\frac 12$.
Handling the gamma factors is easy:
$$
\gathered
\Gamma\bmatrix
N\ze-\frac{N-1}2,\,N\ze+w'+\frac{N+1}2\\N\ze+\frac{N+1}2,\,N\ze+w'-\frac{N-
1}2\endbmatrix
=\frac{\Gamma(-N\ze-\frac{N-1}2)}{\Gamma(-N\ze-w'-\frac{N-1}2)}\cdot
\frac{\Gamma(-N\ze-w'+\frac{N+1}2)}{\Gamma(-N\ze+\frac{N+1}2)}\\ =
\left(-\ze-\frac 12\right)^{w'}\left(-\ze+\frac
12\right)^{-w'}(1+O(N^{-1}))=\left(\frac{\ze+\frac12}{\ze-\frac 
12}\right)^{w'}(1+O(N^{-1})),
\endgathered
$$
as $N\to\infty$, and the estimate is uniform on a small complex neighborhood
of any compact subset $I\subset(-\infty,-\frac12)$.

To complete the proof of the first limit relation \tht{10.1} we need to have
the equality
$$
\lim_{N\to\infty}{}_3F_2\left[\matrix
-N,\,z+w',\,z'+w'\\z+z'+w+w',\,N\ze+w'-\frac{N-
1}2\endmatrix\biggl|\,1\,\right]={}_2F_1\left[\matrix
z+w',\,z'+w'\\z+z'+w+w'\endmatrix\biggl|\frac 1{\frac 12-\ze}\right]
$$
with the uniform convergence.
This limit transition is justified by the following lemma.

\proclaim{Lemma 10.5} Let $A,B,C$ be complex numbers, $C\ne 0,-1,-2,\dots$,
$\{D_N\}_{N=1}^\infty$ be a sequence of complex numbers, $D_N\ne 0,-1,\dots$
for all $N$. Assume that
$$
\lim_{N\to\infty}\frac{-N}{D_N}=q\in(0,1).
$$
Then
$$
\lim_{N\to\infty}{}_3F_2\left[\matrix
-N,\,A,\,B\\C,\,D_N\endmatrix\,\biggl|\,1\right]={}_2F_1\left[\matrix
A,\,B\\C\endmatrix\,\biggl|\,q\right].
$$
The convergence is uniform on any set of sequences $\{D_N\}_{N=1}^\infty$ such
that $\{\frac N{D_N}-\lim\frac N{D_N}\}$ uniformly converges to 0 as
$N\to\infty$ and  $\lim\frac N{D_N}$ is uniformly bounded from $-1$
on this set.
\endproclaim

\demo{Proof}
We have
$$
{}_3F_2\left[\matrix -N,\,A,\,B\\C,\,D_N\endmatrix\,\biggl|\,1\right]=
\sum_{k=0}^N \frac{(-N)_k(A)_k(B)_k}{k!(C)_k(D_N)_k}\,.
$$
Let us show that these sums converge uniformly in $N>N_0$ for some $N_0>0$.
Thanks to our hypothesis, for large enough $N$ we can pick
$q_0\in (q,1)$ and a positive number $d_N$ such that for all our
sequences $\{D_N\}$, $-\Re D_N>d_N>Nq_0^{-1}$. Then (note that $k\le N$)
$$
|(D_N)_k|=|(-D_N)(-D_N-1)\ldots(-D_N-k+1)|\ge d_N(d_N-1)\cdots (d_N-k+1).
$$
Hence,
$$
\left|\frac{(-N)_k}{(D_N)_k}\right|\le
\frac{N(N-1)\cdots(N-k+1)}{d_N(d_N-1)\cdots (d_N-k+1)}\le \left(\frac
N{d_N}\right)^k<q_0^k.
$$
Thus, the sums above for large enough $N$ are
majorated by the convergent series
$$
\sum_{k=0}^\infty \left|\frac{(A)_k(B)_k}{k!(C)_k}\right|\, q_0^k
$$
and, therefore, converge uniformly. This means that to compute the limit as
$N\to\infty$ we can pass to the limit $N\to\infty$ in every term of the sum.
Since for any fixed $k$
$$
\lim_{N\to\infty}\frac{(-N)_k}{(D_N)_k}=q^k,
$$
this yields
$$
\lim_{N\to\infty}{}_3F_2\left[\matrix
-N,\,A,\,B\\C,\,D_N\endmatrix\,\biggl|\,1\right]=\sum_{k\ge
0}\frac{(A)_k(B)_k}{k!(C)_k}\, q^k={}_2F_1\left[\matrix
A,\,B\\C\endmatrix\,\biggl|\,q\,\right].
$$
The fact that we majorated the series by the same convergent series for all
our sequences $\{D_N\}$, and the uniform convergence of the terms of the
series guarantee the needed uniform convergence on the set of sequences.\qed
\enddemo

To prove the first limit relation for $I\subset(\frac 12,+\infty)$ we
just note that by uniqueness of monic orthogonal polynomials with a fixed
weight, $\Rout^{(N)}(x)$ is invariant with respect to the substitution
$$
x\mapsto -x, \qquad (z,z')\longleftrightarrow (w',w),
$$
cf. Lemmas 8.4, 8.5, and so is $\Rout(x)$, because of the transformation
formula $$
{}_2F_1\left[\matrix A,\, B\\C\endmatrix\,\Bigr|\, \ze\right]=(1-\ze)^{-A}
{}_2F_1\left[\matrix A,\, C-B\\C\endmatrix\,\Bigr|\, \frac\ze{\ze
-1}\right].
$$

The proof of the second relation \tht{10.1} is very similar. Note that both
$\Sout^{(N)}(x)$ and $\Sout(x)$ are skew--symmetric with respect to the
substitution above.

The proof of Proposition 10.3 is complete.\qed
\enddemo

\demo{Proof of Proposition 10.4} The argument is quite similar to the proof of
Proposition 10.3 above. Let us evaluate the asymptotics of the right--hand
side of \tht{8.4}. Let us look at the first term. Clearly, the argument for
the second term will be just the same. Gamma--factors give (here $\Re z\in
J\subset(-\frac 12,\frac 12)$)  $$
\multline
\Gamma\left[\matrix
N\ze+\frac{N+1}2,\,-N\ze+\frac{N+1}2,\,N+1+\s\\
-N\ze+z'+\frac{N+1}2,\,N\ze+w+\frac{N+1}2,N+1+z+w'
\endmatrix\right]
\\
=\left(\ze+\frac 12\right)^{-w}\left(\ze-\frac 12\right)^{-z'}(1+O(N^{-1})),
\endmultline
$$
as $N\to\infty$, uniformly on a neighborhood of $J$.
To complete the proof of the first relation \tht{10.2} we need to show that
$$
\multline
\lim_{N\to\infty}{}_3F_2\left[\matrix  z+w'+1,\,
-z'-w,\,-N\ze+z+\frac{N+1}2\\z-z'+1,\,N+1+z+w'\endmatrix\,\Biggl|\,1\right]
\\={}_2F_1\left[\matrix
z+w'+1,\, -z'-w\\z-z'+1\endmatrix\,\Biggl|\,\frac 12 -\ze\right]
\endmultline
$$
uniformly in $\ze$.
This is achieved by the following lemma.

\proclaim{Lemma 10.6}
Let $A,B,C,\delta$ be complex numbers, $C\ne 0,-1,-2,\dots$,
$\{D_N\}_{N=1}^\infty$ be a sequence of complex numbers. Assume that
$$
\lim_{N\to\infty}\frac{D_N}{N}=q\in(0,1).
$$
Then
$$
\lim_{N\to\infty}{}_3F_2\left[\matrix
A,\,B,\,D_N\\N+\delta,\,C\endmatrix\,\biggl|\,1\right]={}_2F_1\left[\matrix
A,\,B\\C\endmatrix\,\biggl|\,q\right].
$$
The convergence is uniform on any set of sequences $\{D_N\}_{N=1}^\infty$ such
that $\{\frac {D_N}N-\lim\frac {D_N}N\}$ uniformly converges to 0 as
$N\to\infty$ and $\lim \frac{D_N}N$ is uniformly bounded from 1 on this set.
\endproclaim

\demo{Proof}
We have
$$
{}_3F_2\left[\matrix A,\,B,\,D_N\\N+\delta,\,C\endmatrix\,\biggl|\,1\right]=
\sum_{k=0}^\infty \frac{(A)_k(B)_k(D_N)_k}{k!(N+\delta)_k(C)_k}\,.
$$
Let us show that these sums converge uniformly in $N>N_0$ for some $N_0>0$.

 Let $d_N$ be the smallest integer greater than $\operatorname{sup}|D_N|$,
where the supremum is taken over all our sequences $\{D_N\}$. Let $l$ be the
largest integer less than $\Re \delta$. Then for large enough $N$ our series is
majorated by the series
$$
\sum_{k=0}^\infty
\left|\frac{(A)_k(B)_k(d_N)_k}{k!(N+l)_k(C)_k}\right|\,.
$$
Using the hypothesis we may assume that for large enough $N$, $d_N<q_0(N+l)$
for some $q_0\in(0,1)$. In particular, $d_N<N+l$. If $k\le N+l-d_N$ then we
have
$$
\gathered
\frac{(d_N)_k}{(N+l)_k}=\frac{d_N(d_N+1)\cdots (d_N+k-1)}{(N+l)(N+l+1)\cdots
(N+l+k-1)}\le \left(\frac{d_N+k-1}{N+l+k-1}\right)^k\\ \le
\left(\frac{d_N+(N+l-d_N)}{N+l+(N+l-d_N)}\right)^k=\left(\frac{N+l}{2N-
d_N+2l}\right)^k\le(2-q_0)^{-k}
\endgathered
$$
for large enough $N$. If $k\ge N+l-d_N$
then we get
$$
\gathered
\frac{(d_N)_k}{(N+l)_k}=\frac{d_N(d_N+1)\cdots (d_N+k-1)}{(N+l)(N+l+1)\cdots
(N+l+k-1)}\\=\frac{d_N(d_N+1)\cdots (N+l-1)}{(d_N+k)(d_N+k+1)\cdots (N+l+k-1)}
\le \left(\frac {N+l}{N+l+k}\right)^{N+l-1-d_N}\\ \le \left( 1+\frac
k{N+l}\right)^{-(N+l-1-d_N)}\le \left( 1+\frac k{N+l} \right)^{-(1-q_0)N}
\endgathered
$$
for large enough $N$. The last
expression is a decreasing function in $N$ (assuming $N+l>0$). Hence, for
$N>N_0$ we get
$$
\frac{(d_N)_k}{(N+l)_k}\le\left( 1+\frac k{N_0+l}
\right)^{-(1-q_0)(N_0+l)}.
$$

Thus, we have proved that $(d_N)_k/(N+l)_k$ does not exceed the maximum of the
$k$th member of a geometric progression with ratio $(2-q_0)^{-1}<1$ and the
inverse of the value at the point $k$ of a polynomial of arbitrarily large
(equal to $[(1-q_0)(N_0+l)]$) degree. Since $(A)_k(B)_k/(k!(C)_k)$ has
polynomial behavior in $k$ for large $k$, this means that we majorated the
series
$$
{}_3F_2\left[\matrix
A,\,B,\,D_N\\N+\delta,\,C\endmatrix\,\biggl|\,1\right]= \sum_{k=0}^\infty
\frac{(A)_k(B)_k(D_N)_k}{k!(N+\delta)_k(C)_k}
$$
by a convergent series with
terms not depending on $N$. Therefore, to compute the limit of this series as
$N\to\infty$, we can take the limit term-wise. Since for any fixed $k\ge 0$
$$
\lim_{N\to\infty}\frac{(D_N)_k}{(N+\delta)_k}=q^k,
$$
we get
$$
\lim_{N\to\infty}\sum_{k=0}^\infty
\frac{(A)_k(B)_k(D_N)_k}{k!(N+\delta)_k(C)_k}=\sum_{k=0}^\infty
\frac{(A)_k(B)_k}{k!(C)_k}\, q^k={}_2F_1\left[\matrix
A,\,B\\C\endmatrix\,\biggl|\,q\right].
$$
As in the proof of Lemma 10.5, the fact that we majorated the series by the
same convergent series for all our sequences $\{D_N\}$, and the uniform
convergence of the terms of the series guarantee the needed uniform
convergence on the set of sequences. \qed
\enddemo

The proof of Proposition 10.4 is complete.\qed
\enddemo

To conclude the proof of Theorem 10.1 we need to get rid of the extra
restrictions $\s\ne 0$ and $z-z'\notin \Z$ imposed in the beginning of the
proof.

Define a function $\Psi:\x\to\C$, which is similar to the function
$\Psi^\n$ introduced in \tht{8.10}, by
$$
\Psi(x)=\cases \psiout(x),&x\in\xout,\\
            \dfrac 1{\psiin(x)}\,,&x\in\xin.
\endcases
$$
Note that for any $x\in\x$, $\Psi(x)$ is an entire function in $(z,z',w,w')$.

\proclaim{Lemma 10.7} Let $(\zw)\in\Dadm'$. The kernel $K(x,y)$ can be
written in the form
$$
K(x,y)=\sqrt{\Psi(x)\Psi(y)}\Kcirc(x,y),
$$
where $\Kcirc(x,y)$ admits a holomorphic continuation in the
parameters to the domain $\Cal D'_0\supset\Dadm'$.
Moreover, for any $(\zw)\in\Cal D'_0$,
$$
\Kcirc(x,y)=\lim_{N\to\infty}N^{1-\Si}\Kcirc^\n(x_N,y_N)
$$
uniformly on compact subsets of $\x\times\x$.
\endproclaim

Recall that the kernel $\Kcirc^\n$ was defined in Lemma 8.8 and $\Cal
D'$ was defined just before this lemma.

\demo{Proof of the lemma} It will be convenient to use a more detailed notation
for the kernels in question. So, we will use the notation
$\Kcirc^\n(x,y\mid\zw)$ instead of $\Kcirc^\n(x,y)$. Next, we define the
kernel $\Kcirc(x,y\mid\zw)$: in the block form,
$$
\gathered
\Kcirc_{\out,\out}(x,y\mid\zw)=
\frac{\Rout(x)\Sout(y)-\Sout(x)\Rout(y)}{x-y}\,,\\
\Kcirc_{\out,\inr}(x,y\mid\zw)=
\Psi(y)\,\frac{\Rout(x)\Rin(y)-\Sout(x)\Sin(y)}{x-y}\,,\\
\Kcirc_{\inr,\out}(x,y\mid\zw)=
\Psi(x)\,\frac{\Rin(x)\Rout(y)-\Sin(x)\Sout(y)}{x-y}\,,\\
\Kcirc_{\inr,\inr}(x,y\mid\zw)=
\Psi(x)\Psi(y)\,\frac{\Rin(x)\Sin(y)-\Sin(x)\Rin(y)}{x-y}\,,
\endgathered
$$
These expressions are well defined if $(\zw)$ is in the subdomain
$$
\Cal D''_0=\{(\zw)\in\Cal D'_0\mid \Si\ne0, \,z-z'\notin\Z\}
$$
By virtue of Propositions 10.2, 10.3, and 10.4,
$$
\lim_{N\to\infty}N^{1-\Si}\Kcirc^\n(x_N,y_N\mid\zw)=
\Kcirc(x,y\mid\zw) \tag10.4
$$
for any fixed $(\zw)\in\Cal D''_0$, uniformly on compact subsets of
$\x\times\x$. Moreover, one can verify that the estimates of
Propositions 10.3 and 10.4 are uniform in $(\zw)$ varying on any
compact subset of the domain $\Cal D''_0$. Thus, the limit relation
\tht{10.4} holds uniformly on compact subsets of $\x\times\x\times\Cal
D''_0$.

On the other hand, we know that the kernel $\Kcirc^\n$ is holomorphic
in $(\zw)$ on the larger domain $\Cal D'_0\supset\Cal D''_0$. It
follows that the additional restrictions $\Si\ne0$ and $z-z'\notin\Z$ can
be removed. Specifically, the right--hand side of \tht{10.4} can be
extended to the domain $\Cal D'_0$ and the limit relation \tht{10.4}
holds on $\x\times\x\times\Cal D'_0$. Indeed, we can avoid the
hyperplanes $\Si=0$ or $z-z'=k$, where $k\in\Z$, by making use of Cauchy's
integral over a small circle in the $z$--plane.

This completes the proof of Lemma 10.7. \qed
\enddemo

Now we can complete the proof of the relation \tht{10.1}. Proposition
10.2 implies that
$$
\lim_{N\to\infty}N^\s\,\Psi^{(N)}(x_N)=\Psi(x) \tag10.5
$$
for any $(\zw)\in\Cal D'_0$, uniformly on compact subsets of $\x$.
Indeed, as is seen from the proof of Proposition 10.2, it does not
use the additional restrictions and holds for any $(\zw)\in\Cal D_0$.
To pass from the functions $\psiout^\n$, $\psiin^\n$, and $\psiout$,
$\psiin$ to the functions $\Psi^\n$ and $\Psi$, we use the assumption
$(\zw)\in\Cal D'_0$ which makes it possible to invert the `inner'
functions for all values of parameters. (Note that if
$(z,z')\in\Zd$ then for $x>\frac 12$ and large enough $N$ both
$\Psi^{(N)}(x_N)$ and $\Psi(x)$ vanish. Similarly, if $(w,w')\in\Zd$
then the vanishing happens for $x<-\frac 12$.)

Since
$$
\gathered
K(x,y)= \Kcirc(x,y)\sqrt{\Psi(x)\Psi(y)},\\
K^{(N)}(x,y)=\Kcirc^\n(x,y)\sqrt{\Psi^{(N)}(x)\Psi^{(N)}(y)},
\endgathered
$$
\tht{10.1} follows from \tht{10.4} and \tht{10.5}. This completes the
proof of Theorem 10.1.
\qed
\enddemo

We conclude this section by a list of properties (without proofs) of the
correlation kernel $K$ and functions $\Rout$, $\Sout$, $\Rin$, $\Sin$. The
proofs can be found in \cite{BD}.

All the results below should be compared with similar results for $K^{(N)}$
and  $\Rout^{(N)}$, $\Sout^{(N)}$, $\Rin^{(N)}$, $\Sin^{(N)}$, which were
proved in the previous sections.

\subhead Symmetries \endsubhead
All four functions  $\Rout,\,\Sout,\,\Rin,\,\Sin$ are invariant with respect
to the transpositions $z\leftrightarrow z'$ and $w\leftrightarrow w'$.

Further, let us denote by $\Cal S$ the following familiar change of the
parameters and the variable: $(z,z',w,w',x)\longleftrightarrow
(w,w',z,z',-x)$. Then
$$
\gather
\Cal S(\psiout)=\psiout,\quad \Cal S(\psiin)=\psiin,\\
\Cal S(\Rout)=\Rout, \quad \Cal S(\Sout)=-\Sout,\quad \Cal S(\Rin)=\Rin, \quad
\Cal S(\Sin)=-\Sin.
\endgather
$$

The functions $\Rout,\,\Sout,\,\Rin,\,\Sin$ and the kernel
$K$ for admissible values of parameters take real values on $\x$.
Moreover, the kernel $K(x,y)$ is J--symmetric, see \S5(f). That is,
$$
\gather
K_{\out,\out}(x,y)=K_{\out,\out}(y,x), \quad
K_{\inr,\inr}(x,y)=K_{\inr,\inr}(y,x),\\
K_{\inr,\out}(x,y)=-K_{\out,\inr}(y,x).
\endgather
$$

\subhead Branching of analytic continuations
\endsubhead
The formulas for $\Rout,\,\Sout,$ $\Rin,$ $\Sin$ above provide analytic
continuations of these functions. We can view $\Rout$ and $\Sout$ as functions
which are analytic and single-valued on $\C\setminus\xin$, and $\Rin$ and
$\Sin$ as functions analytic and single-valued on $\C\setminus\xout$. (Recall
that the Gauss hypergeometric function can be viewed as an analytic and single
valued function on $\C\setminus [1,+\infty)$.)

For a function $F(\ze)$ defined on $\C\setminus \R$ we will denote by $F^+$
and $F^-$ its boundary values:
$$
F^+(x)=F(x+i0), \qquad F^-(x)=F(x-i0).
$$
Then we have
$$
\gather
\text {on  }\xin\qquad
\frac 1{\psiin}\,\frac{\Sout^--\Sout^+}{2\pi i}=\Rin\,,
\quad \frac1{\psiin}\,\frac {\Rout^--\Rout^+}{2\pi i}=\Sin\,,
\\
\text {on  }\xout\qquad
\frac
1{\psiout}\,\frac{\Sin^--\Sin^+}{2\pi i}=\Rout\,,\quad
\frac1{\psiout}\,\frac {\Rin^--\Rin^+}{2\pi i}=\Sout\,.
\endgather
$$

This can also be restated as follows. Let us form a matrix
$$
m=\bmatrix \Rout & -\Sin\\ -\Sout &\Rin\endbmatrix.
$$
Then the matrix $m$
has the jump relation $m_+=m_-v$ on $\x$, where the jump matrix equals
$$
v(x)=\cases \left(\matrix1&2\pi i\,\psiout(x)\\0&1\endmatrix\right),&
x\in\xout,\\
\left(\matrix 1&0\\2\pi i\,\psiin(x)&1\endmatrix\right),&x\in\xin.
\endcases
$$

Furthermore, if $\s>0$ then $m(\ze)\sim I$ as $\ze\to\infty$.

\subhead Differential equations \endsubhead
We use Riemann's notation
$$
P\left(\matrix t_1 & t_2 & t_3\\ a & b & c \\ a' & b' & c'\endmatrix \ \ze
\right) $$
to denote the two--dimensional space of solutions to the second order
Fuchs' equation with singular points $t_1,t_2,t_3$ and exponents
$a,a'$; $b,b'$; $c,c'$, see, e.g., \cite{Er, Vol. 1, 2.6}.

We have
$$
\gathered
\Rout(x)\in P\left(\matrix -\frac12 & \infty & \frac 12\\ w & 0 & z \\ w' &
1-\s & z'\endmatrix\ x \right),\quad \Sout(x)\in
P\left(\matrix -\frac12 & \infty & \frac 12\\ w & 1 & z \\ w' & -\s &
z'\endmatrix\ x\right),\\
\Rin(x)\in P\left(\matrix -\frac12 & \infty & \frac 12\\ -w' & 0 & -z' \\ -w &
1+\s & -z \endmatrix\ x\right),\quad \Sin\in P\left(\matrix -\frac12 &
\infty & \frac 12\\ -w' & 1 & -z' \\ -w & \s & -z\endmatrix\ x\right).
\endgathered
$$

\subhead The resolvent kernel \endsubhead
There exists a limit
$$
L(x,y)=\lim_{N\to\infty}N\cdot L^{(N)}(x_N,y_N)\,,\quad
x,y\in\x.
$$
In the block form corresponding to the splitting
$\x=\xout\sqcup\xin$, the kernel $L(x,y)$ looks as follows:
$$
L=\left[\matrix 0 & \Cal A\\ -{\Cal A}^* & 0\endmatrix\right]\,,
$$
where
$\Cal A$ is a kernel on $\xout\times\xin$ of the form
$$
\Cal A(x,y)=\frac{\sqrt{\psiout(x)\psiin(y)}}
{x-y}\,.
$$

This kernel defines a bounded operator in $L^2(\x,dx)$ if and only if
$|z+z'|<1$ and $|w+w'|<1$. If, in addition, we know that $\s>0$ then we can
prove that $L=K/(1-K)$ or $K=L/(1+L)$ as bounded operators in $L^2(\x,dx)$.

\head 11. Integral parameters $z$ and $w$
\endhead

If one of the parameters $z,z'$ and one of the parameters $w,w'$ are integral
then the measure $P_N$ defined in \S3 is concentrated on a finite set of
signatures, and there is a somewhat simpler way to compute the correlation
kernel of $\Cal P$.

Let us assume that $z=k$ and $w=l$, where $k,l\in\Z$, $k+l\ge 1$. Then
$(z,z',w,w')$ form an admissible quadruple of parameters (see Definition 3.4)
if $z'$ and $w'$ is real and $z'-k>-1$, $w'-l>-1$. We excluded the case $k+l=0$
from our consideration because in this case the measure $P_N$ is concentrated
on one signature.

It is easily seen from the definition of $P_N$ that the measure $P_N$
is now concentrated on the signatures $\la\in\GT_N$ such that
$$
k\ge \la_1\ge\dots\ge\la_N\ge -l.
$$
Note that it may happen that this set does not include the zero signature
because $k$ and $-l$ can be of the same sign.

{}From now on in this section we will assume that $\la$ satisfies the 
inequalities
above. Denote
$$
\x_{k,l}^{(N)}=\left\{-\frac{N-1}2-l,\dots,\frac{N-1}2+k\right\}.
$$
Then
$$
\cl(\la)=\left\{\la_1-1+\frac{N+1}2,\dots,\la_N-
N+\frac{N+1}2\right\}\subset\x_{k,l}^{(N)}.
$$

Let us associate to $\la$ a point configuration $Y(\la)$ in $\x^{(N)}$ as 
follows
$$
Y(\la)=\x_{k,l}^{(N)}\setminus \cl(\la).
$$
Note that $Y(\la)$ defines $\la$ uniquely. Since $|\cl(\la)|=N$, we have 
$|Y(\la)|=k+l$.
Let
$$
Y(\la)=\{y_1,\dots,y_{k+l}\}.
$$
The configuration $Y(\la)$ coincides with the configuration
$X(\la)=\cl(\la)^{\Delta}$ from (4.1) on the set
$\xin^{(N)}\cap\x_{k,l}^{(N)}$.

\proclaim{Proposition 11.1} Let $z=k$, $w=l$ be
integers, $k+l\ge 1$, and $z'>k-1$, $w'>l-1$ be real numbers. Then
$$
P_N(\la)=\const\,\prod_{i=1}^{k+l}\frac{\Gamma(-
y_i+z'+\frac{N+1}2)\Gamma(y_i+w'+\frac {N+1}2)}{\Gamma(-
y_i+k+\frac{N+1}2)\Gamma(y_i+l+\frac {N+1}2)}\,\prod_{1\le i<j\le k+l}(y_i-
y_j)^2.
$$
\endproclaim

\example{Remark 11.2}
Note that for any integer $n$ the shift
$$
k\mapsto k+n,\quad l\mapsto l-n,\quad z'\mapsto z'+n,\quad w'\mapsto w'-n, \quad 
y\mapsto y+n
$$
leaves the measure $P_N$ invariant, cf. Remark 3.7. This means that essentially
the measure depends on three, not four, parameters. If we now set
$l=0$ then $\la$ can be viewed as a Young diagram. Then one can show
that $Y(\la)=\{\frac{N-1}2-\la'_j+j\}_{j=1}^k$\,, where $\la'$ is the
transposed diagram.

\endexample
\demo{Proof of Proposition 11.1}
Set $x_i=\la_i-i+\frac{N+1}2$. Then by Proposition 6.1
$$
P_N(\la)=\const\,\prod\limits_{i=1}^N f(x_i)\prod\limits_{1\le i<j\le N}(x_i-
x_j)^2.
$$
Then, since $\{y_i\}_{i=1}^{k+l}=\x_{k,l}^{(N)}\setminus \{x_i\}_{i=1}^N$,
similarly to Proposition 5.7 we get
$$
P_N(\la)=\const\,\prod\limits_{i=1}^N h(y_i)\prod\limits_{1\le i<j\le k+l}(y_i-
y_j)^2,
$$
where
$$
h(y)=\frac 1{f(y)\prod\limits_{x\in\x_{k,l}^{(N)}\setminus y}(y-x)^2}\,.
$$
Substituting
$$
\prod\limits_{x\in\x_{k,l}^{(N)}\setminus
y}(y-x)^2=\Gamma^2\left(-
y+k+\frac{N+1}2\right)\Gamma^2\left(y+l+\frac{N+1}2\right)
$$
and $f(x)$ from \tht{6.1} we see that
$$
h(y)=\frac{\Gamma(-y+z'+\frac{N+1}2)\Gamma(y+w'+\frac
{N+1}2)}{\Gamma(-y+k+\frac{N+1}2)\Gamma(y+l+\frac {N+1}2)}\,. \qed
$$
\enddemo

Denote by $\Cal P_{k,l}^{(N)}$ the point process consisting of the measure
$P_N(\la)$ on point configurations $Y(\la)$.

Below we will be using Hahn polynomials. These are classical orthogonal
polynomials on a finite set, and we will follow the notation of \cite{NSU}.

\proclaim{Proposition 11.3} For any $n=1,2,\dots$, the $n$th correlation
function of the process $\Cal P_{k,l}^{(N)}$ has the form $$
\rho_n^{(N)}(y_1,\dots,y_n)=\det[K_{k,l}^{(N)}(y_i,y_j)]_{i,j=1}^n.
$$
$K_{k,l}^{(N)}$ is the normalized Christoffel--Darboux kernel for shifted Hahn
polynomials defined as follows:
$$
K_{k,l}^{(N)}(x,y)=\frac
{A_{m-1}}{A_{m}H_{m-1}}\,\frac{\pp_{m}(x)\pp_{m-1}(y)-\pp_{m-1}(x)\pp_{m}(y)}
{x-y}\,\sqrt{h(x)h(y)}\,,
$$
where  $m=k+l$, $h(x)$ is
as above,
$$
\gathered
\pp_{m}(x)=h_m^{(z'-k,w'-l)}\left(x+l+\frac{N-1}2,m+N\right),\\
\pp_{m-1}(x)=h_{m-1}^{(z'-k,w'-l)}\left(x+l+\frac{N-1}2,m+N\right)\\
\endgathered
$$
are Hahn polynomials,
$$
H_{m-1}=\left\Vert h_{m-1}^{(z'-k,w'-l)}\left(x,m+N\right)\right\Vert^2,
$$
and the numbers $A_{m-1}$, $A_m$ are the leading coefficients of
$h_{m-1}^{(z'-k,w'-l)}\left(x,m+N\right)$ and
$h_{m}^{(z'-k,w'-l)}\left(x,m+N\right)$.
\endproclaim

\demo{Proof}
Note that if we shift our phase space $\x_{k,l}^{(N)}$ by $l+\frac{N-1}2$
then the weight function turns into the function
$$
h\left(y-l-\frac{N-1}2\right)=\frac{\Gamma(-y+l+z'+N)\Gamma(y-l+w'+1)}{\Gamma(-
y+k+l+N)\Gamma(y+1)}
$$
on the space
$$
\x_{k,l}^{(N)}+l+\frac{N-1}2=\{0,1,\dots,k+l+N-1\}.
$$

But this is exactly the weight function for the Hahn polynomials
$$
h_n^{(z'-k,w'-l)}(y,k+l+N),\quad n=0,1,2,\dots\,,
$$
see \cite{NSU, 2.4}. Then the claim follows from Proposition 5.1. \qed
\enddemo

Explicit formulas for the Hahn polynomials and their data can be found in
\cite{NSU}.

We know that the processes $\Cal P^{(N)}$ and $\Cal P^{(N)}_{k,l}$ restricted
to the set $\xin^{(N)}\cap\x_{k,l}^{(N)}$ coincide by construction. Same is
true for the correlation kernels, but it is not obvious (recall that the
correlation kernel of a determinantal point process is not defined uniquely,
see \S5(b)).

\proclaim{Proposition 11.4} For any  $x,y\in\xin^{(N)}\cap\x_{k,l}^{(N)}$,
$$
K_{k,l}^{(N)}(x,y)=K_{\inr,\inr}^{(N)}(x,y).
$$
\endproclaim
\demo{Proof} Follows from the relations (here $x\in \x_{k,l}^{(N)}$)
$$
\gathered
\pp_{m}(x)\sqrt{h(x)}=(-1)^{x-k-\frac{N-1}2}\p_{N-1}(x)\sqrt{f(x)},\\
\pp_{m-1}(x)\sqrt{h(x)}=(-1)^{x-k-\frac{N-1}2}\p_{N}(x)\sqrt{f(x)},\\
A_{m-1}=h_{N},\quad A_{m}=h_{N-1},\quad H_{m-1}=h_N.
\endgathered
$$
(The polynomials $\p_{N-1}(x),\p_{N}(x)$ were introduced in \tht{7.1}.)

These relations can be proved either by a direct verification using explicit
formulas (which is rather tedious), or they can be deduced from the following
general fact.

\proclaim{Lemma 11.5 \cite{B5}} Let $\Cal
X=\{x_0,x_1,\dots,x_M\}$ be a finite set of distinct points on the real line,
$u(x)$ and $v(x)$ be two positive functions on $\Cal X$ such that
$$
u(x_k)v(x_k)=\frac 1{\prod_{i\ne k}(x_k-x_i)^2}\,, \quad k=0,1,\dots,M,
$$
and
$P_0,P_1,\dots,P_{M}$ and $Q_0,Q_1,\dots,Q_M$ be the systems of orthogonal
polynomials on $\Cal X$ with respect to the weights $u(x)$ and $v(x)$,
respectively,
$$
\gathered \deg P_i=\deg Q_i=i,\quad \Vert
P_i\Vert^2=p_i,\quad \Vert Q_i\Vert^2=q_i,
\\ P_i=a_ix^i+\text{\rm lower
terms},\quad Q_i=b_ix^i+\text{\rm lower terms}.
\endgathered
$$
Assume that the polynomials are normalized so that $p_i=q_{M-i}$ for all 
$i=0,1,\dots,M$.

Then
$$
\gathered
P_i(x)\sqrt{u(x)}=\epsilon(x)Q_{M-i}(x)\sqrt{v(x)},\quad x\in\Cal X,\\
a_ib_{M-i}=p_i=q_{M-i}, \quad i=0,1,\dots,M,
\endgathered
$$
where
$$
\epsilon(x_k)=\operatorname{sgn}\prod_{i\ne k} (x_k-x_i),\quad k=0,1,\dots,M.
$$
\endproclaim

Taking
$$
M=N+m-1,\quad \Cal X=\x_{k,l}^{(N)},\quad u(x)=f(x),\quad
v(x)=h(x)
$$
we get the needed formulas. The proof of Proposition 11.4 is
complete. \qed
\enddemo

\proclaim{Theorem 11.6}
Assume $z=k$ and $w=l$ are integers, $k+l\ge 1$, $z'$ and $w'$ are real numbers
such that $z'-k>-1$, $w'-l>-1$. The the correlation kernel of the process
$\Cal P$ vanishes if at least one of the arguments is in $\xout$, and on
$\xin\times\xin$ it is equal to the normalized $(k+l)$th Christoffel--Darboux
kernel for the Jacobi polynomials on $\left(-\frac 12,\frac 12\right)$ with the
weight function $(\frac 12-x)^{z'-k}(\frac 12 +x)^{w'-l}$.
\endproclaim

\demo{Proof} One way to prove this statement is to substitute integral $z$ and
$w$ into the formulas of Theorem 10.1. A simpler way, however, is to use the
asymptotic relation
$$
\frac 1{M^n}\,h_{n}^{\al,\be}\left(\left[\frac
{M(1+s)}2\right],M\right)=P_n^{(\al,\be)}(s)+O(M^{-1}), \quad M\to\infty,
$$
where $P_n^{(\al,\be)}$ is the $n$th Jacobi polynomial with parameters
$(\al,\be)$, see, e.g., \cite{NSU, (2.6.3)}.
The estimate is uniform in $s$ belonging to any compact set inside $(-1,1)$.
It is not hard to see that the weight function $h(y)$ as well as the constants
$H_{m-1}$, $A_{m-1}$, $A_m$, see above, converge to the weight function and the
corresponding constants for the Jacobi polynomials.
Then Theorem 9.2  and Proposition 11.3 imply the claim. \qed
\enddemo

\head Appendix
\endhead

The hypergeometric series ${}_3F_2$ evaluated at the unity viewed as a function
of parameters has a large number of two and three--term relations. A lot of
them were discovered by J.~Thomae back in 1879. In 1923,
F.~J.~W.~Whipple introduced a notation which provided a clue to
the numerous formulas obtained by Thomae. An excellent exposition of
Whipple's work was given by W.~N.~Bailey in \cite{Ba, Chapter 3}. We will
be using the notation of \cite{Ba} below.

\subhead Proof of \tht{7.6} \endsubhead
The formula \tht{7.6} coincides with the relation $$Fp(0;4,5)=Fp(0;1,5),$$ see
\cite{Ba, 3.5, 3.6}.

\subhead Proof of Lemma 8.2 \endsubhead The right--hand side of \tht{8.6}
contains two ${}_3F_2$'s. We will use appropriate transformation formulas to
rewrite both of them.

For the first one we employ the relation
$$
\frac{\sin\pi \be_{14}}{\pi \Ga(\al_{014})}\,
Fp(0)=\frac{Fn(1)}{\Ga(\al_{234})\Ga(\al_{245})\Ga(\al_{345})}
-\frac{Fn(4)}{\Ga(\al_{123})\Ga(\al_{125})\Ga(\al_{135})}\,,
\tag A.1
$$
which is \cite{Ba 3.7(1)} with the indices 1 and 2 interchanged, and 3
and 4 interchanged.

By \cite{Ba, 3.5, 3.6}, we have
$$
\gathered
Fp(0)=Fp(0;4,5)=\frac{1}{\Ga(\al_{123})\Ga(\be_{40})\Ga(\be_{50})}\,
{}_3F_2\bmatrix
\al_{145},\,\al_{245},\,\al_{345};\\
\be_{40},\,\be_{50}\endbmatrix\\
=\frac{1}{\Ga(s)\Ga(e)\Ga(f)}\, {}_3F_2\bmatrix
a,\,b,\,c;\\e,\,f\endbmatrix,\\
Fn(1)=Fn(1;2,4)=\frac{1}{\Ga(\al_{124})\Ga(\be_{12})\Ga(\be_{14})}\,
{}_3F_2\bmatrix
\al_{135},\,\al_{013},\,\al_{015};\\
\be_{12},\,\be_{14}\endbmatrix\\
=\frac{1}{\Ga(e-c)\Ga(1+a-b)\Ga(1-b-c+f)}\, {}_3F_2\bmatrix
f-b,\,1-b,\,1-e+a;\\1+a-b,\,1-b-c+f\endbmatrix,\\
Fn(4)=Fp(4;1,2)=\frac{1}{\Ga(\al_{124})\Ga(\be_{41})\Ga(\be_{42})}\,
{}_3F_2\bmatrix
\al_{034},\,\al_{045},\,\al_{345};\\
\be_{41},\,\be_{42}\endbmatrix\\
=\frac{1}{\Ga(e-c)\Ga(1+b+c-f)\Ga(1+a+c-f)}\, {}_3F_2\bmatrix
1-f+c,\,1-s,\,c;\\1+b+c-f,\,1+a+c-f\endbmatrix,
\endgathered
$$
where $s=e+f-a-b-c$.
Thus, \tht{A.1} takes the form
$$
\gathered
{}_3F_2\bmatrix a,\,b,\,c;\\e,\,f\endbmatrix
\\=
\Ga\bmatrix 1-f+a,\,s,\,e,\,f,\,b+c-f\\
e-a,\,b,\,c,\,e-c,\,1+a-b\endbmatrix {}_3F_2\bmatrix
f-b,\,1-b,\,1-e+a;\\1+a-b,\,1-b-c+f\endbmatrix\\
+\Ga\bmatrix 1-f+a,\,e,\,f,\,-b-c+f\\
f-c,\,f-b,\,e-c,\,1+a+c-f\endbmatrix {}_3F_2\bmatrix
1-f+c,\,1-s,\,c;\\1+b+c-f,\,1+a+c-f\endbmatrix,
\endgathered
$$
where we used the identity
$$
\frac \pi{\sin\pi(1-b-c+f)}=\Ga(1-b-c+f)\Ga(b+c-f)=-\Ga(-b-c+f)\Ga(1+b+c-f).
$$
Now set
$$
\gathered
a=N,\, b=-z-w',\, c=-z-w,\, e=u-z+\frac{N+1}2,\, f=-z-z'-w-w',\\
s=e+f-a-b-c=u-z'-\frac{N-1}2.
\endgathered
$$
Also, recall the notation $\s=z+z'+w+w'$. Multiplying the last relation by
$$
\Ga\bmatrix
u+\frac{N+1}2,\,u-z-\frac{N-1}2\\u-\frac{N-1}2,\,u-z+\frac{N+1}2
\endbmatrix
$$
we get
$$
\gathered
\Ga\bmatrix
u+\frac{N+1}2,\,u-z-\frac{N-1}2\\u-\frac{N-1}2,\,u-z+\frac{N+1}2\endbmatrix
{}_3F_2\bmatrix N,\,-z-w',\,-z-w;\\u-z+\frac{N+1}2,\,-\s\endbmatrix\\
=\Ga\bmatrix
u+\frac{N+1}2,\,1+N+\s,\,u-z'-\frac{N-1}2,\,-\s,\,-z+z'\\
u-\frac{N-1}2,\,-z-w',\,-z-w,\,u+w+\frac{N+1}2,\,1+N+z+w'\endbmatrix\\
\times {}_3F_2\bmatrix
-z'-w,\,z+w'+1,\,-u+z+\frac{N+1}2;\\1+N+z+w',\,1+z-z'\endbmatrix\\
+\bigl\{\text{  a similar expression with $z$ and $z'$ interchanged  }\bigr\}.
\endgathered
\tag A.2
$$

This is the transformation for the first term in \tht{8.6}.

As for the second term, we use \cite{Ba, 3.2(2)} which reads
$$
\gathered
{}_3F_2\bmatrix a,\,b,\,c;\\e,\,f\endbmatrix
=\Ga\bmatrix 1-a,\,e,\,f,\,c-b\\ e-b,\, f-b,\,1+b-a,\,c\endbmatrix
{}_3F_2\bmatrix b,\,b-e+1,\,b-f+1;\\1+b-c,\,1+b-a\endbmatrix\\+
\bigl\{\text{  a similar expression with $b$ and $c$ interchanged  }\bigr\}.
\endgathered
$$
Set
$$
a=-N+1,\,b=z+w'+1,\,c=z'+w'+1,\,e=\s+2,\,f=u+w'-\frac{N-3}2.
$$
We get
$$
\gathered
{}_3F_2\bmatrix -N+1,\,z+w'+1,\,z'+w'+1;\\ \s+2,\,u+w'-\frac{N-3}2\endbmatrix
\\
=\Ga\bmatrix N,\,\s+2,\,u+w'-\frac{N-3}2,\,z'-z\\ z'+w+1,\,
u-z-\frac{N-1}2,\,N+1+z+w',\,z'+w'+1\endbmatrix
\\ \times {}_3F_2\bmatrix
z+w'+1,\,-z'-w,\,-u+z+\frac{N+1}2;\\1+z-z',\,1+N+z+w'\endbmatrix\\ +
\bigl\{\text{  a similar expression with $z$ and $z'$ interchanged  }\bigr\}.
\endgathered
$$

Let us multiply this by the prefactor of the hypergeometric function in the
second term of \tht{8.6}. Recalling the formula for $h_{N-1}=h(N-1,z,z',w,w')$,
see \tht{7.2}, and canceling some gamma--factors we obtain that the second
term of \tht{8.6} equals
$$
\gathered
\Ga\bmatrix u+\frac{N+1}2,\,-u+\frac{N+1}2,\, 1+N+\s,\,z+w+1,\,z+w'+1,\,z'-z\\
-u+z'+\frac{N+1}2,\,u+w+\frac{N+1}2,\,
\s+1,\,1+N+z+w'\endbmatrix\\
\times
\frac{\sin\pi(u-z-\frac{N-1}2)}\pi
\,F(u)\,{}_3F_2\bmatrix
z+w'+1,\,-z'-w,\,-u+z+\frac{N+1}2;\\1+z-z',\,1+N+z+w'\endbmatrix\\ +
\bigl\{\text{  a similar expression with $z$ and $z'$ interchanged  }\bigr\},
\endgathered
\tag A.3
$$
where we also used the identity
$$
\Ga(u-z-\tfrac{N-1}2)\Ga(-u+z+\tfrac{N+1}2)=\dfrac\pi{\sin\pi(u-z-\frac{N-
1}2)}\, .
$$
(Recall that $F(u)$ was defined right before Lemma 8.2.)

Now, in order to get \tht{8.6} we have to add \tht{A.2} and \tht{A.3}. Since
both expression have two parts with the second parts different from the first
parts by switching $z$ and $z'$, it suffices to transform the sum of the first
parts. We immediately see that the hypergeometric functions entering the first
parts of \tht{A.2} and \tht{A.3} are identical. By factoring out the
${}_3F_2$'s and some of the gamma--factors, and using the identity
$\Ga(\tau)\Ga(1-\tau)=\pi/{\sin\pi \tau}$ several times, we obtain that the sum
of the first parts of \tht{A.2} and \tht{A.3} equals
$$
\multline
\Ga\bmatrix u+\frac{N+1}2,\,-u+\frac{N+1}2,\, 1+N+\s,\,z+w+1,\,z+w'+1,\,z'-z\\
-u+z'+\frac{N+1}2,\,u+w+\frac{N+1}2,\,
\s+1,\,1+N+z+w'\endbmatrix\\ \times {}_3F_2\bmatrix
z+w'+1,\,-z'-w,\,-u+z+\frac{N+1}2;\\1+z-z',\,1+N+z+w'\endbmatrix
\endmultline
$$
multiplied by
$$
\frac
1\pi\left(-\frac{\sin\pi(u-\frac{N-1}2)\sin\pi(z+w')\sin\pi(z+w)}
{\sin\pi(u-z'-\frac{N-1}2)\sin\pi\s}
+\sin\pi\left(u-z-\tfrac{N-1}2\right)\,F(u)\right).
\tag A.4
$$

Observe that $F(u)$ as a function in $u$ is a linear combination of
$1/\sin\pi(-u+z+\frac{N+1}2)$ and $1/{\sin\pi(-u+z'+\frac{N+1}2)}$.
Thus, \tht{A.4} is a meromorphic function. It is easily verified that
all the singularities of \tht{A.4} are removable and \tht{A.4} is an
entire function. Moreover, since the ratios $\sin(u+\al)/\sin(u+\be)$ are
periodic with period $2\pi$ and are bounded as $\Im u\to\pm\infty$ (for
arbitrary $\al,\be\in \C$), $F(u)$ is bounded on the entire complex plane. By
Liouville's theorem, $F(u)$ does not depend on $u$. Substituting
$u=\tfrac{N-1}2$ we see that $\tht{A.4}$ is equal to $-\sin\pi z/\pi$.

This immediately implies that the sum of \tht{A.2} and \tht{A.3} is equal to
the right--hand side of \tht{8.4}, and the first part of Lemma 8.2 is proved.

Now let us look at the formula \tht{8.7}. Note that the hypergeometric
functions in \tht{8.7} can be obtained from those in \tht{8.6} by the
following shift:
$$
N\mapsto N+1,\quad z\mapsto z-\frac 12,\quad z'\mapsto z'-\frac 12,\quad
w\mapsto w-\frac 12,\quad w'\mapsto w'-\frac 12.
$$
We use for them exactly the same transformation formulas as we used for
\tht{8.6}. By computations very similar to the above, we obtain that the first
term of \tht{8.7} is equal to
$$
\gathered
-\Ga\bmatrix u+\frac{N+1}2,\,-u+\frac{N+1}2,\, N+1,\,\s,\,z'-z\\
-u+z'+\frac{N-1}2,\, u+w+\frac{N+1}2,\,N+1+z+w',\,z'+w,\,z'+w'\endbmatrix \\
\times \frac{\sin\pi(u-\frac{N-1}2)\sin\pi(z+w')\sin\pi(z+w)}
{\pi\sin\pi(u-z'-\frac{N-1}2)\sin\pi\s}\\ \times
{}_3F_2\bmatrix
-z'-w+1,\,z+w',\,-u+z+\frac{N+1}2;\\N+1+z+w',\,1+z-z'\endbmatrix\\
-\bigl\{\text{  a similar expression with $z$ and $z'$ interchanged  }\bigr\},
\endgathered
\tag A.5
$$
while the second term of \tht{8.7} is equal to
$$
\gathered
\Ga\bmatrix u+\frac{N+1}2,\,-u+\frac{N+1}2,\, N+1,\,\s,\,z'-z\\
-u+z'+\frac{N-1}2,\, u+w+\frac{N+1}2,\,N+1+z+w',\,z'+w,\,z'+w'\endbmatrix \\
\times \frac {\sin\pi(u-z-\frac{N-1}2)}\pi\, F(u)\,{}_3F_2\bmatrix
-z'-w+1,\,z+w',\,-u+z+\frac{N+1}2;\\N+1+z+w',\,1+z-z'\endbmatrix\\
+\bigl\{\text{  a similar expression with $z$ and $z'$ interchanged  }\bigr\},
\endgathered
\tag A.6
$$
Adding \tht{A.5} and \tht{A.6} and using the fact that \tht{A.4} is equal to
$-\sin\pi z/\pi$, we arrive at the right--hand side of \tht{8.5}.\qed

\subhead Proof of Lemma 8.4 \endsubhead
We start with deriving a convenient transformation formula for ${}_3F_2$.
\cite{Ba 3.7(6)} with indices 4 and 5 interchanged reads
$$
\gathered
\frac{\sin\pi\be_{40}Fp(0)}{\pi\Ga(\al_{045})\Ga(\al_{034})\Ga(\al_{024})\Ga(\al
_{014})}\\=
-
\frac{Fn(0)}{\Ga(\al_{345})\Ga(\al_{245})\Ga(\al_{145})\Ga(\al_{134})\Ga(\al_{23
4})\Ga(\al_{124})}+
K_0Fn(4),
\endgathered
\tag A.7
$$
where
$$
K_0=\frac {\sin\pi\al_{145}\sin\pi\al_{245}\sin\pi\al_{345}+
\sin\pi\al_{123}\sin\pi\be_{40}\sin\pi\be_{50}}{\pi^3}\,.
$$

By \cite{Ba, 3.6} we have
$$
\gathered
Fp(0)=Fp(0;4,5)=\frac{1}{\Ga(\al_{123})\Ga(\be_{40})\Ga(\be_{50})}\,
{}_3F_2\bmatrix
\al_{145},\,\al_{245},\,\al_{345};\\
\be_{40},\,\be_{50}\endbmatrix\\
=\frac{1}{\Ga(s)\Ga(e)\Ga(f)}\, {}_3F_2\bmatrix
a,\,b,\,c;\\e,\,f\endbmatrix,\\
Fn(0)=Fn(0;4,5)=\frac{1}{\Ga(\al_{045})\Ga(\be_{04})\Ga(\be_{05})}\,
{}_3F_2\bmatrix
\al_{023},\,\al_{013},\,\al_{012};\\
\be_{04},\,\be_{05}\endbmatrix\\
=\frac{1}{\Ga(1-s)\Ga(2-e)\Ga(2-f)}\, {}_3F_2\bmatrix
1-a,\,1-b,\,1-c;\\2-e,\,2-f\endbmatrix,\\
Fn(4)=Fn(4;0,3)=\frac{1}{\Ga(\al_{034})\Ga(\be_{40})\Ga(\be_{43})}\,
{}_3F_2\bmatrix
\al_{124},\,\al_{145},\,\al_{245};\\
\be_{40},\,\be_{43}\endbmatrix\\
=\frac{1}{\Ga(1-f+c)\Ga(e)\Ga(1+a+b-f)}\, {}_3F_2\bmatrix
e-c,\,a,\,b;\\e,\,1+a+b-f\endbmatrix,\\
K_0=\frac {\sin\pi a\sin\pi b\sin\pi c +\sin\pi s\sin\pi e\sin\pi
f}{\pi^3}\,,
\endgathered
$$
where $s=e+f-a-b-c$. Then \tht{A.7} takes the form
$$
\gathered
{}_3F_2\bmatrix
a,\,b,\,c;\\e,\,f\endbmatrix\\
=-\frac{\pi}{\sin\pi e}\,\Ga\bmatrix
1-f+c,\,1-f+b,\,1-f+a,\,s,\,e,\,f\\
c,\,b,\,a,\,e-b,\,e-a,\,e-c,\,2-e,\,2-f
\endbmatrix{}_3F_2\bmatrix
1-a,\,1-b,\,1-c;\\2-e,\,2-f\endbmatrix
\\
+\frac{\pi K_0}{\sin\pi e}\,\Ga\bmatrix
1-s,\,1-f+b,\,1-f+a,\,s,\,f\\
1+a+b-f\endbmatrix{}_3F_2\bmatrix
e-c,\,a,\,b;\\e,\,1+a+b-f\endbmatrix,
\endgathered
\tag A.8
$$
with $K_0$ as above.

Let us apply \tht{A.8} to both hypergeometric functions in \tht{8.6}. For the
first one we set
$$
\gathered
a=N,\quad b=-z-w',\quad c=-z-w,\quad e=-\s,\quad f=u-z+\frac{N+1}2\,,\\
s=e+f-a-b-c=u-z'-\frac{N-1}2.
\endgathered
$$
Then
$K_0=-\sin\pi\s\sin\pi(u-z'-\frac{N-1}2)\sin\pi(u-
z+\frac{N+1}2)/\pi^3$,\footnote{Here we
used that $N\in\Z$ and $\sin\pi a=\sin\pi N=0$.}
and rewriting the sines in $K_0$ as products of gamma-functions, we have $$
\gathered
{}_3F_2\bmatrix N,\,-z-w',\,-z-w;\\u-z+\frac{N+1}2,\,-\s\endbmatrix
\\= \Ga\bmatrix
-u-w-\frac{N-1}2,\,-u-w'-\frac{N-1}2,\,-u+z+\frac{N+1}2,\,u-z'-\frac{N-1}2,\,-
\s,\,u-z+\frac{N+1}2\\
-z-w,\,-z-w',\,N,\,-z'-w,\,-\s-N,\,-z'-w',\,2+\s,\,-u+z-\frac{N-3}2
\endbmatrix\\ \times\frac{\pi}{\sin\pi\s}\,{}_3F_2\bmatrix
-N+1,\,z+w'+1,\,z+w+1;\\ \s+2,\,-u+z-\frac{N-3}2\endbmatrix
\\
+\Ga\bmatrix
-u-w'-\frac{N-1}2,\,-u+z+\frac{N+1}2\\
-u-w'+\frac{N+1}2,\,-u+z-\frac{N-1}2\endbmatrix{}_3F_2\bmatrix
-z'-w',\,N,\,-z-w';\\-\s,\,-u-w'+\frac{N+1}2\endbmatrix.
\endgathered
\tag A.9
$$

For the second one we set
$$
\gathered
a=-N+1,\quad b=z+w'+1,\quad c=z'+w'+1,\quad e=\s+2,\quad
f=u+w'-\frac{N-3}2,\\s=e+f-a-b-c=u+w+\frac{N+1}2\,.
\endgathered
$$
Then $K_0=\sin\pi\s\sin\pi(u+w+\frac{N+1}2)\sin\pi(u+w'-\frac{N-3}2)/\pi^3$.
Observe that the first term in the right--hand side of \tht{A.8} will now
vanish thanks to the $\Ga(a)$ in the denominator (remember that
$N\in\{1,2,\dots\}$). Rewriting sines as products of gamma--functions again, we
get
$$
\gathered
{}_3F_2\bmatrix
-N+1,\,z+w'+1,\,z'+w'+1;\\\s+2,\,u+w'-\frac{N-3}2\endbmatrix\\
=
\Ga\bmatrix
-u+z+\frac{N+1}2,\,-u-w'-\frac{N-1}2\\-u-w'+\frac{N-1}2,\,
-u+z-\frac{N-3}2\\ \endbmatrix{}_3F_2\bmatrix
z+w+1,\,-N+1,\,z+w'+1;\\\s+2,\,-u+z-\frac{N-3}2 \endbmatrix.
\endgathered
\tag A.10
$$

Now let us substitute \tht{A.9} and \tht{A.10} into \tht{8.6}. Observe that
the two ${}_3F_2$'s from the right--hand sides of \tht{A.9} and \tht{A.10} are
obtained from the ${}_3F_2$'s from \tht{8.6} by the change
$$
u\mapsto -u, \quad (z,z',w,w')\mapsto (w',w,z',z).
\tag A.11
$$

Our goal is to show that the prefactors of these ${}_3F_2$'s after substitution
will be symmetric to the prefactors of ${}_3F_2$'s in \tht{8.6} with respect
to  \tht{A.11}.

One prefactor is easy to handle. The coefficient of ${}_3F_2\bmatrix
-z'-w',\,N,\,-z-w';\\-\s,\,-u-w'+\frac{N+1}2\endbmatrix$ after the
substitution of \tht{A.9} and \tht{A.10} into \tht{8.6} equals
$$
\gathered
\Ga\bmatrix u+\frac{N+1}2,\,u-z-\frac{N-1}2\\ u-\frac{N-1}2,\,
u-z+\frac{N+1}2\endbmatrix\,\Ga\bmatrix
-u-w'-\frac{N-1}2,\,-u+z+\frac{N+1}2\\
-u-w'+\frac{N+1}2,\,-u+z-\frac{N-1}2\endbmatrix\\=
\Ga\bmatrix -u+\frac{N+1}2,\,-u+w'-\frac{N-1}2\\ -u-\frac{N-1}2,\,
-u+w'+\frac{N+1}2\endbmatrix,
\endgathered
$$
which is symmetric to $\Ga\bmatrix u+\frac{N+1}2,\,u-z-\frac{N-1}2\\
u-\frac{N-1}2,\, u-z+\frac{N+1}2\endbmatrix$ with respect to \tht{A.11}.

As for the other prefactor, the verification is more involved. Namely, we need
to prove the following equality:
$$
\gathered
\Ga\bmatrix
-u-w-\frac{N-1}2,\,-u-w'-\frac{N-1}2,\,-u+z+\frac{N+1}2,\,u-z'-\frac{N-1}2,\,-
\s,\,u-z+\frac{N+1}2\\
-z-w,\,-z-w',\,N,\,-z'-w,\,-\s-N,\,-z'-w',\,2+\s,\,-u+z-\frac{N-3}2
\endbmatrix
\\ \times
\Ga\bmatrix u+\frac{N+1}2,\,u-z-\frac{N-1}2\\ u-\frac{N-1}2,\,
u-z+\frac{N+1}2\endbmatrix\,\frac{\pi}{\sin\pi\s}
\\+
\Gamma\bmatrix u+\frac{N+1}2,\,-u+\frac{N+1}2\\
-u+z+\frac{N+1}2,\,-u+z'+\frac{N+1}2,\,u+w+\frac{N+1}2,\,
u+w'-\frac{N-3}2\endbmatrix
\\ \times
\Ga\bmatrix
-u+z+\frac{N+1}2,\,-u-w'-\frac{N-1}2\\-u-w'+\frac{N-1}2,\,
-u+z-\frac{N-3}2\endbmatrix
\frac{1}{h(N-1,z,z',w,w')}\,\frac{1}{\sin\pi\s\,\sin\pi(z'-z)}
\\ \times
\left(\dfrac{\sin\pi(z+w)\sin\pi(z+w')\sin\pi z'}
{\sin\pi(-u+z'+\frac{N+1}2)}-\dfrac{\sin\pi(z'+w)\sin\pi(z'+w')\sin\pi
z}{\sin\pi(-u+z+\frac{N+1}2)}\right)
\\=
\Gamma\bmatrix -u+\frac{N+1}2,\,u+\frac{N+1}2\\
u+w'+\frac{N+1}2,\,u+w+\frac{N+1}2,\,-u+z'+\frac{N+1}2,\,
-u+z-\frac{N-3}2\endbmatrix
\\ \times
\frac{1}{h(N-1,w',w,z',z)}\,\frac{1}{\sin\pi\s\,\sin\pi(w-w')}
\\ \times
\left(\dfrac{\sin\pi(z'+w')\sin\pi(z+w')\sin\pi w}
{\sin\pi(u+w+\frac{N+1}2)}-\dfrac{\sin\pi(z'+w)\sin\pi(z+w)\sin\pi
w'}{\sin\pi(u+w'+\frac{N+1}2)}\right).
\endgathered
$$
After massive cancellations\footnote{which also rely on the fact that
$N\in\Z$, cf. the previous footnote.} this equality is reduced to the following
trigonometric identity (here $y=u-\frac{N-1}2$): $$
\gathered
\frac{\sin\pi
y\sin\pi(z+w)\sin\pi(z'+w)\sin\pi(z+w')\sin\pi(z'+w')}
{\sin\pi(y-z)\sin\pi(y-z')\sin\pi(y+w)\sin\pi(y+w')}\\
=\frac{\sin\pi(z'+w)\sin\pi(z'+w')\sin\pi z}{\sin\pi(z-z')\sin\pi(y-z)}
+\frac{\sin\pi(z+w)\sin\pi(z+w')\sin\pi z'}{\sin\pi(z'-z)\sin\pi(y-z')}\\
+\frac{\sin\pi(z+w')\sin\pi(z'+w')\sin\pi w}{\sin\pi(w-w')\sin\pi(y+w)}
+\frac{\sin\pi(z+w)\sin\pi(z'+w)\sin\pi w'}{\sin\pi(w'-w)\sin\pi(y+w')}\,.
\endgathered
$$

One way to prove this identity is to view both sides as meromorphic functions
in $y$. Then it is easily verified that the difference of the left-hand side
and the right--hand side is an entire function. Moreover, both sides are
periodic with period $2\pi$ and bounded for $|\Im y|$ large enough. This
implies that both sides are identically equal. The proof of Lemma 8.4 is
complete.

\subhead On the proof of Lemma 8.5 \endsubhead
This proof is very similar to that of Lemma 8.4 above. The needed
transformation formulas for the ${}_3F_2$'s are obtained from \tht{A.9} and
\tht{A.10} by the shift $(N,z,z',w,w')\mapsto (N+1,z-\frac12,z'-\frac
12,w-\frac 12,w'-\frac 12)$. After substituting the resulting expressions into
\tht{8.7} we collect the coefficients of ${}_3F_2$'s and compare them with
what we want. As in the proof of Lemma 8.4, one of the desired equalities
follows immediately, while the other is reduced to the same trigonometric
identity.

\subhead On analytic continuation of the series ${}_3F_2(1)$
\endsubhead
Here we prove that the function
$$
\frac{1}{\Ga(e)\Ga(f)\Ga(e+f-a-b-c)}\,{}_3F_2\left[\matrix a,\,b,\,c\,\\
e,\,f\endmatrix\,\Bigr|\,1\right]
$$
can be analytically continued to an entire function in 5 complex variables
$a,b,c,e,f$. We stated this claim in the beginning of \S7 and used it
in \S8.

Apply the transformation formula
$$
{}_3F_2\left[\matrix a,\,b,\,c\,\\
e,\,f\endmatrix\,\Bigr|\,1\right]=
\Ga\bmatrix e,\,f,\,s \\ a,\,s+b,\,s+c\endbmatrix
{}_3F_2\left[\matrix e-a,\,f-a,\,s\,\\
s+b,\,s+c\endmatrix\,\Bigr|\,1\right],
$$
where $s=e+f-a-b-c$. It allows us to conclude that the function in
question continues to the domain $\Re(a)>0$ (other parameters being
arbitrary).

Likewise, we can continue to the domain $\Re(b)>0$ and also to the
domain $\Re(c)>0$. Then one can apply a general theorem about
`forced' analytic continuation of holomorphic 
functions on tube domains: see, e.g., \cite{H, Theorem 2.5.10}.

\vskip 2 true cm

{\smc A.~Borodin}: School of Mathematics, Institute for Advanced
Study, Einstein Drive, Princeton NJ 08540, U.S.A.

E-mail address:
{\tt borodine\@math.upenn.edu}

{\smc G.~Olshanski}: Dobrushin Mathematics Laboratory, Institute for
Information Transmission Problemes, Bolshoy Karetny 19, 101447
Moscow GSP-4, RUSSIA.

E-mail address: {\tt olsh\@iitp.ru, olsh\@online.ru}

\newpage

$$\epsffile{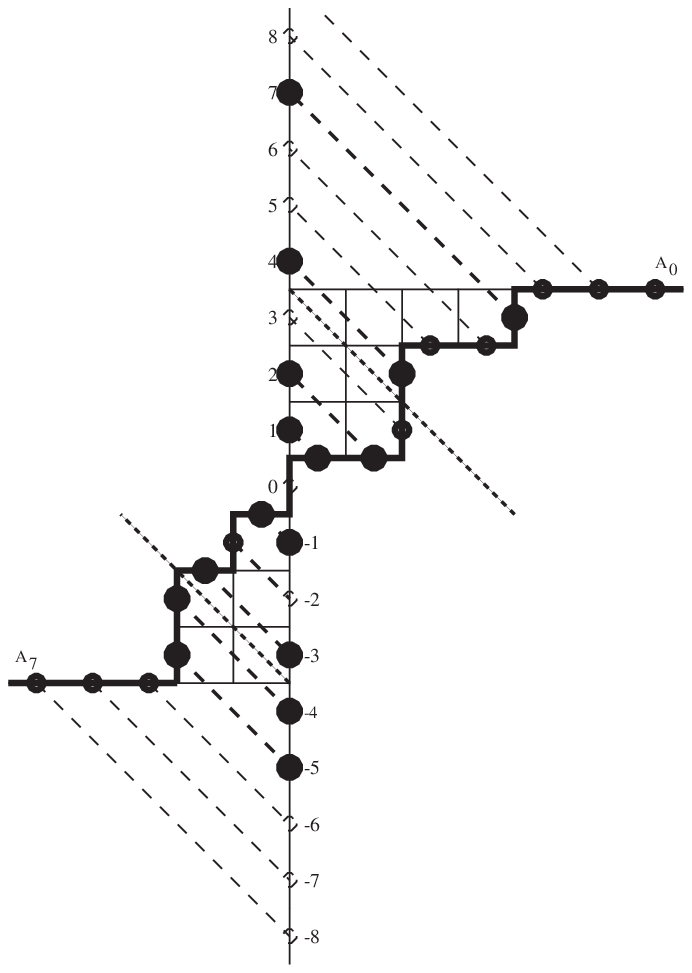}$$
\centerline{Figure (Proposition 4.1)}


\Refs

\widestnumber\key{NMPZ}

\ref\key AF
\by A.~Alastuey and P.~J.~Forrester
\paper Correlations in two--component log--gas systems
\jour J. Statist. Phys.
\vol 81
\yr 1995
\pages 579--627
\endref

\ref\key As
\by R.~Askey
\paper An integral of Ramanujan and orthogonal polynomials
\jour J. Indian Math. Soc.
\vol 51
\yr 1987
\pages 27--36
\endref

\ref\key Ba
\by W.~N.~Bailey
\book Generalized Hypergeometric Series
\publ Cambridge
\yr 1935
\endref

\ref\key BC
\by R.~Beals and R.~R.~Coifman
\paper Scattering and inverse scattering for first order systems
\jour Comm. Pure Appl. Math.
\vol 37
\yr 1984
\pages 39--90
\endref

\ref\key BDT
\by R.~Beals, P.~Deift, C.~Tomei
\book Direct and inverse scattering on the line
\bookinfo Mathematical surveys and monographs
\publ Amer. Math. Soc.
\vol 28
\yr 1988
\endref


\ref\key B1
\by A.~Borodin
\paper Characters of symmetric groups and correlation functions of
point processes
\jour  Funct. Anal. Appl.
\vol 34
\yr 2000
\issue 1
\pages 10--23
\endref

\ref\key B2
\bysame
\paper Harmonic analysis on the infinite symmetric group and the
Whittaker kernel
\jour St.~Petersburg Math. J.
\vol 12
\yr 2001
\issue 5
\endref

\ref\key B3
\bysame
\paper Riemann--Hilbert problem and the discrete Bessel kernel
\jour Intern. Math. Research Notices
\yr 2000
\issue 9
\pages 467--494; {\tt math/9912093}
\endref

\ref\key B4
\bysame
\paper Biorthogonal ensembles
\jour Nucl. Phys. B
\vol 536
\yr 1998
\pages 704--732;
{\tt math/9804027}
\endref

\ref\key B5
\bysame
\paper Duality of orthogonal polynomials on a finite set
\pages {\tt math/0101125}
\endref

\ref\key BD
\by A.~Borodin and P.~Deift
\paper Fredholm determinants, Jimbo-Miwa-Ueno tau-functions, and
representation theory
\pages {\tt math/0111007}
\endref

\ref\key BOO
\by A.~Borodin, A.~Okounkov and G.~Olshanski
\paper Asymptotics of Plancherel measures for symmetric groups
\jour J. Amer. Math. Soc.
\vol 13
\yr 2000
\pages 491--515; {\tt math/9905032}
\endref

\ref\key BO1
\by A.~Borodin and G.~Olshanski
\paper Point processes and the infinite symmetric group
\jour Math. Research Lett.
\vol 5
\yr 1998
\pages 799--816; {\tt math/9810015}
\endref

\ref\key BO2
\bysame
\paper Distributions on partitions, point processes and the hypergeometric
kernel
\jour Comm. Math. Phys.
\vol 211
\yr 2000
\issue 2
\pages 335--358; {\tt math/9904010}
\endref

\ref\key BO3
\bysame
\paper Z--Measures on partitions, Robinson--Schensted--Knuth
correspondence, and
$\beta=2$ random matrix ensembles
\jour  Mathematical Sciences Research Institute Publications
\vol 40
\yr 2001
\pages 71--94; {\tt math/9905189}
\endref

\ref\key BO4
\bysame
\paper Infinite random matrices and ergodic measures
\jour to appear in Comm. Math. Phys.
\pages{\tt math-ph/0010015}
\endref

\ref
\key Boy
\by R.~P.~Boyer
\paper Infinite traces of AF-algebras and characters of $U(\infty)$
\jour J.\ Operator Theory
\vol 9
\yr 1983
\pages 205--236
\endref

\ref\key CJ1
\by F.~Cornu, B.~Jancovici
\paper On the two-dimensional Coulomb gas
\jour  J. Statist. Phys.
\vol 49
\yr 1987
\issue 1-2
\pages 33--56
\endref

\ref\key CJ2
\bysame
\paper The electrical double layer: a solvable model
\jour Jour. Chem. Phys.
\vol 90
\yr 1989
\pages 2444
\endref


\ref\key DVJ
\by D.~J.~Daley, D.~Vere--Jones
\book An introduction to the theory of point processes
\bookinfo Springer series in statistics
\publ Springer
\yr 1988
\endref

\ref\key De
\by P.~Deift
\paper Integrable operators
\inbook Differential operators and spectral theory: M. Sh. Birman's
70th anniversary collection (V.~Buslaev, M.~Solomyak, D.~Yafaev,
eds.)
\bookinfo
American Mathematical Society Translations, ser. 2, v. 189
\publ Providence, R.I.: AMS
\yr 1999
\endref

\ref\key Dy
\by F.~J.~Dyson
\paper Statistical theory of the energy levels of complex systems I,
II, III
\jour J. Math. Phys.
\vol 3
\yr 1962
\pages 140--156, 157--165, 166--175
\endref

\ref \key Ed
\by A.~Edrei
\paper On the generating function of a doubly--infinite,
totally positive sequence
\jour Trans.\ Amer.\ Math.\ Soc.\
\vol 74 \issue 3 \pages 367--383 \yr 1953
\endref

\ref\key Er
\by A.~Erdelyi (ed.)
\book Higher transcendental functions, {\rm Vols. 1, 2}
\publ Mc Graw--Hill
\yr 1953
\endref


\ref\key F1
\by P.~J.~Forrester
\paper Positive and negative charged rods alternating along a line:
exact results
\jour J. Statist. Phys.
\vol 45 \yr 1986 \issue 1--2
\pages 153--169
\endref

\ref\key F2
\bysame
\paper Solvable isotherms for a two--component system of charged rods
on a line
\jour J. Statist. Phys.
\vol 51 \yr 1988 \issue 3--4
\pages 457--479
\endref

\ref\key F3
\bysame
\paper Exact results for correlations in a two--component log--gas
\jour J. Statist. Phys.
\vol 54 \yr 1989 \issue 1--2
\pages 57-79
\endref

\ref\key G
\by M.~Gaudin
\paper L'isotherme critique d'un plasma sur r\'eseau $(\beta=2,\;d=2,\;n=2)$
\jour  J. Physique
\vol 46
\issue 7
\yr 1985
\pages 1027--1042
\endref

\ref\key H
\by L.~H\"ormander
\book An introduction to complex analysis in several variables
\publ D.~van Nostrand
\publaddr Princeton, NJ
\yr 1966
\endref

\ref\key IIKS
\by A.~R.~Its, A.~G.~Izergin, V.~E.~Korepin, N.~A.~Slavnov
\paper Differential equations for quantum correlation functions
\jour Intern. J. Mod. Phys.
\vol B4
\yr 1990
\pages 10037--1037
\endref

\ref\key J
\by K.~Johansson
\paper Discrete orthogonal polynomial ensembles and the Plancherel
measure
\jour Ann. Math. (2)
\vol 153
\yr 2001
\issue 1
\pages 259--296; {\tt math/9906120}
\endref

\ref\key Ke
\by S.~V.~Kerov
\paper Elementary introduction to the $z$--measures
\paperinfo unpublished manuscript
\endref

\ref\key KOO
\by S.~Kerov, A.~Okounkov, G.~Olshanski
\paper The boundary of Young graph with Jack edge multiplicities
\jour Intern. Math. Res. Notices
\yr 1998
\issue 4
\pages 173--199
\endref

\ref \key KOV
\by S.~Kerov, G.~Olshanski, A.~Vershik
\paper Harmonic analysis on the infinite symmetric group. A deformation
of the regular representation
\jour Comptes Rend. Acad. Sci. Paris, S\'er. I
\vol 316
\yr 1993
\pages 773--778; detailed version in preparation
\endref

\ref \key KBI
\by  V.~E.~Korepin, N.~M.~Bogoliubov, A.~G.~Izergin
\book Quantum inverse scattering method and correlation functions
\publ Cambridge University Press
\yr 1993
\endref


\ref\key Len
\by A.~Lenard
\paper Correlation functions and the uniqueness of the state in classical
statistical mechanics
\jour Comm. Math. Phys
\vol 30
\yr 1973
\pages 35--44
\endref

\ref\key Les1
\by P.~A.~Lesky
\paper Unendliche und endliche Orthogonalsysteme von Continuous
Hahnpolynomen
\jour Results in Math.
\vol 31
\yr 1997
\pages 127--135
\endref

\ref\key Les2
\bysame
\paper Eine Charakterisierung der kontinuierlichen und diskreten
klassischen Orthogonalpolynome
\paperinfo Preprint 98--12, Mathematisches Institut A, Universitaet
Stuttgart (1998)
\endref

\ref \key Me
\by M.~L.~Mehta
\book Random matrices
\publ 2nd edition, Academic Press, New York
\yr 1991
\endref

\ref\key MT
\by L.~M.~Milne--Thomson
\book The calculus of finite differences
\publ London, Macmillan \& Co
\yr 1933
\endref

\ref\key NW
\by T.~Nagao, M.~Wadati
\paper Correlation functions of random matrix ensembles related to
classical orthogonal polynomials
\jour  J. Phys. Soc. Japan
\vol 60
\issue 10
\yr 1991
\pages 3298-3322
\endref

\ref \key Ner
\by Yu.~A.~Neretin
\paper Hua type integrals over unitary groups and over projective
limits of unitary groups
\jour Duke Math. J., to appear; {\tt math-ph/0010014}
\endref

\ref\key NSU
\by A.~F.~Nikiforov, S.~K.~Suslov and V.~B.~Uvarov
\book Classical orthogonal polynomials of a discrete variable
\bookinfo Springer Series in Computational Physics
\publ Springer
\publaddr New York
\yr1991
\endref

\ref\key NMPZ
\by S.~Novikov, S.~V.~Manakov, L.~P.~Pitaevskii, V.~E.~Zakharov
\book
Theory of Solitons: The Inverse Scattering Method
\bookinfo  Contemporary Soviet Mathematics
\publ Consultants Bureau [Plenum]
\publaddr New York--London
\yr 1984
\endref

\ref \key OkOl
\by A.~Okounkov and G.~Olshanski
\paper Asymptotics of Jack polynomials as the number of variables
goes to infinity
\jour Intern. Math. Res. Notices
\yr 1998
\issue 13
\pages 641--682
\endref

\ref \key Ol1
\by G.\ Olshanski
\paper Unitary representations of infinite-dimensional
pairs $(G,K)$ and the formalism of R.\ Howe
\jour Soviet Math. Doklady
\vol 27
\issue 2
\yr 1983
\pages 290--294
\endref

\ref \key Ol2
\bysame
\paper Unitary representations of infinite-dimensional
pairs $(G,K)$ and the formalism of R.\ Howe
\inbook Representation of Lie Groups and Related Topics
\eds A.\ Vershik and D.\ Zhelobenko
\bookinfo Advanced Studies in Contemporary Math. {\bf 7}
\publ Gordon and Breach Science Publishers
\publaddr New York etc.
\yr 1990
\pages 269--463
\endref

\ref\key Ol3
\bysame
\paper The problem of harmonic analysis on the
infinite--dimensional unitary group, {\tt math/0109193}
\endref

\ref \key OV
\by G.\ Olshanski and A.\ Vershik
\paper Ergodic unitary invariant measures on the space
of infinite Hermitian matrices
\inbook Contemporary Mathematical Physics
\eds R.~L.~Dobrushin, R.~A.~Minlos, M.~A.~Shubin,
A.~M.~Vershik
\bookinfo American Mathematical Society Translations, Ser.~2, Vol.~175
\publ  Amer.\ Math.\ Soc.
\publaddr Providence
\yr 1996
\pages 137--175
\endref

\ref\key  P.I
\by G.~Olshanski
\paper Point processes and the infinite symmetric group. Part I: The
general formalism and the density function
\paperinfo Preprint, 1998, {\tt math/9804086}
\endref

\ref\key P.II
\by A.~Borodin
\paper Point processes and the infinite symmetric group. Part II:
Higher correlation functions
\paperinfo Preprint,
 1998, {\tt math/9804087}
\endref

\ref\key P.III
\by A.~Borodin and G.~Olshanski
\paper Point processes and the infinite symmetric group. Part III:
Fermion point processes
\paperinfo Preprint, 1998, {\tt math/9804088}
\endref

\ref\key P.IV
\by A.~Borodin
\paper Point processes and the infinite symmetric group. Part IV:
Matrix Whittaker kernel
\paperinfo Preprint, 1998, {\tt math/9810013}
\endref

\ref\key P.V
\by G.~Olshanski
\paper Point processes and the infinite symmetric group. Part V:
Analysis of the matrix Whittaker kernel
\paperinfo Preprint, 1998, {\tt math/9810014}
\endref

\ref \key Pi1
\by D.~Pickrell
\paper Measures on infinite dimensional Grassmann manifold
\jour J.~Func.\ Anal.\
\vol 70
\yr 1987 
\pages 323--356
\endref


\ref\key Pi2
\bysame
\paper Mackey analysis of infinite classical motion groups
\jour Pacific J. Math.
\vol 150
\yr 1991
\pages 139--166
\endref

\ref
\key Sh
\by A.~Shiryaev
\book Probability
\publ Springer-Verlag
\publaddr New York
\yr 1996
\endref

\ref\key So
\by A.~Soshnikov
\paper Determinantal random point fields
\jour Russian Math. Surveys
\vol 55
\yr 2000
\issue 5
\pages 923--975; {\tt math/0002099}
\endref

\ref\key Th
\by E.~Thoma
\paper Die unzerlegbaren, positive-definiten Klassenfunktionen
der abz\"ahlbar unendlichen, symmetrischen Gruppe
\jour Math.~Zeitschr.
\vol 85
\yr 1964
\pages 40-61
\endref

\ref\key TW
\by C.~A.~Tracy and H.~Widom
\paper Universality of the distribution functions of random matrix theory
\inbook Integrable systems: from classical to quantum (Montr\'eal,
QC, 1999)
\bookinfo CRM Proc. Lecture Notes
\vol 26
\publ Amer. Math. Soc.
\yr 2000
\pages 251--264; {\tt solv-int/9901003}
\endref

\ref\key VK1
\by A.~M.~Vershik, S.~V.~Kerov
\paper Asymptotic theory of characters of the symmetric group
\jour Funct. Anal. Appl.
\vol 15
\yr 1981
\pages 246--255
\endref

\ref \key VK2
\bysame
\paper Characters and factor representations of the
infinite unitary group
\jour Soviet Math.\ Doklady
\vol 26
\pages 570--574
\yr 1982
\endref

\ref \key Vo
\by D.~Voiculescu
\paper Repr\'esentations factorielles de type {\rm II}${}_1$ de
$U(\infty)$
\jour J.\ Math.\ Pures et Appl.\
\vol 55 \pages 1--20 \yr 1976
\endref

\ref\key Wa
\by A.~J.~Wassermann
\paper Automorphic actions of compact groups on operator algebras
\paperinfo Thesis, University of Pennsylvania
\yr 1981
\endref

\ref \key Zh
\by D.~P.~Zhelobenko
\book Compact Lie groups and their representations
\publ Nauka, Moscow, 1970 (Russian); English translation: Transl.
Math. Monographs {\bf 40}, Amer. Math. Soc., Providence, R.I., 1973
\endref

\endRefs
\enddocument